\documentclass[10pt]{amsart}
\usepackage[letterpaper,textheight=249mm,vcentering,left=2.5cm,right=2.5cm,headsep=10pt]{geometry}
\usepackage{amsmath,amssymb,amsthm,mathtools,mathrsfs}
\usepackage{booktabs,array,enumitem}
\usepackage[hidelinks]{hyperref}
\hypersetup{pdftitle={Families of regular non-Berwald Landsberg surfaces},
 pdfauthor={Ming Li and Kefeng Liu},
 pdfsubject={Version 4.4; seed families, their geometry, and local exact return}}
\allowdisplaybreaks[2]
\hypersetup{bookmarksdepth=1}

\makeatletter
\renewcommand{\@pnumwidth}{3em}
\makeatother
\newtheorem{theorem}{Theorem}[section]
\newtheorem{lemma}[theorem]{Lemma}
\newtheorem{proposition}[theorem]{Proposition}
\newtheorem{corollary}[theorem]{Corollary}
\theoremstyle{definition}
\newtheorem{definition}[theorem]{Definition}
\newtheorem{example}[theorem]{Example}
\theoremstyle{remark}\newtheorem{remark}[theorem]{Remark}
\newcommand{\dd}{\mathrm d}
\newcommand{\Id}{\operatorname{Id}}
\newcommand{\Div}{\operatorname{div}}
\newcommand{\Curl}{\operatorname{curl}}
\newcommand{\supp}{\operatorname{supp}}
\newcommand{\Op}{\operatorname{Op}}

\newcommand{\R}{\mathbb R}
\newcommand{\C}{\mathbb C}

\numberwithin{equation}{section}
\title[Regular non-Berwald Landsberg surfaces]{Families of regular non-Berwald Landsberg surfaces}
\author{Ming Li}
\author{Kefeng Liu}
\address{Mathematical Science Research Center, Chongqing University of Technology}
\email[Ming Li]{mingli@cqut.edu.cn}
\email[Kefeng Liu]{kefengliu@cqut.edu.cn}
\date{October 12, 2026; revised November 23, 2026}
\begin{document}
\raggedbottom
\begin{abstract}
We construct families of smooth, strongly convex, non-Berwald Landsberg metrics
of zero flag curvature on open subsets of the plane. Each metric is regular
in every nonzero tangent direction and is nonreversible. The indicatrices are
obtained from the inverse differential of a periodically forced oscillator.
The construction is indexed by smooth balanced periodic seeds satisfying
a simple characteristic condition at a nonzero initial value. It includes
single-frequency seeds of every even frequency at least four, mixtures
of interacting frequencies, and seed neighborhoods of arbitrarily large
finite dimension within the balanced class. A normalized
profile of the canonical Cartan invariant distinguishes fixed seed orbits
at sufficiently small amplitudes. In particular, the construction yields
countably infinitely many pairwise locally nonhomothetic metrics.

The analytic problem is to make the actual return map equal to the identity
on an open set of nonzero initial values. We show that exact return on
a fixed open set forces three balance equations on the seed. Under a
positive drift condition a finite formal approximation gives
an arbitrarily high prescribed power of a small amplitude in the return error.
We normalize the complete two-branch linearized operator and reduce its scalar
part to a real transport direction. A local inverse, followed by an exact
correction of the full vector remainder, yields estimates with a fixed loss
of derivatives. A smoothing iteration then eliminates the actual return error
on a fixed open set of positive radius while retaining the seed as the leading
nonlinearity.
\end{abstract}
\maketitle
\markboth{Regular non-Berwald Landsberg surfaces}{Regular non-Berwald Landsberg surfaces}
\tableofcontents
\clearpage
\part{Construction and geometry}\label{part:geometry}
\section{Introduction and main results}\label{v4:introduction}

A Finsler metric assigns a strongly convex Minkowski norm to each tangent
space of a smooth manifold. Its canonical spray is positively homogeneous
of degree two in the tangent vector, but need not be quadratic. The Berwald
tensor consists of its third fibre derivatives. A Berwald metric has
vanishing Berwald tensor; the Landsberg condition requires a contraction
of that tensor to vanish. Whether every regular Landsberg metric is Berwald
is the Landsberg--Berwald problem, also called the unicorn problem. In this paper, we construct a family of smooth regular non-Berwald Landsberg surfaces. 

The geometric construction described in this paper was announced at the Conference on
Riemann--Finsler Geometry and Related Topics, held in Chongqing on
November 13--17, 2026. Mao subsequently posted the first version of his
preprint \cite{Mao} on arXiv on November 19, 2026.

Throughout
this paper regularity means smoothness and strong convexity in every
nonzero tangent direction. We use the conventions of \cite{BCS}; background
on the problem is given in \cite{Bao2007,Crampin2024}.
The qualification concerning all directions is essential. Almost regular
examples may have singular extremal directions, as in the
$(\alpha,\beta)$ constructions discussed in \cite{Shen2009,Elgendi2026}. Generalized
Landsberg structures on a three-dimensional manifold, studied through
coframings in \cite{Bryant1995,SSS2010}, must still be realized by closed strongly convex indicatrices
in the tangent spaces of a surface. In the construction below the remaining
condition is exact return of a periodically forced oscillator on an open
set of initial values. A single periodic trajectory does not supply that
condition.

Our higher-dimensional paper \cite{LiLiuAnalytic} concerns fibrewise
real-analytic Landsberg metrics in base dimension at least four. The
geometric construction from the oscillator arose independently in our work
on surfaces. We ask which periodic
seeds can be realized by regular Landsberg metrics and which features of
the resulting geometry retain the choice of seed. The normalized quadratic
seed used in Mao's construction \cite[Theorem 3.1]{Mao},
$1+\sqrt{70/13}\cos4s$, is the first single-frequency member of the class
considered here. Higher frequencies and their mixtures provide a setting
in which these questions can be studied systematically.

An averaging argument proposed by Szab\'o \cite{Szabo2008} contains a gap
concerning nonlinear parallel transport, as explained by Matveev
\cite{Matveev2009} and acknowledged in Szab\'o's correction
\cite{SzaboCorrection}. Preservation of the fibrewise Hessian metric
must therefore be distinguished from linearity of parallel transport.

The identification of the canonical connection is particularly important
for flat constructions. The example proposed in \cite{Zhou2014} was shown
to be Berwald in \cite{ElgendiYoussef2021}. The rigidity result in
\cite[Theorem 4.6]{ElgendiMuzsnay2024} concerns nonzero flag curvature
with an additional derivative condition, whereas the metrics below have
zero flag curvature. Related questions about Berwald curvature, holonomy,
and Landsberg metrizability are studied in \cite{ElgendiMuzsnay2025}.
The closed three-dimensional case is treated in the preprint
\cite{MaoZhou2026}. Our construction is local on a surface and is smooth in every
nonzero tangent direction.

\subsection{The geometric result and the meaning of the parameters}
The seeds are real smooth $\pi$-periodic functions of mean one. Three real
balance equations cancel the first nontrivial return obstruction, and a
positive coefficient obtained from the seed gives solvability of the
higher formal moment equations. We denote this class by $\mathcal A^\infty$.
For the local analytic construction we use the subclass
$\mathcal B^\infty$: at some nonzero initial value, the associated scalar
symbol has exactly one simple characteristic line. The definitions are
given in Sections~\ref{v4:seed-section} and~\ref{v4:char-section}.
They concern explicit functions of one variable and can be checked before
the closure equation is solved.

\begin{theorem}[Families of regular Landsberg surfaces]\label{v4:main-theorem}
For every $h\in\mathcal B^\infty$, one can choose a nonzero initial value
$x_h$, an open ball $U_h\Subset\R^2\setminus\{0\}$ about $x_h$, and
$\varepsilon_h>0$ such that, for every $0<|\varepsilon|<\varepsilon_h$,
there exists a smooth strongly convex Finsler metric $F_{\varepsilon,h}$
on $U_h$ with the following properties:
\begin{enumerate}[label=\textup{(\roman*)}]
\item It is regular on the entire slit tangent bundle $TU_h\setminus0$.
\item Its Landsberg tensor and flag curvature vanish.
\item Its Berwald tensor is nonzero at every base point in some tangent
direction. Its tangent norms are nonreversible and non-Riemannian.
\item It is obtained from an actual smooth periodic oscillator control
$H_{\varepsilon,h}$ whose return map is exactly the identity on $U_h$.
On a fixed sampling cylinder this control satisfies
$H_{\varepsilon,h}(q,s)=\varepsilon h(s)q^2+o(|\varepsilon|)$ in $C^2$.
\end{enumerate}
\end{theorem}

Here a family means a collection indexed by the admissible seeds and
sufficiently small nonzero amplitudes. The ball and the amplitude threshold
may depend on the seed. For each fixed finite order, the formal preparation
can be chosen smoothly near a reference seed; the local linear estimates
also hold on a common low-norm neighborhood. The theorem asserts existence
of the completed metrics for each parameter, without a joint smooth
selection or an extension through $\varepsilon=0$.

\begin{corollary}[Infinitely many distinct local geometries]
\label{v41:intro-nonhomothetic}
The metrics in Theorem~\ref{v4:main-theorem} contain a countably infinite
collection such that no nonempty base open subset of one is homothetic
to a base open subset of another. In particular they give countably
infinitely many pairwise nonisometric germs of regular non-Berwald
Landsberg surfaces.
\end{corollary}

The corollary uses one sufficiently small amplitude for each
even-frequency seed. The uniform separation of their Fourier modes
allows the geometric distinction to hold for every pair at once
(Corollary~\ref{v41:infinite-nonhomothetic}).

The amplitude also records a choice of scale. For $\varepsilon>0$, put
$\widetilde H(q,s)=\varepsilon H_{\varepsilon,h}(q/\varepsilon,s)$
and $D_\varepsilon x=\varepsilon x$. The rescaled flow is
$D_\varepsilon\Psi_{H_{\varepsilon,h}}(x/\varepsilon,s)$, so its
differential at $\varepsilon x$ equals the original differential at $x$.
The two indicatrices therefore have the same coordinate vectors at these
points, and the corresponding metric satisfies
$D_\varepsilon^*\widetilde F=\varepsilon F_{\varepsilon,h}$ on $U_h$.
The rescaled control closes on $\varepsilon U_h$ and can itself depend
on $\varepsilon$.

The seed remains visible in the leading Berwald tensor. In the fixed initial
coordinates, let $e_s=(\cos s,\sin s)^{\mathsf T}$ and
$m_s=(-\sin s,\cos s)^{\mathsf T}$. On compact subsets of the chosen base,
the tensor evaluated at a fixed direction has the expansion
\begin{equation}\label{v4:intro-leading-tensor}
 B^a{}_{ijk}(x,re_s)
 =\frac{2\varepsilon h(s)}r\,
      m_s^a m_{s,i}m_{s,j}m_{s,k}+o(|\varepsilon|),\qquad r>0.
\end{equation}
The error is uniform in $x$ and $s$ for fixed $r$. With the order-three
control estimate of Theorem~\ref{v4:exact-closure}, the error for the
chosen realizations is $O_{h,r}(\varepsilon^2)$. The seed is also the
uniform limit of an exact geometric invariant: the angular profile of
the horizontal derivative $\mathcal I_1$ of the Cartan scalar, divided by
its fibre average. This normalized profile is preserved by local
isometries and constant homotheties, up to an angular translation or reflection
(Proposition~\ref{v41:profile-invariant}). Thus two fixed seeds outside the
same translation--reflection orbit give locally nonhomothetic metrics
when their respective amplitudes are sufficiently small
(Proposition~\ref{v41:seed-separation}).

The class contains the seeds
\[
 h_4(s)=1+\sqrt{70/13}\cos4s,\qquad
 h_6(s)=1+\sqrt{210/11}\cos6s,
\]
and, more generally, the single-frequency seeds of every even frequency
$n\ge4$. For $n\ge6$, the initial value may be any prescribed nonzero
point; the control and its closing neighborhood are then constructed for
that chosen point. The characteristic condition is
open relative to $\mathcal A^\infty$ in a suitable finite $C^r$ topology.
For any $J\ge1$, a sufficiently small neighborhood of $h_4$
in the coefficient ellipsoid
\begin{equation}\label{v4:intro-ellipsoid}
 h=1+\sum_{j=1}^J(a_j\cos4js+b_j\sin4js),\qquad
 \sum_{j=1}^J
  \frac{16j^2-3}{2(16j^2-1)(16j^2-9)}(a_j^2+b_j^2)=\frac13
\end{equation}
therefore gives a $(2J-1)$-dimensional set of allowed seed parameters.
These dimensions refer to the seed spaces, not to a moduli space of
pairwise nonisometric metrics.
The two-frequency $(4,10)$ ellipsoid gives another explicit example:
its single-frequency endpoints have relative neighborhoods consisting
of admissible seeds, including members with both coefficients nonzero
(Corollary~\ref{v43:char-four-ten-endpoints}). There are also seeds
with frequencies $2,4,6$ for which distinct interactions cancel in
the balance equations (Proposition~\ref{v44:coupled-realization}). The full coefficient
ellipsoids constructed in Section~\ref{v4:seed-section} satisfy the formal
balance conditions; their entire membership in $\mathcal B^\infty$
remains to be determined.

\subsection{From indicatrices to the return equation}
Mu, Miao and Huang relate regularity in a rotationally symmetric flat
Finsler construction to isochronous centers of a planar vector field
\cite{MuMiaoHuang2025}. In the periodically forced setting below, we
prescribe a common return time on an open set of initial values. The
oscillator's special dependence on the state variable also gives the
Landsberg condition.
Take the initial-value plane of
\begin{equation}\label{v4:intro-oscillator}
 p_s=q,\qquad q_s=-p+H(q,s),\qquad \Psi_H(x,0)=x
\end{equation}
as the base surface. At a fixed initial point $x$, the inverse differential
of the flow traces the candidate unit curve
\[
 u(x,s)=(D_x\Psi_H(x,s))^{-1}e_1.
\]
If the original flow returns exactly for every initial point of an open
set $U$, differentiation of that identity and time periodicity give a
smooth periodic variational frame. Near the harmonic oscillator the curve
winds once. The variational equation then gives positive curvature and
one intersection with every positive ray. Defining $F(x,ru(x,s))=r$
produces a norm on every nonzero tangent vector over $U$.

The resulting connection must be identified with the canonical connection
of this norm. We do this by computing the energy Hessian and the
Euler--Lagrange equations. The full Berwald tensor is
\begin{equation}\label{v4:intro-exact-B}
 B^a{}_{ijk}(x,ru(x,s))
 =\frac{H_{qq}(q_H(x,s),s)}r\,
       w^a(x,s)\beta_i(x,s)\beta_j(x,s)\beta_k(x,s),
\end{equation}
where $w=(D_x\Psi_H)^{-1}e_2$ and $\beta=\dd_xq_H$.
Its Landsberg contraction vanishes. The horizontal distribution is
integrable and the flag curvature is zero. Smallness of the closure
correction preserves a nonzero coefficient in \eqref{v4:intro-exact-B}.
These calculations occupy the first part of the paper.

We take $U$ near a nonzero initial value. The auxiliary oscillator
trajectories from $U$ lie in a common flow tube, whose projection to the
control cylinder is the domain on which $H$ is needed. Exact return
is imposed for every initial point in $U$.

\subsection{How exact return is obtained}
Mao treats the fourfold seed by the same oscillator return equation and
reconstructs the metric from the dual indicatrix \cite[Section 8]{Mao}.
At the level of the formal return, the two proofs share a mixed curl
identity \cite[Section 3]{Mao}. For local exact return, both arguments
use a periodic comparison operator, normalize its two columns, and
correct the full remainder after solving a transport equation. A
smoothing iteration then gives exact return
\cite[Sections 2 and 4--7]{Mao}. For a variable seed
we determine the balance equations, prove surjectivity of the formal
moment maps, and identify the characteristic condition that permits
the local inverse at a nonzero initial value.

The seed balance equations and a one-dimensional maximum principle give
formal corrections of every order. We use only a finite truncation.
Rescaling it by a small nonzero amplitude produces a genuine smooth
control whose actual return error is of arbitrarily high prescribed order
in that amplitude. The truncation order is chosen after the derivative
loss and the finite parameter costs of the inverse have been determined.
Convergence of a formal series is not needed.

For the linear correction we first build a periodic comparison family from
the actual flow. Its endpoint integral has two frequency branches. The
complete target columns agree at zero amplitude; their normalized difference
therefore extends to that amplitude, including its finite-frequency terms.
At a seed in $\mathcal B^\infty$, the resulting scalar characteristic set
has a real transport direction near the selected nonzero initial value.
After a finite elliptic reduction, a short-interval transport inverse is
combined with the full vector equation. A positive Sobolev gain of the
complete error allows one final localization to remove it exactly.
The symbol estimates are uniform for the unrestricted periodic corrections
used in the iteration. The derivative losses and the powers of the
amplitude are fixed before the finite approximation is selected.

The comparison family is generally not an oscillator flow. Its discrepancy
from the true endpoint differential contains the current return residual.
The smoothing iteration retains this discrepancy, the background smoothing
error, the source smoothing error, and the nonlinear Taylor remainder.
It produces one sequence converging in every finite smooth norm on a base
ball of positive limiting radius. Its limit makes the original return
map exactly the identity. The $C^2$ size of the correction is small compared
with the amplitude, which supplies the geometric nonvanishing assertions.

\subsection{Contributions, scope, and organization}
Exact return on a fixed open set forces the three balance equations
for every quadratic seed with the asymptotics of Theorem~\ref{v4:main-theorem}
(Corollary~\ref{v44:balance-necessary}). A maximum principle gives
surjectivity of the moment maps for seeds with positive drift. A
highest-frequency calculation gives surjectivity for trigonometric
polynomial seeds of highest frequency at least four even when the drift
changes sign. These are statements about formal preparation.
The seed-dependent characteristic symbol
then supplies a criterion for local exact return, including a total
angular-variation test and single-frequency seeds of every even frequency
at least four. The normalized Cartan-derivative profile distinguishes fixed
seed orbits and yields countably many pairwise locally nonhomothetic
surface germs. The geometric estimates also prove nonreversibility of
every tangent norm in these small-amplitude constructions.

The formal class $\mathcal A^\infty$ and the characteristic subclass
$\mathcal B^\infty$ serve different steps of the argument; their equality
is not asserted. The completed control is $2\pi$-periodic and retains the
stated leading asymptotics. The half-period symmetry and exact central
jets used in the finite preparation are not required of that control.
The resulting examples give a starting point for studying how seed
variation changes a regular Landsberg surface. Determining the full
realizable seed class and the geometric equivalences among its
realizations are further questions.

Part~\ref{part:geometry} gives the geometric construction and its tensors.
Part~\ref{part:seeds} defines the seed families, constructs finite-order
approximations, and verifies the nonzero-point characteristic condition.
Part~\ref{part:analysis} constructs the local right inverse in
Theorem~\ref{v4:local-inverse}. In Part~\ref{part:realization} this
inverse and the finite approximation give exact return on a fixed open
ball (Theorem~\ref{v4:exact-closure}). Geometric reconstruction then
gives Theorem~\ref{v4:main-theorem}, and the angular profile separates
the metrics in Corollary~\ref{v41:intro-nonhomothetic}.
The appendices record the fixed-scale tools
and complete-kernel estimates used in these arguments.
\section{Exact return and the construction of indicatrices}
\label{v4:geom-chapter}

We construct the indicatrices from the variational flow of the oscillator and compute the resulting geometry.
Here $U$ is a connected open set of initial values; it need not contain the origin.
We first show what exact return on $U$ supplies, without making an existence assertion for the closing control.
The analytic construction in the later parts will supply this control on an open set of nonzero initial values.

\subsection{Geometric objectives and basic notions}

On an open set $U\subset\R^2$, a regular Finsler metric is a function
$F:TU\to[0,\infty)$ that is smooth on the slit tangent bundle $TU\setminus0$,
continuous on the zero section, and satisfies
\[
 F(x,\lambda y)=\lambda F(x,y)\quad(\lambda>0),\qquad
 F(x,y)>0\quad(y\ne0).
\]
Its energy and fundamental tensor are
\[
 E(x,y)=\tfrac12F(x,y)^2,\qquad
 g_{ij}(x,y)=\partial_{y^i}\partial_{y^j}E(x,y)
\]
and we require $g_{ij}(x,y)$ to be positive definite for every $y\ne0$.
This is the strong convexity condition.
At a fixed base point, $F_x$ is a possibly asymmetric Minkowski norm; its indicatrix is
\[
 \Sigma_x=\{y\in T_xU:F(x,y)=1\}.
\]
Reversibility, $F(x,-y)=F(x,y)$, is an additional condition and is not part of our definition of a metric.
The construction below includes every nonzero direction over the chosen open set in the base.

The starting point of the geometric construction is a family of indicatrices
$s\mapsto u(x,s)\in T_xU$, from which length is defined along each positive ray:
the vector $ru(x,s)$ has length $r$.
For this definition to work, each indicatrix must be smooth and closed,
meet every positive ray exactly once, and give a radial energy with positive definite Hessian in linear fiber coordinates.
We derive these properties from the variational equation of the oscillator.

Furthermore, the variation of the indicatrices with the base point must be compatible with the metric's own canonical connection.
This connection is determined by the Euler--Lagrange equation of the energy.
Writing the geodesic equation as $\ddot x^a+2G^a(x,\dot x)=0$,
the Berwald and Landsberg tensors are respectively
\[
 B^a{}_{ijk}=\partial_{y^i}\partial_{y^j}\partial_{y^k}G^a,
 \qquad
 L_{ijk}=-\tfrac12E_{y^a}B^a{}_{ijk}.
\]
The condition $B=0$ means that the spray is a quadratic polynomial in the fiber variables; this is the Berwald condition.
The condition $L=0$ requires only the displayed contraction to vanish.
Our objective is a metric with $L=0$, zero flag curvature, and $B\ne0$.
We compute these quantities from the energy.

\subsection{The oscillator, closure condition, and geometric realization theorem}

Take the initial-value plane of the oscillator as the base, with initial point $x=(x^1,x^2)$.
At time $s$, the solution map $\Psi_s:x\mapsto(p(x,s),q(x,s))$
sends an initial point to the state plane.
Its inverse differential pulls the first standard basis vector of the state plane back to a vector at the base point $x$:
\[
 u(x,s)=(D_x\Psi_s)^{-1}\binom10.
\]
The curve obtained by fixing $x$ and varying $s$ is precisely the proposed indicatrix $\Sigma_x$.
Thus the trajectory $s\mapsto\Psi(x,s)$ in the state plane
and the indicatrix $s\mapsto u(x,s)$ in the tangent space are different objects;
they are linked by the variational matrix obtained by differentiating with respect to the initial point.

Exact return of every initial point in $U$ at the common time $2\pi$ allows differentiation in the initial point,
yielding a closed variational frame and a full family of indicatrices.

\begin{theorem}[Geometric realization of exact return]\label{v4:reconstruction}
Let $U\subset\R^2$ be connected and open. Let
$W\subset\R_q\times\bigl(\R/(2\pi\mathbb Z)\bigr)$ be open, and let
$H\in C^\infty(W;\R)$.
Suppose that the equation
\begin{equation}\label{v4:geom-oscillator}
 p_s=q,\qquad q_s=-p+H(q,s),\qquad
 \Psi(x,0)=x,\quad \Psi=(p,q)^{\mathsf T}
\end{equation}
has solutions from every $x\in U$ through the full interval $[0,2\pi]$,
on a common open flow tube whose projection to the $(q,s)$ cylinder is
contained in $W$. Suppose also that
\begin{equation}\label{v4:geom-exact}
 \Psi(x,2\pi)=x\qquad(x\in U).
\end{equation}
Set $A=D_x\Psi$ and $u=A^{-1}(1,0)^{\mathsf T}$.
Assume that the closed curve $u(x_0,\cdot)$ has winding number one about zero
for one $x_0\in U$. Then
\begin{equation}\label{v4:geom-defF}
 F(x,ru(x,s))=r\quad(r>0),\qquad F(x,0)=0
\end{equation}
defines a regular, strongly convex Finsler metric on $TU$, smooth in all nonzero directions.
Here smoothness means joint smoothness on $TU\setminus0$, with continuity required on the zero section.
Its fundamental tensor is $g_{ru}=\alpha\otimes\alpha+\beta\otimes\beta=A^{\mathsf T}A$,
where $\alpha=\dd_xp$ and $\beta=\dd_xq$.
Its flag curvature and Landsberg tensor both vanish.
In natural fiber coordinates, its Berwald tensor is
\begin{equation}\label{v4:geom-Bmain}
 B^a{}_{ijk}(x,ru)=
 \frac{H_{qq}(q(x,s),s)}{r}\,
 w^a(x,s)\beta_i(x,s)\beta_j(x,s)\beta_k(x,s),
\end{equation}
where $w=A^{-1}(0,1)^{\mathsf T}$ and $\beta=\dd_xq$.
In particular, if $H_{qq}(q(x,s),s)\ne0$ for some $x\in U$ and some $s$,
then the metric is non-Berwald.
\end{theorem}

The control is defined near every sampled velocity and time.
The base set $U$ consists of initial values, while the flow tube contains
the intermediate states and a neighborhood of them. These states may leave
$U$ during the period. The theorem imposes no normalization of $H$ at $q=0$.
The winding assumption will follow from comparison with the harmonic oscillator
in Lemma~\ref{v4:geom-small-control}.
The construction and Hessian calculation below prove the Finsler assertions;
Section~\ref{v4:geom-tensors} proves the assertions about the canonical spray and its curvature.

\begin{remark}[Periodic coefficients and closure]
Periodic coefficients alone do not imply exact return. For example, take
$H(q,s)=\varepsilon q^5$, $\varepsilon>0$. Solutions with sufficiently small
initial values exist through time $2\pi$, but satisfy
\[
 \frac{\dd}{\dd s}\frac{p^2+q^2}{2}=\varepsilon q^6.
\]
The integral over one period is positive for every nonzero initial value:
if $q$ vanished identically, the oscillator equation would also give $p=0$.
Thus these trajectories do not return to their initial values.
\end{remark}

\subsection{The variational matrix and frames pulled back from state coordinates}

\begin{lemma}[Smooth dependence on a common flow tube]\label{v4:geom-smooth-flow}
For a smooth control $H$, the solution of \eqref{v4:geom-oscillator} is jointly
smooth in initial value and time on every open common flow domain.
The derivatives at different orders belong to the same solution.
\end{lemma}
\begin{proof}
Write the integral equation on a short time interval within a relatively
compact flow tube. Its linearization in the continuous path variable is the
identity minus the Volterra operator $v\mapsto\int D_zf\,v$, where
$f(s,z)=(q,-p+H(q,s))$. Choose the interval so that its length times a bound
for $D_zf$ is less than one. The operator is invertible by its Neumann series,
and the smooth implicit function theorem gives dependence on the initial
point at each finite differentiability order. Uniqueness identifies these
solutions. A finite number of short intervals covers any compact time
interval in the tube, and their solution maps compose. The differential
equation supplies time derivatives and all mixed derivatives.
\end{proof}

Write the initial point as $x=(x^1,x^2)$ and the state coordinates as $(p,q)$.
At fixed $s$, the differential of $x\mapsto\Psi(x,s)$ is
\[
 A(x,s)=
 \begin{pmatrix}
 \partial_1p&\partial_2p\\
 \partial_1q&\partial_2q
 \end{pmatrix}.
\]
Differentiating \eqref{v4:geom-oscillator} with respect to the initial point gives
\begin{equation}\label{v4:geom-variation}
 A_s=
 \begin{pmatrix}0&1\\-1&I\end{pmatrix}A,\qquad
 A(x,0)=\Id,\qquad
 I(x,s)=H_q(q(x,s),s).
\end{equation}
Here $I$ is a scalar function; the identity matrix is denoted by $\Id$.
The determinant satisfies
\[
 \partial_s\det A=I\det A,\qquad
 \det A(x,s)=
 \exp\left(\int_0^s I(x,t)\,\dd t\right)>0.
\]
The variational matrix is therefore invertible at every time.

Define the time-dependent coframe and its dual frame by
\begin{equation}\label{v4:geom-frames}
 \alpha=\dd_xp,\qquad \beta=\dd_xq,\qquad
 u=A^{-1}\binom10,\qquad w=A^{-1}\binom01.
\end{equation}
They satisfy
\begin{equation}\label{v4:geom-duality}
 \alpha(u)=1,\quad \alpha(w)=0,\quad
 \beta(u)=0,\quad \beta(w)=1.
\end{equation}
At fixed $s$, the functions $p,q$ are local coordinates on the base,
and $u,w$ are the corresponding coordinate vector fields.
Their commutator in the base variables is therefore zero.
The identity also follows directly from $\dd_x\alpha=\dd_x\beta=0$:
\[
 0=\dd_x\alpha(u,w)
 =u(\alpha(w))-w(\alpha(u))-\alpha([u,w])
 =-\alpha([u,w]).
\]
The same calculation with $\beta$ yields $[u,w]=0$.

Differentiating the inverse matrix in time gives
$(A^{-1})_s=-A^{-1}\left(\begin{smallmatrix}0&1\\-1&I\end{smallmatrix}\right)$,
and hence
\begin{equation}\label{v4:geom-frameODE}
 u_s=w,\qquad w_s=-u-Iw,\qquad [u,w]=0.
\end{equation}
On the other hand, differentiating the original oscillator directly with respect to the initial point gives
\begin{equation}\label{v4:geom-coframeODE}
 \alpha_s=\beta,\qquad
 \beta_s=-\alpha+I\beta,\qquad
 \dd_x\alpha=\dd_x\beta=0.
\end{equation}
These two systems supply all the input for the geometric computations below.

For comparison, the harmonic oscillator $H=0$ has
\begin{equation}\label{v4:geom-harmonic-frame}
 C_s=\begin{pmatrix}\cos s&\sin s\\-\sin s&\cos s\end{pmatrix},
 \quad e_s=(\cos s,\sin s)^{\mathsf T},\quad
 m_s=(-\sin s,\cos s)^{\mathsf T},
\end{equation}
with $\Psi(x,s)=C_sx$, $A=C_s$, $u=e_s$, and $w=m_s$ at every base point.
This comparison is made on the chosen initial-value set; it does not require
an actual metric fiber over the origin.

\subsection{Closure, winding number, and coverage of all directions}

Actual return gives smooth matching across the parameter seam.
The winding assumption and the orientation in \eqref{v4:geom-frameODE}
then ensure that the indicatrix meets every positive ray exactly once.

\begin{lemma}[Parametrization of all directions by the indicatrix]\label{v4:geom-directions}
Under the hypotheses of Theorem~\ref{v4:reconstruction},
$s\mapsto u(x,s)$ is a smooth closed curve for every $x\in U$.
The direction map
\[
 s\longmapsto\frac{u(x,s)}{|u(x,s)|}
\]
is an orientation-preserving diffeomorphism from the parameter circle $\R/(2\pi\mathbb Z)$ onto the Euclidean circle of directions.
\end{lemma}
\begin{proof}
For $s$ near zero, the maps $s\mapsto\Psi(x,2\pi+s)$ and
$s\mapsto\Psi(x,s)$ solve the same nonautonomous equation, because its
coefficients are $2\pi$-periodic. Their initial values agree by
\eqref{v4:geom-exact}, so uniqueness gives
\[
 \Psi(x,2\pi+s)=\Psi(x,s).
\]
Local existence and smooth dependence make this an identity of smooth maps
on a neighborhood in $(x,s)$ of each point $(x,0)$.
Repeated translation therefore extends the solution periodically to all $s$;
the extensions agree on overlaps by uniqueness.
Differentiating this identity gives periodicity of $A,u,w$ and every mixed
derivative. Thus the parameter seam is smooth. Here it is essential that
\eqref{v4:geom-exact} holds on an open set of initial points.

By \eqref{v4:geom-frames}, $u,w$ are the two columns of $A^{-1}$, hence
\[
 \det(u,u_s)=\det(u,w)=\frac1{\det A}>0.
\]
In particular, $u\ne0$ and $u_s\ne0$.
Locally write $u=R(\cos\theta,\sin\theta)$, where $R=|u|>0$.
Then
\[
 \theta_s=\frac{\det(u,u_s)}{|u|^2}>0.
\]
Thus the direction angle is strictly increasing. Its increment over one period is $2\pi k(x)$,
where $k(x)$ is a positive integer.
It can be written as the continuous function
\[
 k(x)=\frac1{2\pi}\int_0^{2\pi}
       \frac{\det(u,u_s)}{|u|^2}\,\dd s.
\]
Because $U$ is connected, an integer-valued continuous function there is constant;
the winding assumption gives $k(x_0)=1$, so $k(x)=1$ everywhere.
A strictly monotone direction map that winds once is bijective,
and its derivative is nonzero, so it is a diffeomorphism.
\end{proof}

\begin{lemma}[A small-control criterion for winding number one]\label{v4:geom-small-control}
Let $\overline U$ be compact, put $R=\max_{\overline U}|x|$, and choose
$b>R$. Suppose that $H$ is smooth on a neighborhood of
$[-b,b]\times\R/(2\pi\mathbb Z)$ and
\[
 2\pi\|H\|_{C^0}<b-R,\qquad
 2\pi\|H_q\|_{C^0}<\log2
\]
on that cylinder. Then all trajectories from $\overline U$ exist for a full
period and $\sup_{x,s}|u(x,s)-e_s|<1$. If exact return holds on $U$,
its closed indicatrices have winding number one.
\end{lemma}
\begin{proof}
Variation of constants gives, as long as the trajectory lies in $|q|<b$,
\[
 \Psi(x,s)=C_sx+\int_0^s C_{s-t}\binom{0}{H(q(x,t),t)}\dd t.
\]
Since $C_s$ is orthogonal, its distance from $C_sx$ is at most
$2\pi\|H\|_{C^0}$. Thus the trajectory remains strictly inside $|z|<b$;
continuation gives the claimed full interval. The strict margin also gives
a common open flow tube around $\overline U$.

Write $A=C_sD$. The variational equation gives
\[
 D_s=C_s^{-1}\begin{pmatrix}0&0\\0&I\end{pmatrix}C_sD,
 \qquad D(x,0)=\Id.
\]
The coefficient has operator norm at most $\eta:=\|H_q\|_{C^0}$.
The inverse equation and Gronwall's inequality give
$\|D^{-1}-\Id\|\le e^{2\pi\eta}-1<1$.
As $u=D^{-1}e_s$, the same bound applies to $|u-e_s|$.
If exact return holds, $u(x,\cdot)$ is closed by the preceding proof.
The straight homotopy $(1-t)e_s+tu(x,s)$ never vanishes, since its distance
from $e_s$ is less than one. It is a homotopy of closed nonvanishing curves,
so their winding numbers agree and equal one.
\end{proof}

Every nonzero tangent vector therefore has a unique representation
$y=ru(x,s)$, with $r>0$ and $s$ understood modulo $2\pi$.
The map
\begin{equation}\label{v4:geom-polar}
 \mathcal P:U\times(0,\infty)\times\R/(2\pi\mathbb Z)
 \longrightarrow TU\setminus0,\qquad
 (x,r,s)\longmapsto(x,ru(x,s))
\end{equation}
has Jacobian determinant
$r\det(u,w)>0$ in a local angular coordinate.
Its local inverses are therefore smooth, and bijectivity makes them agree to form a global inverse.
Formula~\eqref{v4:geom-defF} consequently defines a jointly smooth function on $TU\setminus0$,
which is positive and positively homogeneous of degree one.

\subsection{Computing the actual fiber Hessian from the radial energy}

Fix a base point $x$ and set $E(y)=\tfrac12F(y)^2$.
Since $E(ru)=r^2/2$, differentiation along the radial direction and along the indicatrix gives
\[
 \dd E_u(u)=1,\qquad \dd E_u(w)=0.
\]
Thus $\dd E_u=\alpha$. Write $g_y=D_y^2E$.
Homogeneity of degree two gives
\[
 g_u(u,\cdot)=\dd E_u,\qquad
 g_u(u,u)=1,\qquad g_u(u,w)=0.
\]
Now differentiate the identity $\dd E_{u(s)}(w(s))=0$:
\[
 0=g_u(w,w)+\dd E_u(w_s)
   =g_u(w,w)+\alpha(-u-Iw).
\]
The duality relations give $g_u(w,w)=1$.
Consequently,
\begin{equation}\label{v4:geom-metric}
 g_{ru}=g_u=\alpha\otimes\alpha+\beta\otimes\beta.
\end{equation}
The first equality follows from the degree-zero homogeneity of $g$.
The vectors $u,w$ form an orthonormal basis for $g_{ru}$, so the actual fiber Hessian is positive definite.

At fixed $x$, we have $\dd y=u\,\dd r+rw\,\dd s$;
substituting into \eqref{v4:geom-metric} gives
\begin{equation}\label{v4:geom-cone}
 g=\dd r^2+r^2\dd s^2.
\end{equation}
This is an expression in nonlinear polar coordinates; the components $g_{ij}$ in linear fiber coordinates may vary.
In particular, the Hessian arc-length element on the indicatrix $r=1$ is $\dd s$,
and its circumference is exactly $2\pi$.
For each fixed $x$, the map
\[
 ru(x,s)\longmapsto r(\cos s,\sin s)
\]
is a Hessian isometry from the punctured fiber to the punctured Euclidean plane.
The isometry is generally nonlinear in the linear structure of $T_xU$.
Likewise, $u_s=w$ and $u_{ss}=-u-Iw$ imply
\[
 \det(u_s,u_{ss})=\det(w,-u-Iw)=\det(u,w)>0,
\]
so the indicatrix also has everywhere positive signed curvature in natural linear fiber coordinates.

To include the zero section, observe that on every compact subdomain of the base,
the continuous nonvanishing function $|u(x,s)|$ has a positive lower bound $c$.
Since $|y|=r|u|$, we obtain $F(x,y)\le c^{-1}|y|$;
thus setting $F(x,0)=0$ gives continuity.
Along a line avoiding the origin, the second derivative of $E$ is controlled by its positive definite Hessian, so it is convex.
On a line through the origin with direction $v$,
\[
 E(tv)=
 \begin{cases}
  t^2E(v),&t\ge0,\\
  t^2E(-v),&t\le0.
 \end{cases}
\]
The one-sided first derivatives at zero both vanish, so the joined function remains convex.
Hence the unit ball $\{F\le1\}$ is a convex body with the origin in its interior,
and $F$ is its Minkowski functional and satisfies the triangle inequality.
This verifies all the fiber conditions required for a regular Finsler metric.
\section{The canonical connection and curvature tensors}\label{v4:geom-tensors}

The metric constructed in Section~\ref{v4:geom-chapter} has a connection
determined by its energy. We compute this connection and its curvature,
then prove stability of non-Berwaldness and nonreversibility under the
small closing corrections used later.

\subsection{Two tangent-bundle coordinate systems and the canonical spray}

We first recall how the energy determines the connection.
For the velocity $y=\dot x$ of a curve, the Euler--Lagrange equation expands as
\[
 0=\frac{\dd}{\dd t}E_{y^j}(x,\dot x)-E_{x^j}(x,\dot x)
  =g_{ja}\ddot x^a+E_{x^ky^j}y^k-E_{x^j}.
\]
Positive definiteness of $g$ allows us to solve for the acceleration. Thus
\[
 G^a=\tfrac12g^{aj}(E_{x^ky^j}y^k-E_{x^j}),\qquad
 N_i^a=\partial_{y^i}G^a,\qquad
 \delta_i=\left.\partial_{x^i}\right|_y-N_i^a\partial_{y^a}.
\]
The functions $G$ and $N$ are homogeneous in $y$ of degrees two and one, respectively.
The horizontal fields $\delta_i$ are thus determined by the energy.

In the base derivatives below, the distinction between coordinate systems is essential:
$\partial_{x^i}|_y$ in natural coordinates holds the fiber components fixed,
whereas $\partial_{x^i}|_{r,s}$ in polar coordinates holds $r,s$ fixed.
Taking the total differential of $y^a=ru^a(x,s)$ gives
\begin{equation}\label{v4:geom-dy}
 \dd y^a=u^a\dd r+rw^a\dd s
             +r(\partial_i u^a)\dd x^i.
\end{equation}
Here $\partial_i u^a$ always keeps the parameter $s$ fixed.
Contraction with the respective coframe components $\alpha_a,\beta_a$ gives
\begin{align}
 \dd r&=\alpha_a\dd y^a-r\alpha_a\partial_i u^a\dd x^i,
 \label{v4:geom-dr}\\
 \dd s&=\frac{\beta_a}{r}\dd y^a-\beta_a\partial_i u^a\dd x^i.
 \label{v4:geom-ds}
\end{align}
The fiber derivatives are therefore
\begin{equation}\label{v4:geom-vertical}
 \partial_{y^a}=\alpha_a\partial_r+\frac{\beta_a}{r}\partial_s,
\end{equation}
and the base coordinate vectors satisfy
\begin{equation}\label{v4:geom-basechange}
 \left.\partial_{x^i}\right|_{r,s}
 =\left.\partial_{x^i}\right|_y
       +r\partial_i u^a\,\partial_{y^a}.
\end{equation}

Since $E=r^2/2$, equation~\eqref{v4:geom-dr} gives
\begin{equation}\label{v4:geom-energyderivatives}
 E_{y^j}=r\alpha_j,\qquad
 E_{x^j}=-r^2\alpha_a\partial_j u^a.
\end{equation}
The second expression is a partial derivative at fixed natural fiber coordinate $y$.
These energy derivatives determine the canonical connection as follows.

\begin{proposition}[The canonical spray and nonlinear connection]\label{v4:geom-connection}
The spray and connection coefficients determined by the energy Euler--Lagrange equation are
\begin{equation}\label{v4:geom-sprayN}
 G^a(x,ru)=-\frac{r^2}{2}u^i\partial_i u^a,\qquad
 N_j^a=\partial_{y^j}G^a=-r\partial_j u^a.
\end{equation}
The horizontal coordinate fields in polar coordinates are
\begin{equation}\label{v4:geom-horizontal}
 \delta_j=\left.\partial_{x^j}\right|_y-N_j^a\partial_{y^a}
          =\left.\partial_{x^j}\right|_{r,s}.
\end{equation}
\end{proposition}
\begin{proof}
Fix $r>0$ and $s$, and first consider a candidate curve in the base satisfying
$\dot x=r u(x,s)$.
Its velocity parameters $r,s$ are constant by construction.
By \eqref{v4:geom-energyderivatives},
\[
 \frac{\dd}{\dd t}E_{y^j}(x,\dot x)
 =r^2u^i\partial_i\alpha_j.
\]
Since $\alpha=\dd_xp$,
we have $\partial_i\alpha_j=\partial_j\alpha_i$;
differentiating $\alpha_i u^i=1$ with respect to $x^j$ yields
\[
 u^i\partial_i\alpha_j
 =u^i\partial_j\alpha_i
 =-\alpha_i\partial_j u^i.
\]
The candidate curve thus satisfies
$\frac{\dd}{\dd t}E_{y^j}=E_{x^j}$.
Every nonzero initial velocity has a unique representation $ru(x,s)$,
and strong convexity ensures that the Euler--Lagrange equation determines the acceleration and has a locally unique solution.
Hence these curves give all geodesics with nonzero initial velocity.
Along such a curve,
$\ddot x^a=r^2u^i\partial_i u^a$.
Comparison with $\ddot x^a+2G^a=0$ gives the first formula.

To compute $N_j^a$, temporarily write $K^a=u^i\partial_i u^a$.
Using $u_s=w$ and $[u,w]=0$,
\[
 \partial_sK^a
 =w^i\partial_i u^a+u^i\partial_iw^a
 =2w^i\partial_i u^a.
\]
Thus \eqref{v4:geom-vertical} gives
\[
 \begin{aligned}
 \partial_{y^j}G^a
 &=-r\alpha_j u^i\partial_i u^a
   -r\beta_jw^i\partial_i u^a\\
 &=-r(\alpha_j u^i+\beta_jw^i)\partial_i u^a
 =-r\partial_j u^a.
 \end{aligned}
\]
The last step is the dual-basis decomposition
$\alpha_j u^i+\beta_jw^i=\delta_j^i$.
Substituting into the definition of the horizontal fields
and then using \eqref{v4:geom-basechange} proves \eqref{v4:geom-horizontal}.
\end{proof}

\subsection{Geometric descriptions of geodesics and nonlinear parallel transport}

The proof of Proposition~\ref{v4:geom-connection} also describes every geodesic:
fix $r>0$ and $s_0$, and integrate the base vector field $ru(\,\cdot\,,s_0)$.
In the state coordinates at the fixed time $s_0$, namely $\Psi_{s_0}=(p(\,\cdot\,,s_0),q(\,\cdot\,,s_0))$,
this geodesic satisfies
\[
 \frac{\dd}{\dd t}p(x(t),s_0)=r\alpha(u)=r,\qquad
 \frac{\dd}{\dd t}q(x(t),s_0)=r\beta(u)=0.
\]
Thus the geodesics in a given direction family become horizontal straight lines in the corresponding state coordinates.
The state coordinates depend on the fixed direction parameter $s_0$; $t$ is the geodesic parameter.

Along an arbitrary prescribed base curve $x(t)$,
the parallel transport equation is
$\dot y^a+N_i^a\dot x^i=0$.
Substituting $y=ru(x,s)$ and \eqref{v4:geom-sprayN}
cancels the terms involving base derivatives and leaves
$u\dot r+rw\dot s=0$.
Linear independence of $u,w$ gives
$\dot r=\dot s=0$.
Parallel transport therefore preserves $F=r$
and, in the polar coordinates for the fiber Hessian, preserves
$\dd r^2+r^2\dd s^2$.
Explicitly, parallel transport from $x_0$ to $x_1$ is
\begin{equation}\label{v4:geom-transport}
 \mathcal T_{x_0,x_1}\bigl(ru(x_0,s)\bigr)=ru(x_1,s).
\end{equation}
As long as the base curve remains in $U$, the right-hand side depends only on its endpoints
and is defined simultaneously for all nonzero vectors.
It preserves Finsler length and pulls back the Hessian metric of one fiber to that of the other.
This gives a geometric interpretation of the Landsberg property; the next subsection also verifies $L=0$ directly.
The transport in \eqref{v4:geom-transport} is generally nonlinear in natural fiber coordinates.

\subsection{The Cartan, Berwald, and Landsberg tensors and flag curvature}

In natural fiber coordinates, the Cartan tensor is defined by
$C_{ijk}=\tfrac12\partial_{y^k}g_{ij}$.
Because $E$ is homogeneous of degree two and $g$ of degree zero,
every contraction of $C$ with the radial vector vanishes.
Differentiating $g_{u(s)}(w(s),w(s))=1$, with both the evaluation point and the arguments varying, gives
\[
 0=2C_u(w,w,w)+2g_u(w_s,w)
   =2C_u(w,w,w)-2I.
\]
Therefore
\begin{equation}\label{v4:geom-Cartan}
 C_{ijk}(x,ru)=\frac{I(x,s)}r\,\beta_i\beta_j\beta_k.
\end{equation}
The base derivatives satisfy
\begin{equation}\label{v4:geom-dI}
 \partial_i I(x,s)=H_{qq}(q(x,s),s)\beta_i(x,s),
\end{equation}
so the shapes of nearby fibers can still vary nontrivially.

The Berwald tensor is the third fiber derivative of the spray,
or equivalently the second fiber derivative of the connection coefficients.
By \eqref{v4:geom-sprayN} and \eqref{v4:geom-vertical},
\begin{equation}\label{v4:geom-firstverticalN}
 \partial_{y^j}N_i^a
 =-\alpha_j\partial_i u^a-\beta_j\partial_iw^a.
\end{equation}
The right-hand side is independent of $r$.
Its $s$ derivative expands in full as follows:
\[
 \begin{aligned}
 \partial_s(\partial_{y^j}N_i^a)
 ={}&-\beta_j\partial_i u^a-\alpha_j\partial_iw^a\\
 &-(-\alpha_j+I\beta_j)\partial_iw^a\\
 &-\beta_j[-\partial_i u^a-I\partial_iw^a
                         -(\partial_iI)w^a]\\
 ={}&\beta_j(\partial_iI)w^a.
 \end{aligned}
\]
All terms not containing $\partial_iI$ in the preceding expansion cancel exactly.
Using \eqref{v4:geom-vertical} and \eqref{v4:geom-dI} once more gives
\[
 B^a{}_{ijk}
 =\partial_{y^k}\partial_{y^j}N_i^a
 =\frac{\beta_k\beta_j}{r}(\partial_iI)w^a
 =\frac{H_{qq}}r\,w^a\beta_i\beta_j\beta_k.
\]
This proves the full tensor formula \eqref{v4:geom-Bmain}.

For the Landsberg tensor, we use the convention
$L_{ijk}=-\tfrac12E_{y^a}B^a{}_{ijk}$.
Since $E_{y^a}=r\alpha_a$ and $\alpha(w)=0$,
we immediately obtain
\[
 L_{ijk}=0.
\]
To identify the curvature of this same canonical spray, put
\[
 \mathcal R^a{}_{ij}=\delta_iN_j^a-\delta_jN_i^a,
 \qquad
 [\delta_i,\delta_j]=-\mathcal R^a{}_{ij}\partial_{y^a}.
\]
Equation~\eqref{v4:geom-horizontal} gives $\mathcal R^a{}_{ij}=0$.
The degree-two homogeneity of $G$ gives $y^iN_i^a=2G^a$;
differentiating it in $y^b$ gives
$y^i\partial_{y^b}N_i^a=N_b^a$. It follows that
\begin{align}
 y^i\mathcal R^a{}_{ki}
 &=2\partial_{x^k}G^a-y^i\partial_{x^i}N_k^a
   +2G^b\partial_{y^b}N_k^a-N_k^bN_b^a \notag\\
 &=2\partial_{x^k}G^a-y^i\partial_{x^i}\partial_{y^k}G^a
   +2G^b\partial_{y^b}\partial_{y^k}G^a
   -(\partial_{y^k}G^b)(\partial_{y^b}G^a).
 \label{v4:geom-jacobi-identity}
\end{align}
The last expression is the coordinate Jacobi curvature $R^a{}_k$ of the
canonical spray. Hence $R_y(v)=0$ for every $y\ne0$ and every $v$.
For linearly independent $y,v$, the flag curvature is
\begin{equation}\label{v4:geom-flag-identity}
 K(x,y,\operatorname{span}\{y,v\})
 =\frac{g_y(R_y(v),v)}
 {g_y(y,y)g_y(v,v)-g_y(y,v)^2}=0.
\end{equation}
The denominator is positive by positive definiteness and strict
Cauchy--Schwarz. This proves vanishing of the flag curvature in all directions.

\subsection{The seed in the leading geometry of the closed metrics}

The exact tensor formulas give a way of preserving the desired geometry
without prescribing a jet at the zero initial value. The smallness
needed here concerns two derivatives of the actual control, on the entire
velocity interval covered by the flow.

\begin{theorem}[Geometric stability of a quadratic seed]\label{v4:geometric-stability}
Let $K\Subset\R^2\setminus\{0\}$ be compact, let $b>\max_K|x|$, and let
$h$ be a nonzero real smooth $\pi$-periodic function. For sufficiently small
nonzero $\varepsilon$, let $H_{\varepsilon,h}$ be real smooth and
$2\pi$-periodic on a neighborhood of
$[-b,b]\times\R/(2\pi\mathbb Z)$. Suppose that
\begin{equation}\label{v4:geom-seed-approximation}
 \rho_h(\varepsilon):=
 \frac{\|H_{\varepsilon,h}(q,s)-\varepsilon h(s)q^2\|_{C^2}}
      {|\varepsilon|}
 \longrightarrow0
\end{equation}
on this cylinder as $\varepsilon\to0$. Let $U_{\varepsilon,h}$ be a nonempty
connected open subset of the interior of $K$, and assume that
$P_{H_{\varepsilon,h}}=\Id$ there.
Then, for all sufficiently small nonzero $\varepsilon$, the reconstruction
of Theorem~\ref{v4:reconstruction} gives a regular, strongly convex Landsberg
metric on $U_{\varepsilon,h}$ with zero flag curvature. Its Berwald tensor
is nonzero in some direction at every base point. Every tangent norm is
nonreversible and, in particular, non-Riemannian.

The smallness threshold may depend on $h,K,b$ and the rate at which
$\rho_h(\varepsilon)$ tends to zero. No continuity or smooth dependence of
the completed metrics on $\varepsilon$ or $h$ is part of this assertion.
\end{theorem}

A specific modulus for this dependence is
\begin{equation}\label{v43:geom-remainder-modulus}
 \rho_h^*(a)=\sup_{0<|\varepsilon|\le a}\rho_h(\varepsilon),
 \qquad \rho_h^*(a)\longrightarrow0\quad(a\downarrow0).
\end{equation}
Here $a$ is small enough that the controls are defined for every amplitude
in the supremum. The estimates below hold once
$|\varepsilon|+\rho_h^*(|\varepsilon|)$ is sufficiently small, with the
required bound depending on $h,K,b$.

\begin{proof}
Write $H_\varepsilon=H_{\varepsilon,h}$ and omit $h$ in the other subscripts.
By \eqref{v4:geom-seed-approximation}, $\|H_\varepsilon\|_{C^2}=O(|\varepsilon|)$.
Lemma~\ref{v4:geom-small-control} supplies the common full-period flow tube
and the winding number one. More quantitatively, variation of constants
and the variational equation give, uniformly for $x\in K$ and
$0\le s\le2\pi$,
\begin{equation}\label{v4:geom-harmonic-estimates}
 \begin{aligned}
 \Psi_\varepsilon(x,s)&=C_sx+O(|\varepsilon|),&
 A_\varepsilon(x,s)&=C_s+O(|\varepsilon|),\\
 u_\varepsilon(x,s)&=e_s+O(|\varepsilon|),&
 w_\varepsilon(x,s)&=m_s+O(|\varepsilon|),\\
 \alpha_\varepsilon(x,s)&=e_s^{\mathsf T}+O(|\varepsilon|),&
 \beta_\varepsilon(x,s)&=m_s^{\mathsf T}+O(|\varepsilon|).
 \end{aligned}
\end{equation}
For example, the integral expression in the proof of that lemma gives
the first estimate. Substituting $I_\varepsilon=H_{\varepsilon,q}=O(|\varepsilon|)$
in the equation for $D=C_s^{-1}A_\varepsilon$ gives the second estimate,
its inverse estimate, and hence the other four. The equation $u_s=w$
also controls the angular derivative. These are estimates for each actual
flow, and require no differentiation with respect to $\varepsilon$.

The control approximation further gives
\begin{equation}\label{v4:geom-I-approximation}
 \begin{split}
 I_\varepsilon(x,s)
    &=2\varepsilon h(s)(m_s\cdot x)+o(|\varepsilon|),\\
 H_{\varepsilon,qq}(q_\varepsilon(x,s),s)
    &=2\varepsilon h(s)+o(|\varepsilon|),
 \end{split}
\end{equation}
uniformly on the same set. The second remainder is bounded by
$|\varepsilon|\rho_h(\varepsilon)$, and the first by
$C|\varepsilon|(|\varepsilon|+\rho_h(\varepsilon))$; the additional term in
the first estimate comes from
$q_\varepsilon-m_s\cdot x=O(|\varepsilon|)$.
Consequently, choose $s_0$ with $h(s_0)\ne0$.
For all small nonzero $\varepsilon$, the second expression is nonzero
at every $x\in K$ when $s=s_0$. Formula~\eqref{v4:geom-Bmain} then proves
non-Berwaldness, since $w$ and $\beta$ are nonzero.

For the reversibility argument, the tensors must be compared at opposite
vectors in the same linear tangent space. Opposite parameter values of the
perturbed indicatrix need not be opposite vectors. Fix instead
$y=re_\theta$, with $r>0$. The direction map of the indicatrix is increasing
and, by \eqref{v4:geom-harmonic-estimates}, uniformly close to $s\mapsto e_s$.
Thus there is a unique parameter
$\sigma_\varepsilon(x,\theta)=\theta+O(|\varepsilon|)$ modulo $2\pi$ and a
radius $r_\varepsilon=r(1+O(|\varepsilon|))$ such that
\[
 re_\theta=r_\varepsilon
 u_\varepsilon(x,\sigma_\varepsilon(x,\theta)).
\]
The bounds are uniform on $K$ and in $\theta$; they follow directly by
normalizing $u_\varepsilon=e_s+O(|\varepsilon|)$ and using the positive
angular derivative, which is bounded below for small $|\varepsilon|$.
Inserting these parameters in \eqref{v4:geom-Cartan} and
\eqref{v4:geom-Bmain} yields, in the fixed coordinate bases,
\begin{align}
 C_\varepsilon(x,re_\theta)
 &=\frac{2\varepsilon}{r}h(\theta)(m_\theta\cdot x)
       (m_\theta^{\mathsf T})^{\otimes3}+o(|\varepsilon|)/r,
       \label{v4:geom-C-leading}\\
 B_\varepsilon(x,re_\theta)
 &=\frac{2\varepsilon}{r}h(\theta)
       m_\theta\otimes(m_\theta^{\mathsf T})^{\otimes3}
       +o(|\varepsilon|)/r.
       \label{v4:geom-B-leading}
\end{align}
Indeed, all changes in $h,m$ and the radial factor multiply a leading
factor $\varepsilon$, so their $O(|\varepsilon|)$ changes contribute
$O(\varepsilon^2)$; the original remainders remain uniform $o(|\varepsilon|)$.

Since $h(\theta+\pi)=h(\theta)$ and $m_{\theta+\pi}=-m_\theta$,
the leading Cartan expression is even under $\theta\mapsto\theta+\pi$.
To see that its coefficient cannot vanish on a whole fiber, define
\[
 a_h=\min_{|v|=1}\max_\theta|h(\theta)m_\theta\cdot v|.
\]
This minimum is positive. Otherwise, at a minimizing unit vector $v$,
$m_\theta\cdot v$ would vanish wherever $h(\theta)\ne0$. The latter set
contains an open interval, whereas the former trigonometric function
cannot vanish on an interval. Continuity and compactness justify taking
the positive minimum. Put $r_K=\min_{x\in K}|x|>0$. Then
\[
 \max_\theta|h(\theta)m_\theta\cdot x|\ge a_h r_K
 \qquad(x\in K).
\]
By \eqref{v4:geom-C-leading}, at a suitable direction in every fiber,
\[
 C_\varepsilon(x,re_\theta)+C_\varepsilon(x,-re_\theta)
 =\frac{4\varepsilon}{r}h(\theta)(m_\theta\cdot x)
        (m_\theta^{\mathsf T})^{\otimes3}+o(|\varepsilon|)/r
 \ne0
\]
for sufficiently small nonzero $\varepsilon$, uniformly in $x\in K$.
If a tangent norm were reversible, its energy and fiber Hessian would be
even in $y$, and its Cartan tensor, being half the derivative of that Hessian,
would be odd. This contradicts the preceding inequality. Each tangent norm
is therefore nonreversible. A Riemannian norm is reversible, so none of
these fibers is Riemannian.
\end{proof}

The leading tensor records the seed in the prescribed initial-value coordinates.
To distinguish completed metrics intrinsically, we now use the canonical
structure equations to define a scalar profile on each indicatrix.

\subsection{The Cartan structure equations}

On the unit tangent bundle $r=1$, define
\begin{equation}\label{v41:canonical-coframe}
 \omega^1=-\beta,\qquad
 \omega^2=\alpha,\qquad
 \omega^3=-\dd s.
\end{equation}
The corresponding three vector fields are
$-w^i\delta_i,u^i\delta_i,-\partial_s$.
They are componentwise dual to the three forms.
Here frames in two different spaces must be distinguished:
$e_1=-w,e_2=u$ is an orthonormal frame of the pullback bundle $\pi^*TU$,
and both vectors at $(x,u(x,s))$ belong to $T_xU$;
whereas $-w^i\delta_i,u^i\delta_i,-\partial_s$ belong to the tangent space of the unit tangent bundle,
with the first two being horizontal lifts and the third the vertical direction along the indicatrix.
This distinction identifies the structure equations with the preceding tensor computations in natural coordinates.
By \eqref{v4:geom-coframeODE},
\[
 \dd\alpha=\dd s\wedge\beta,\qquad
 \dd\beta=\dd s\wedge(-\alpha+I\beta),
\]
and therefore
\begin{align*}
 \dd\omega^1&=I\omega^1\wedge\omega^3
                      +\omega^2\wedge\omega^3,\\
 \dd\omega^2&=-\omega^1\wedge\omega^3,\qquad
 \dd\omega^3=0.
\end{align*}
Using the standard structure equations for surfaces,
\[
 \dd\omega^1=-\mathcal I\omega^1\wedge\omega^3
                      +\omega^2\wedge\omega^3,\quad
 \dd\omega^2=-\omega^1\wedge\omega^3,\quad
 \dd\omega^3=\mathcal K\omega^1\wedge\omega^2
                      -\mathcal J\omega^1\wedge\omega^3,
\]
we read off
$\mathcal I=-I,\mathcal J=0,\mathcal K=0$.
These are the surface invariants corresponding to the tensor formulas above.
The sign of $\mathcal I$ comes from choosing the transverse unit vector to be $-w$;
it does not affect the tensor formulas in natural coordinates above.

This frame also expresses the distinction between the Landsberg and non-Berwald properties through the two horizontal directions.
Subscripts $1,2$ denote differentiation along $-w^i\delta_i,u^i\delta_i$, respectively.
Since $\delta_i$ keeps $s$ fixed and $\dd_x I=H_{qq}\beta$, we have
\begin{equation}\label{v41:cartan-transverse-derivative}
 \mathcal I_2=-H_{qq}\beta(u)=0,\qquad
 \mathcal I_1= H_{qq}\beta(w)=H_{qq}.
\end{equation}
The first equality gives the zero value in the surface identity $\mathcal J=\mathcal I_2$;
the second gives precisely the transverse coefficient in \eqref{v4:geom-Bmain}:
$B_{ru}(e_1,e_1,e_1)=(\mathcal I_1/r)e_1$.

The remaining derivative $\mathcal I_3$ has a simple evolution along
geodesics, also observed in \cite[Remark~8.11]{Mao}. Write $E_j$ for
the fields dual to $\omega^j$. Here $E_3=-\partial_s$ and $\mathcal I=-H_q$,
so $\mathcal I_3=H_{qs}+H_{qq}(-p+H)$ along the oscillator. On a
unit geodesic with angle $s$ fixed, $p$ increases at speed one and
$q$ is constant. Thus, writing $t$ for geodesic time, we obtain
\begin{equation}\label{v44:geom-affine-geodesic}
 \mathcal I_3(t)=\mathcal I_3(0)-t\,\mathcal I_1(0).
\end{equation}
This identity also holds on every regular Landsberg surface of zero
flag curvature. Indeed, its structure equations give
$[E_2,E_1]=0$, $[E_2,E_3]=-E_1$ and $E_2\mathcal I=\mathcal J=0$.
Hence $E_2\mathcal I_1=0$ and $E_2\mathcal I_3=-\mathcal I_1$ along
the unit geodesic flow generated by $E_2$. The computation uses local
oriented frames and continues along the geodesic by choosing its
transverse orientation continuously.

Likewise, zero flag curvature does not require the mixed curvature of the Chern connection to vanish.
With the index convention $\nabla e_a=\theta_a{}^b\otimes e_b$,
the Chern connection forms of a surface \cite[Chapter~4]{BCS} in this frame are
\[
 \Theta=
 \begin{pmatrix}-\mathcal I\omega^3&-\omega^3\\\omega^3&0\end{pmatrix}
 =\begin{pmatrix}-I\,\dd s&\dd s\\-\dd s&0\end{pmatrix}.
\]
Every entry is a multiple of $\dd s$, so $\Theta\wedge\Theta=0$.
Using $\dd I=H_{qq}\beta+(\partial_s I)\dd s$ gives
\[
 \Omega=\dd\Theta-\Theta\wedge\Theta
 =\begin{pmatrix}-H_{qq}\,\omega^1\wedge\omega^3&0\\0&0\end{pmatrix}.
\]
Thus the horizontal--horizontal curvature vanishes, whereas the horizontal--vertical mixed component is nonzero wherever $H_{qq}\ne0$.

\subsection{An intrinsic angular profile and separation of seeds}
\label{v41:profile-section}

Write $SF=\{(x,y):F(x,y)=1\}$ and let $E_1$ be the vector field dual to
$\omega^1$ in the canonical coframe. Thus $E_1=-w^i\delta_i$ and
$\mathcal I_1=E_1\mathcal I$ is a smooth function on $SF$.
The reconstruction identifies $SF$ with $U\times\R/(2\pi\mathbb Z)$ by
$(x,s)\mapsto(x,u(x,s))$. Since $g_u(u_s,u_s)=1$, the parameter $s$
measures Hessian arclength on each indicatrix. The positive angular
density $|\omega^3|$ consequently has total mass $2\pi$ there.

The structure equations give $\dd\omega^3=0$ because
$\mathcal K=\mathcal J=0$. Thus the canonical angle form has local
primitives. In the present construction the global circle-valued
coordinate $s$ and its period $2\pi$ come from the full indicatrix
parametrization, and \eqref{v41:canonical-coframe} gives $\omega^3=-\dd s$.
Define
\begin{equation}\label{v41:profile-definition}
 \overline{\mathcal I_1}(x)=\frac1{2\pi}
          \int_{S_xF}\mathcal I_1\,|\omega^3|,
 \qquad
 \mathscr P_F(x,y)=\frac{\mathcal I_1(x,y)}{\overline{\mathcal I_1}(x)}.
\end{equation}
This is a smooth function on the part of $SF$ where the denominator is
nonzero. We write $\mathscr P_F(x,s)$ for its expression in the preceding
parametrization. A local homothety means a diffeomorphism
$\varphi:V\to V'$ between nonempty connected base open sets such that
$\varphi^*F'=\lambda F$ for a positive constant $\lambda$.
Its unit-bundle lift is
\[
 \Phi(x,y)=\bigl(\varphi(x),\lambda^{-1}d\varphi_x(y)\bigr).
\]
An isometry is the case $\lambda=1$.

\begin{proposition}[The canonical profile]\label{v41:profile-invariant}
Let $F,F'$ be reconstructed metrics as above. A local homothety has a
constant orientation sign $\delta\in\{1,-1\}$ and satisfies
\begin{equation}\label{v41:coframe-homothety}
 \Phi^*\omega'^1=\delta\lambda\omega^1,\qquad
 \Phi^*\omega'^2=\lambda\omega^2,\qquad
 \Phi^*\omega'^3=\delta\omega^3.
\end{equation}
Consequently, for one constant $c\in\R/(2\pi\mathbb Z)$,
\begin{equation}\label{v41:profile-homothety}
 s'\circ\Phi=\delta s+c,
 \qquad
 \mathcal I'_1\circ\Phi=\lambda^{-1}\mathcal I_1,
 \qquad
 \mathscr P_{F'}\circ\Phi=\mathscr P_F
\end{equation}
where the profiles are defined.

For the metrics in Theorem~\ref{v4:geometric-stability}, with
$\operatorname{avg}h=1$, the denominator is nonzero for sufficiently small
nonzero $\varepsilon$, and
\begin{equation}\label{v41:profile-asymptotic}
 \eta_h(\varepsilon):=
 \sup_{x\in U_{\varepsilon,h}}
 \|\mathscr P_{F_{\varepsilon,h}}(x,\cdot)-h\|_{C^0(S^1)}
 \longrightarrow0.
\end{equation}
\end{proposition}
\begin{proof}
Use the base orientations of the initial-value coordinates and put
$\delta=\operatorname{sgn}\det D\varphi$. The determinant is continuous
and nowhere zero on the connected set $V$, so $\delta$ is constant.
At a unit vector $y$, the vector $e_2=y$ and the transverse unit vector
$e_1$ form the Hessian-orthonormal frame with the orientation convention
used above. The second
semibasic form is the Hilbert form, and the first is the metric dual of
$e_1$. The homothety identity and degree-zero homogeneity of the Hessian give
\[
 g'_{\lambda^{-1}d\varphi(y)}(d\varphi(v),d\varphi(z))
       =\lambda^2g_y(v,z).
\]
Hence $e'_2\circ\Phi=\lambda^{-1}d\varphi(e_2)$ and
$e'_1\circ\Phi=\delta\lambda^{-1}d\varphi(e_1)$.
Pairing these vectors with the Hessians proves the first two formulas in
\eqref{v41:coframe-homothety}.

The third form annihilates the canonical horizontal distribution and restricts
to oriented unit angular length on each indicatrix. A constant multiple of
the energy has the same Euler--Lagrange spray. Thus a constant homothety
preserves its horizontal distribution. On vertical vectors the normalized
lift preserves the Hessian length and changes the orientation by $\delta$.
These two properties determine the third form and prove its transformation
in~\eqref{v41:coframe-homothety}. This also shows why the two orientation
signs in that formula are the same.

Since $\omega^3=-\dd s$ and $\omega'^3=-\dd s'$, the differential of
$s'\circ\Phi-\delta s$ is zero. The unit bundle over the connected set $V$
is connected, so this circle-valued difference is constant. Pulling back the
first structure equation gives $\mathcal I'\circ\Phi=\delta\mathcal I$.
The dual-frame transformation is
$\Phi_*E_1=\delta\lambda E'_1$. Differentiating the scalar identity along
$E_1$ therefore gives
\[
 \delta\lambda(\mathcal I'_1\circ\Phi)=\delta\mathcal I_1.
\]
The angular density is preserved, so its fiber average transforms by the
same factor $\lambda^{-1}$. Division proves the profile identity.

Finally, equations~\eqref{v41:cartan-transverse-derivative} and
\eqref{v4:geom-I-approximation} give, uniformly on the allowed base set and
the full angular circle,
\[
 \mathcal I_1=2\varepsilon h(s)+o(|\varepsilon|),\qquad
 \overline{\mathcal I_1}=2\varepsilon+o(|\varepsilon|).
\]
For small nonzero $\varepsilon$, the absolute value of the denominator is
at least $|\varepsilon|$. Dividing these estimates proves
\eqref{v41:profile-asymptotic}.
\end{proof}

For the realization in Theorem~\ref{v4:exact-closure}, we can prescribe
$k=2$, $a=3$ before choosing the finite approximation. Its first
correction term in \eqref{v4:formal-scaled-control} has order
$\varepsilon^3$. Equation~\eqref{v4:newton-closure-conclusion}
therefore gives $H_{\varepsilon,h}=\varepsilon hq^2+O_h(|\varepsilon|^3)$
in $C^2$ on the sampling cylinder. By
\eqref{v41:cartan-transverse-derivative}, uniformly in $x,s$,
\begin{equation}\label{v44:geom-quantitative-profile}
 \begin{aligned}
 \mathcal I_1(x,s)&=2\varepsilon h(s)+O_h(|\varepsilon|^3),&
 \overline{\mathcal I_1}(x)&=2\varepsilon+O_h(|\varepsilon|^3),\\
 \eta_h(\varepsilon)&=O_h(\varepsilon^2).
 \end{aligned}
\end{equation}
The last estimate follows by division. In the fixed natural directions
of \eqref{v4:geom-B-leading}, the angular shift and the changes of
frame in \eqref{v4:geom-harmonic-estimates} are $O_h(|\varepsilon|)$.
Multiplication by the leading coefficient $2\varepsilon h$ gives
an $O_h(\varepsilon^2)/r$ remainder for that tensor.
The constants here may depend on the chosen base set and finite preparation.

The orbit of $\mathscr P_F(x,\cdot)$ under angular translation and reflection
is thus an exact invariant of a completed metric up to local homothety.
The prescribed seed is its small-amplitude limit. Define the corresponding
orbit distance between two seeds by
\begin{equation}\label{v41:seed-orbit-distance}
 d_{O(2)}(h,k)=\min_{\delta=\pm1,\ c\in S^1}
       \|h-k(\delta\,\cdot+c)\|_{C^0(S^1)}.
\end{equation}
The minimum exists by compactness. It is positive precisely when the two
seeds are not related by a translation or reflection.

\begin{proposition}[Separation of fixed seeds]\label{v41:seed-separation}
Let $h,k$ be mean-one seeds with $d_{O(2)}(h,k)>0$, and let their completed
metrics satisfy the hypotheses of Theorem~\ref{v4:geometric-stability}.
There are positive amplitude thresholds for the two families such that,
below those thresholds, no nonempty base open subset of one metric is
homothetic to a base open subset of the other. In particular the two
metrics have no isometric germs, with all nonzero fiber directions included.
\end{proposition}
\begin{proof}
Suppose such a homothety exists and restrict to connected open sets.
For any point $x$ of its source, Proposition~\ref{v41:profile-invariant}
and the triangle inequality give
\[
 d_{O(2)}(h,k)
 \le \|h-\mathscr P_{F_{\varepsilon,h}}(x,\cdot)\|_{C^0}
    +\|\mathscr P_{F_{\varepsilon',k}}(\varphi(x),\cdot)-k\|_{C^0}
 \le\eta_h(\varepsilon)+\eta_k(\varepsilon').
\]
Choose the two amplitude thresholds so that each error is less than one
third of the positive orbit distance. Equation~\eqref{v41:profile-asymptotic}
then yields a contradiction. Uniformity in the base point excludes
homotheties between arbitrary nonempty base open subsets.
\end{proof}

The thresholds in this proposition depend on the chosen pair of seeds.
Thus the conclusion separates fixed distinct seed orbits at sufficiently
small amplitudes; it does not identify a continuous seed parameter space
with a moduli space at one fixed amplitude. Section~\ref{v4:realization-chapter}
uses a uniformly separated sequence of seed orbits to produce infinitely
many pairwise nonhomothetic examples.

The remaining parts construct periodic controls with exact return on an
open set of nonzero initial values and with the smallness required by
Theorem~\ref{v4:geometric-stability}. The preceding geometric statements
then apply to the completed controls.
\part{Seeds and finite-order preparation}\label{part:seeds}
\section{Balanced seeds and their examples}\label{v4:seed-section}

We first construct the periodic functions from which the metrics will be obtained.
The balance equations arise from the first nonzero return coefficient of the oscillator.
A positive coefficient associated with the seed then makes the higher-order moment
equations solvable. These two requirements produce a large class of formal seeds.
Exact return on a neighborhood of a nonzero initial value uses the further
characteristic condition in Section~\ref{v4:char-section}.

\subsection{The oscillator, rotating coordinates, and homogeneous expansion}
\label{v4:seed-sec-expansion}

Let $s\in\R/(2\pi\mathbb Z)$ and consider
\begin{equation}\label{v4:seed-oscillator}
 p'=q,\qquad q'=-p+H(q,s),
 \qquad H(q,s)=h(s)q^2+O(q^4).
\end{equation}
Here $h$ is a real-valued smooth function satisfying $h(s+\pi)=h(s)$.
The remainder is smooth near $q=0$, and its derivatives of orders zero through three with respect to $q$ all vanish there.
The cubic return formula in this section requires no additional symmetry of the remainder;
when the half-period map is used, the full $H$ will explicitly be required to be $\pi$-periodic.

Set
\begin{equation}\label{v4:seed-rotation-frame}
 C_s=\begin{pmatrix}\cos s&\sin s\\-\sin s&\cos s\end{pmatrix},
 \qquad e_s=\binom{\cos s}{\sin s},
 \qquad m_s=\binom{-\sin s}{\cos s}.
\end{equation}
The two columns of $C_s^{-1}$ are therefore $e_s,m_s$.
For initial value $x=(x^1,x^2)^{\mathsf T}$, denote the flow in the original coordinates by $\Psi(x,s)$
and define the flow in rotating coordinates by
\[
 \Phi(x,s)=C_s^{-1}\Psi(x,s),\qquad \Phi(x,0)=x.
\]
Direct differentiation cancels the linear oscillator part completely and gives
\begin{equation}\label{v4:seed-rotated-equation}
 \partial_s\Phi(x,s)
 =m_s H(m_s\cdot\Phi(x,s),s).
\end{equation}
Introduce the linear functions
\begin{equation}\label{v4:seed-linear-coordinates}
 \ell_s(x)=m_s\cdot x=-x^1\sin s+x^2\cos s,
 \qquad r_s(x)=e_s\cdot x=x^1\cos s+x^2\sin s.
\end{equation}
They satisfy $\ell_s'=-r_s$ and $r_s'=\ell_s$.
These linear functions of the initial value are the linear parts of the solutions $q(x,s),p(x,s)$.

The homogeneous quadratic vector field is
\begin{equation}\label{v4:seed-quadratic-field}
 X(s,x)=h(s)m_s\ell_s(x)^2.
\end{equation}
Since $H$ has no cubic term, the rotating flow can be written as
\[
 \Phi(x,s)=x+U_2(x,s)+U_3(x,s)+O(|x|^4),
 \qquad U_j(x,0)=0,
\]
where $U_j$ is a vector polynomial in $x$, homogeneous of degree $j$.
Substituting into \eqref{v4:seed-rotated-equation} and comparing degrees gives
\begin{equation}\label{v4:seed-homogeneous-odes}
 \partial_sU_2=X(s,x),\qquad
 \partial_sU_3=D_xX(s,x)U_2(x,s).
\end{equation}
The second equation does not involve $U_3$ itself: the derivative of the quadratic field is linear,
so its product with $U_3$ already has degree four.
Likewise, $X(s,x+U_2)-X(s,x)-DX(s,x)U_2$ has degree at least four.

Since $C_{2\pi}=\Id$, the full return $P(x)=\Psi(x,2\pi)$ equals $\Phi(x,2\pi)$.
Its homogeneous quadratic and cubic terms are therefore
\begin{align}
 B_2(x)&=\int_0^{2\pi}X(s,x)\,\dd s,
 \label{v4:seed-B2-integral}\\
 B_3(x)&=\int_0^{2\pi}D_xX(s,x)
        \left(\int_0^sX(t,x)\,\dd t\right)\dd s.
 \label{v4:seed-B3-integral}
\end{align}
The flow equation imposes the time ordering $0\le t\le s\le2\pi$ in the cubic integral.

\begin{lemma}[The quadratic return vanishes automatically]\label{v4:seed-lemma-B2}
 For every $\pi$-periodic $h$, one has $B_2=0$.
\end{lemma}
\begin{proof}
 We have $m_{s+\pi}=-m_s$ and $\ell_{s+\pi}=-\ell_s$.
 Together with $h(s+\pi)=h(s)$, this gives
 $X(s+\pi,x)=-X(s,x)$.
 Splitting \eqref{v4:seed-B2-integral} into two half-periods gives two cancelling integrals.
\end{proof}

\subsection{The three basic time modes in complex coordinates}
\label{v4:seed-sec-complex-modes}

Identify the real vector $(v^1,v^2)$ with the complex number $v^1+iv^2$, and set
\[
 z=x^1+ix^2,\qquad Jx=(-x^2,x^1),
\]
Then $Jx$ corresponds to $iz$.
Here $J$ is counterclockwise rotation, opposite to the rotation of the original linear oscillator.
From \eqref{v4:seed-rotation-frame} we obtain
\[
 m_s^1+im_s^2=ie^{is},\qquad
 \ell_s=\frac{i}{2}\left(e^{is}\bar z-e^{-is}z\right).
\]
The complex representation of the quadratic field is therefore
\begin{equation}\label{v4:seed-complex-field}
 X^{\C}(s,z)
 =h(s)\left(
 -\frac{i}{4}e^{3is}\bar z^2
 +\frac{i}{2}e^{is}z\bar z
 -\frac{i}{4}e^{-is}z^2\right).
\end{equation}
Thus, although $h$ contains only even time frequencies, the rotated quadratic field contains only odd time frequencies.
The absence of the zero frequency is the Fourier counterpart of the vanishing quadratic return.

We adopt the convention
\begin{equation}\label{v4:seed-fourier-convention}
 h(s)=\sum_{k\in2\mathbb Z}\widehat h_k e^{iks},
 \qquad
 \widehat h_k=\frac1{2\pi}\int_0^{2\pi}h(s)e^{-iks}\,\dd s,
 \qquad \widehat h_{-k}=\overline{\widehat h_k}.
\end{equation}
First assume that $h$ is a trigonometric polynomial, so all computations involve finite sums.
The general smooth case will follow by continuity at the end of this section.
Write $X^{\C}$ as
\[
 X^{\C}(s,z)=\sum_{\nu\in2\mathbb Z+1}X_\nu^{\C}(z)e^{i\nu s}.
\]
By \eqref{v4:seed-complex-field}, each coefficient is explicitly
\begin{equation}\label{v4:seed-X-frequency}
 X_\nu^{\C}(z)
 =-\frac{i}{4}\widehat h_{\nu-3}\bar z^2
  +\frac{i}{2}\widehat h_{\nu-1}z\bar z
  -\frac{i}{4}\widehat h_{\nu+1}z^2.
\end{equation}
All absent Fourier coefficients are understood to be zero.
To compute the actual two-dimensional derivative, we must also retain the complex conjugate of $X^{\C}$:
\begin{equation}\label{v4:seed-conjugate-field}
 \overline{X^{\C}(s,z)}
 =h(s)\left(
 \frac{i}{4}e^{-3is}z^2
 -\frac{i}{2}e^{-is}z\bar z
 +\frac{i}{4}e^{is}\bar z^2\right).
\end{equation}
In complex notation, the real differential is
\begin{equation}\label{v4:seed-real-derivative-complex}
 DX^{\C}[V]
 =\partial_zX^{\C}\,V^{\C}
  +\partial_{\bar z}X^{\C}\,\overline{V^{\C}}.
\end{equation}
Both terms are needed because $X^{\C}$ depends on $z$ and $\bar z$.

\subsection{How time integration becomes frequency pairing}
\label{v4:seed-sec-frequency-pairing}

Returning to real vector notation, write $X(s,x)=\sum_{\nu\ne0}X_\nu(x)e^{i\nu s}$.
Since there is no zero frequency, the first homogeneous equation integrates directly to
\begin{equation}\label{v4:seed-U2-frequency}
 U_2(x,s)
 =\sum_{\mu\ne0}\frac{e^{i\mu s}-1}{i\mu}X_\mu(x).
\end{equation}
Substitution into the cubic return integral gives
\[
 B_3
 =\sum_{\nu,\mu\ne0}\frac{DX_\nu X_\mu}{i\mu}
 \left(\int_0^{2\pi}e^{i(\nu+\mu)s}\,\dd s
       -\int_0^{2\pi}e^{i\nu s}\,\dd s\right).
\]
The second integral always vanishes; only when $\mu=-\nu$ does the first equal $2\pi$.
Hence
\begin{equation}\label{v4:seed-B3-frequency}
 B_3(x)=2\pi\sum_{\nu\ne0}
            \frac{DX_\nu(x)X_{-\nu}(x)}{-i\nu}.
\end{equation}
The minus sign in the denominator comes from the inner frequency $\mu=-\nu$.
Formula \eqref{v4:seed-B3-frequency} preserves the time ordering of the double integral
and reduces the entire integration problem to algebraic operations on finitely many cubic polynomials.

We next expand $X_\nu$ in terms of the Fourier coefficients of the seed.
Let $k$ be the seed frequency carried by the outer field and $l$ that carried by the inner field.
By \eqref{v4:seed-complex-field}, \eqref{v4:seed-conjugate-field}, and
\eqref{v4:seed-real-derivative-complex}, only four cubic monomials can occur.
Frequency conservation respectively requires
\begin{equation}\label{v4:seed-selection-rule}
 \begin{array}{c|cccc}
 \text{monomial}&\bar z^3&z\bar z^2&z^2\bar z&z^3\\ \hline
 k+l&-4&-2&0&2
 \end{array}.
\end{equation}
For each surviving pair, the inner time frequency is the negative of the
outer one. Formula~\eqref{v4:seed-B3-frequency} therefore determines its
integration factor. We record the resulting weights and then list all
their individual sources.

\begin{lemma}[Pairing weights for the four monomials]\label{v4:seed-lemma-weight-table}
 Write the coefficient of each monomial in $B_3^{\C}/(2\pi i)$ as
 $\sum_k\widehat h_k\widehat h_{t-k}R_t(k)$.
 Then the original weights and their symmetrizations are as follows:
\begin{equation}\label{v4:seed-weight-table}
\begin{array}{c|c|c}
 \text{monomial, }t&R_t(k)&\frac12\bigl(R_t(k)+R_t(t-k)\bigr)\\[2pt]\hline
 \bar z^3,\ -4
 &\frac1{8(k+3)}+\frac1{8(k+1)}&0\\[4pt]
 z\bar z^2,\ -2
 &-\frac1{4(k+3)}-\frac3{8(k+1)}-\frac1{8(k-1)}
 &\frac1{4(k-1)(k+3)}\\[4pt]
 z^2\bar z,\ 0
 &\frac1{8(k+3)}+\frac3{8(k+1)}+\frac1{4(k-1)}
 &\frac{3-k^2}{2(k^2-1)(k^2-9)}\\[4pt]
 z^3,\ 2
 &-\frac1{8(k+1)}-\frac1{8(k-1)}
 &\frac1{4(k+1)(k-3)}
\end{array}
\end{equation}
 All denominators are nonzero because $k$ is even.
\end{lemma}
\begin{proof}
 Let
\[
 \mathsf A=-\frac{i}{4}\bar z^2,\qquad
 \mathsf B=\frac{i}{2}z\bar z,\qquad
 \mathsf C=-\frac{i}{4}z^2
\]
 be the polynomial parts of the three modes in
 \eqref{v4:seed-complex-field}, with frequency shifts $3,1,-1$,
 respectively. In the middle column below the differentiated factor
 is the outer mode; a bar on the inner factor means complex conjugation.
 The seed frequencies satisfy the selection rule
 \eqref{v4:seed-selection-rule}. These are all twelve nonzero terms
 produced by the real differential \eqref{v4:seed-real-derivative-complex}.
\begin{equation}\label{v43:seed-weight-sources}
\begin{array}{c|c|c}
 \text{target}&\text{outer derivative and inner mode}&
       \text{weight in }B_3^{\C}/(2\pi i)\\ \hline
 \bar z^3&(\partial_{\bar z}\mathsf A)\overline{\mathsf C}
       &1/[8(k+3)]\\
 &(\partial_z\mathsf B)\mathsf A&1/[8(k+1)]\\[2pt]\hline
 z\bar z^2&(\partial_{\bar z}\mathsf A)\overline{\mathsf B}
       &-1/[4(k+3)]\\
 &(\partial_z\mathsf B)\mathsf B&-1/[4(k+1)]\\
 &(\partial_{\bar z}\mathsf B)\overline{\mathsf C}&-1/[8(k+1)]\\
 &(\partial_z\mathsf C)\mathsf A&-1/[8(k-1)]\\[2pt]\hline
 z^2\bar z&(\partial_{\bar z}\mathsf A)\overline{\mathsf A}
       &1/[8(k+3)]\\
 &(\partial_{\bar z}\mathsf B)\overline{\mathsf B}&1/[4(k+1)]\\
 &(\partial_z\mathsf B)\mathsf C&1/[8(k+1)]\\
 &(\partial_z\mathsf C)\mathsf B&1/[4(k-1)]\\[2pt]\hline
 z^3&(\partial_{\bar z}\mathsf B)\overline{\mathsf A}&-1/[8(k+1)]\\
 &(\partial_z\mathsf C)\mathsf C&-1/[8(k-1)]
\end{array}
\end{equation}
 For example, $(\partial_{\bar z}\mathsf B)\overline{\mathsf A}$
 has coefficient $(i/2)(i/4)=-1/8$. Its outer time frequency is
 $k+1$ and its inner time frequency is $-k-1$, giving the integrated
 coefficient $-i/[8(k+1)]$. Division by $i$ gives the displayed
 weight. Summing the entries for each target yields the middle
 column of \eqref{v4:seed-weight-table}.

 Since $\widehat h_k\widehat h_{t-k}$ is invariant under $k\leftrightarrow t-k$,
 any $R_t$ can be replaced by its symmetrization without changing the sum.
 For example, when $t=2$,
\[
 \frac12\bigl(R_2(k)+R_2(2-k)\bigr)
 =\frac1{16}\left(\frac1{k-3}-\frac1{k+1}\right)
 =\frac1{4(k+1)(k-3)}.
\]
 For $t=-4$, the interchange replaces $k+3,k+1$ by $-k-1,-k-3$, respectively,
 and the two terms cancel exactly.
 For $t=0$, symmetrization gives
\[
 -\frac{3}{8(k^2-9)}-\frac1{8(k^2-1)}
 =\frac{3-k^2}{2(k^2-1)(k^2-9)}.
\]
 The remaining case $t=-2$ follows by the same common-denominator calculation, or from the expression for $t=2$ by setting $k\mapsto-k$.
\end{proof}

\subsection{A closed formula for the general cubic return}
\label{v4:seed-sec-cubic-formula}

Define the real even multiplier
\begin{equation}\label{v4:seed-multiplier}
 \mathfrak m(k)=\frac{3-k^2}{(k^2-1)(k^2-9)},\qquad k\in2\mathbb Z,
\end{equation}
and the real and complex quadratic forms
\begin{align}
 Q(h)&=\sum_{k\in2\mathbb Z}\mathfrak m(k)|\widehat h_k|^2,
 \label{v4:seed-Q}\\
 S(h)&=\sum_{k\in2\mathbb Z}
       \frac{\widehat h_k\widehat h_{2-k}}{(k+1)(k-3)}.
 \label{v4:seed-S}
\end{align}
We use $\mathfrak m$ to distinguish the scalar Fourier multiplier from the vector $m_s$.
The condition $S=0$ is equivalent to the two real conditions $\operatorname{Re}S=\operatorname{Im}S=0$.

\begin{theorem}[The cubic return of a general half-period seed]\label{v4:seed-theorem-cubic}
 For any real smooth $\pi$-periodic seed $h$ in \eqref{v4:seed-oscillator},
 the full return satisfies $B_2=0$ and
\begin{equation}\label{v4:seed-cubic-final}
 B_3(x)=\pi\left[Q(h)|x|^2+
        \operatorname{Re}\bigl(S(h)z^2\bigr)\right]Jx.
\end{equation}
 In real coordinates, the same cubic term is
\begin{equation}\label{v4:seed-cubic-real}
 B_3(x)=\pi\bigl[(Q+\Re S)(x^1)^2-2(\Im S)x^1x^2
                   +(Q-\Re S)(x^2)^2\bigr]Jx .
\end{equation}
 Thus $B_3=0$ if and only if $Q(h)=S(h)=0$.
\end{theorem}
\begin{proof}
 First apply Lemma \ref{v4:seed-lemma-weight-table} to trigonometric polynomials.
 The coefficient of $z^2\bar z$ is $\pi iQ(h)$.
 The coefficient of $z^3$ is $\pi iS(h)/2$.
 For $z\bar z^2$, changing the summation index in the table by $k\mapsto-k$ and using real-valuedness gives
\[
 \sum_k\frac{\widehat h_k\widehat h_{-2-k}}{(k-1)(k+3)}
 =\overline{S(h)}.
\]
 Its coefficient is therefore $\pi i\overline{S(h)}/2$.
 The coefficient of $\bar z^3$ vanishes.
 The four terms combine into
\[
 B_3^{\C}
 =\pi i\left(Qz^2\bar z+\frac{S}{2}z^3
                         +\frac{\bar S}{2}z\bar z^2\right)
 =\pi iz\left(Q|z|^2+\operatorname{Re}(Sz^2)\right),
\]
 which is precisely \eqref{v4:seed-cubic-final}.

 For general smooth $h$, take the real symmetric Fourier truncations $h_N$.
 Since $h_N\to h$ in every fixed $C^r$ norm, the coefficients of the double integral
 \eqref{v4:seed-B3-integral} converge.
 On the other hand, the series defining $Q,S$ converge absolutely and are continuous on $L^2$.
 For $S$, this follows from the discrete Cauchy--Schwarz inequality and
 $\sup_{k\in2\mathbb Z}|(k+1)^{-1}(k-3)^{-1}|<\infty$.
 Letting $N\to\infty$ preserves the formula.

 Finally, $|x|^2,\operatorname{Re}z^2,\operatorname{Im}z^2$
 form a basis of the real homogeneous quadratic polynomials.
 If the quadratic polynomial in braces vanishes identically, its three real coefficients vanish individually; the converse is immediate.
\end{proof}

Formula \eqref{v4:seed-cubic-final} shows that the cubic error is tangent to each circle centered at the origin.

\subsection{Why cubic controls cannot remove these obstructions}
\label{v4:seed-sec-cubic-control}

Although we choose $H=hq^2+O(q^4)$, it is necessary to check that this choice does not artificially create the obstructions.
Consider the broader expansion
\[
 H(q,s)=h(s)q^2+g(s)q^3+O(q^4).
\]
The additional contribution of the cubic control to the cubic return is
\begin{equation}\label{v4:seed-direct-cubic-gradient}
 \int_0^{2\pi}g(s)m_s\ell_s(x)^3\,\dd s
 =\nabla_x\left(\frac14\int_0^{2\pi}g(s)\ell_s(x)^4\,\dd s\right).
\end{equation}
It is a gradient and therefore has zero curl.
For any real homogeneous quadratic polynomial $A(x)$, direct calculation gives
\begin{equation}\label{v4:seed-curl-tangential}
 \Curl(AJx)
 =\partial_1(Ax^1)-\partial_2(-Ax^2)
 =2A+x\cdot\nabla A=4A.
\end{equation}
Here we use the convention $\Curl V=\partial_1V^2-\partial_2V^1$.

\begin{corollary}[Balance at an equilibrium and on a fixed return domain]
\label{v4:seed-cor-necessary}\label{v44:balance-necessary}
Let $h$ be a real smooth $\pi$-periodic function. Either of the following
conditions implies $Q(h)=S(h)=0$.
\begin{enumerate}[label=\textup{(\roman*)}]
\item The cubic term of the return at the equilibrium vanishes for a
control $H(q,s)=h(s)q^2+g(s)q^3+O(q^4)$.
\item For a sequence of nonzero amplitudes $\varepsilon\to0$, the smooth
controls $H_\varepsilon=\varepsilon h(s)q^2+\rho_\varepsilon(q,s)$
have exact return on the same nonempty open set $U$. They are defined on
a fixed cylinder $(-b,b)\times\R/(2\pi\mathbb Z)$ with $U\subset\{|x|<b\}$,
and $\|\rho_\varepsilon\|_{C^2}=o(|\varepsilon|)$ on compact subcylinders.
\end{enumerate}
\end{corollary}
\begin{proof}
In (i), taking the curl removes \eqref{v4:seed-direct-cubic-gradient}.
Equations~\eqref{v4:seed-cubic-final} and \eqref{v4:seed-curl-tangential} give
\[
 \Curl B_3=4\pi\bigl[Q(h)|x|^2+\operatorname{Re}(S(h)z^2)\bigr].
\]
The three coefficients of this quadratic polynomial determine $Q$ and $S$.

For (ii), fix a closed ball $K$ of positive radius in $U$. Write
$X_s(x)=h(s)m_s\ell_s(x)^2$ and $R_{\varepsilon,s}(x)=m_s\rho_\varepsilon(\ell_s(x),s)$.
On a fixed neighborhood of $K$, the rotating vector field
$Y_s=\varepsilon X_s+R_{\varepsilon,s}$ has $C^2$ norm $O(|\varepsilon|)$,
uniformly in $s$. The integral equation on $[0,2\pi]$ then gives
$\Phi_s(x)-x=O(|\varepsilon|)$ there for $x\in K$. Taylor expansion of
$Y_s(\Phi_s(x))$ around $x$ and a second substitution in the integral
equation give, uniformly on $K$,
\[
 P_{H_\varepsilon}(x)-x
 =\int_0^{2\pi}Y_s(x)\,\dd s+
   \int_0^{2\pi}D Y_s(x)\int_0^s Y_t(x)\,\dd t\,\dd s
   +O(|\varepsilon|^3).
\]
The integral of $X_s$ is zero. In the double integral the terms containing
$R_{\varepsilon,s}$ or its derivative are $o(|\varepsilon|^2)$, since
$\|R_{\varepsilon,s}\|_{C^1}=o(|\varepsilon|)$. The remaining term is
$\varepsilon^2 B_3$ by \eqref{v4:seed-B3-integral}. The single integral
of $R_{\varepsilon,s}$ is the gradient of
\[
 G_\varepsilon(x)=\int_0^{2\pi}\int_0^{\ell_s(x)}
                       \rho_\varepsilon(q,s)\,\dd q\,\dd s .
\]
Take the distributional curl on the interior of $K$, divide by
$\varepsilon^2$, and pass to the limit. The left side is zero by exact
return, the gradient has zero curl, and the uniform $o(|\varepsilon|^2)$
error tends to zero after pairing its rescaled curl with each test
function. Thus $\Curl B_3$ vanishes on an open ball. Being a polynomial,
it vanishes everywhere, which proves (ii).
\end{proof}

\subsection{Half-periodicity controls only the parity of the first nonzero term}
\label{v4:seed-sec-leading-parity}

We now assume additionally that the full $H(q,s)$ is $\pi$-periodic and that
$H(0,s)=H_q(0,s)=0$.
Let $f$ be the flow map from $0$ to $\pi$.
By time periodicity, the full return is $P=f\circ f$, and
$Df(0)=C_\pi=-\Id$.
Thus $f$ and $P$ commute exactly.

\begin{lemma}[The first error has odd degree]\label{v4:seed-lemma-leading-odd}\label{v4:formal-odd-leading}
 If $P(x)=x+B_j(x)+O(|x|^{j+1})$, where $j\ge2$ and $B_j$ is nonzero,
 then $j$ is odd.
\end{lemma}
\begin{proof}
 From $f(P(x))=P(f(x))$, comparing homogeneous terms of degree $j$ gives
 $-B_j(x)=B_j(-x)$.
 Since $B_j$ is homogeneous of degree $j$, we also have $B_j(-x)=(-1)^jB_j(x)$.
 Nonvanishing forces $(-1)^j=-1$.
\end{proof}

\begin{corollary}[The low-order starting point for the formal recursion]\label{v4:seed-cor-order-five}
 For a full $\pi$-periodic control $H=hq^2+O(q^4)$, if $Q(h)=S(h)=0$, then
\[
 P(x)=x+O(|x|^5).
\]
\end{corollary}
\begin{proof}
 Theorem \ref{v4:seed-theorem-cubic} eliminates the quadratic and cubic terms.
 A nonzero quartic term would be the first nonzero term, contradicting the preceding lemma.
\end{proof}

\begin{remark}[The leading-term statement does not imply that all even-degree terms vanish]
 Write $f(x)=-x+A_2(x)+O(|x|^3)$ and
 $P(x)=x+B_3(x)+B_4(x)+O(|x|^5)$.
 Expanding the commutation relation one order further gives
\[
 2B_4=DA_2\,B_3-DB_3\,A_2.
\]
 Thus, when $B_3$ has not yet vanished, $B_4$ can be nonzero.
 What we use in the recursion is that, once all lower-order terms have been eliminated at each step,
 the next possible leading term has odd degree.
\end{remark}

\subsection{The coefficient operators and the balanced seed class}
\label{v4:seed-chapter-family}
\label{v4:seed-sec-cd}

On $\pi$-periodic functions, the operator
$(\partial_s^2+1)(\partial_s^2+9)$ has no kernel.
Indeed, a zero eigenvalue can occur only at frequencies $\pm1,\pm3$,
which do not belong to the $\pi$-periodic function space.
There is therefore a unique $\pi$-periodic solution $f$ of
\begin{equation}\label{v4:seed-f-equation}
 f''''+10f''+9f=9h.
\end{equation}
Define
\begin{equation}\label{v4:seed-cd-definition}
 c=\frac{f''+3f}{9},\qquad d=-\frac{2f'}9.
\end{equation}
The Fourier coefficients are explicitly
\begin{equation}\label{v4:seed-cd-fourier}
 \widehat c_k=\mathfrak m(k)\widehat h_k,
 \qquad
 \widehat d_k=\frac{-2ik}{(k^2-1)(k^2-9)}\widehat h_k.
\end{equation}
Thus $c=Lh$, where $L$ has the multiplier \eqref{v4:seed-multiplier}.
The operator $L$ is a real self-adjoint Fourier multiplier of order $-2$; the map $h\mapsto d$ has order $-3$.
In particular, smooth or real analytic $h$ produces $c,d$ in the same regularity class.

\begin{lemma}[Basic identities for the auxiliary coefficients]\label{v4:seed-lemma-cd-identity}
 For each fixed initial value $x$, we have
\begin{equation}\label{v4:seed-cd-identity}
 (\partial_s^2+4)\bigl(c(s)\ell_s(x)-d(s)r_s(x)\bigr)
 =h(s)\ell_s(x).
\end{equation}
 Moreover,
\begin{equation}\label{v4:seed-mean-and-Q}
 \operatorname{avg}c=\frac13\operatorname{avg}h,
 \qquad Q(h)=\operatorname{avg}(h c),
 \qquad
 \operatorname{avg}u:=\frac1{2\pi}\int_0^{2\pi}u(s)\,\dd s.
\end{equation}
\end{lemma}
\begin{proof}
 Using $\ell'=-r,r'=\ell$ and differentiating term by term gives
\[
 (\partial_s^2+4)(c\ell-dr)
 =(c''-2d'+3c)\ell-(2c'+d''+3d)r.
\]
 Substituting \eqref{v4:seed-cd-definition}, the first parenthesis becomes
 $(f''''+10f''+9f)/9=h$, and the second vanishes identically.
 This proves \eqref{v4:seed-cd-identity}.
 The mean formula follows from $\widehat c_0=\widehat h_0/3$.
 The last formula follows from Parseval's identity and $\widehat c_k=\mathfrak m(k)\widehat h_k$.
\end{proof}

The signs of the multiplier are crucial:
\begin{equation}\label{v4:seed-multiplier-sign}
 \mathfrak m(0)=\frac13,\qquad
 \mathfrak m(2)=\frac1{15},\qquad
 \mathfrak m(k)<0\quad (|k|\ge4).
\end{equation}
Its partial-fraction decomposition is
\begin{equation}\label{v4:seed-partial-fractions}
 \mathfrak m(k)
 =-\frac1{4(k^2-1)}-\frac3{4(k^2-9)}.
\end{equation}
This explains why high-frequency directions can cancel the positive quadratic form in the constant direction
and also provides a direct route to the exact summations below.

\begin{definition}[Balanced seeds]\label{v4:seed-class-A}
Let $L$ denote the Fourier multiplier $\widehat{Lh}_k=\mathfrak m(k)\widehat h_k$.
We set
\[
 \mathcal A^\infty=\{h\in C^\infty(\R/\pi\mathbb Z;\R):
 \operatorname{avg}h=1,\ Q(h)=S(h)=0,\ \min_s Lh(s)>0\}.
\]
The trigonometric polynomial members form $\mathcal A^{\mathrm{trig}}$.
The conditions are imposed on the quadratic seed, not on the final correcting control.
The smaller class used for exact return away from the origin is defined in
Definition~\ref{v4:seed-class-B} by an additional characteristic condition.
\end{definition}

For each member $c=Lh$ has a positive minimum. Positivity persists under sufficiently
small $L^2$ perturbations, since $L:L^2\to H^2\subset C^0$ is continuous.
The perturbations must also satisfy the three real balance equations.
Their transversality on the normalized hyperplane is proved in
Proposition~\ref{v4:formal-transversality}.

\subsection{Normalization and elementary transformations}
\begin{proposition}[Mean normalization and signs]\label{v4:seed-prop-normalization}
 If a real $\pi$-periodic $h$ satisfies $Q=S=0$ and $Lh>0$, then
 $a:=\operatorname{avg}h>0$ and $h/a\in\mathcal A^\infty$
 (with the corresponding regularity class understood).
 Moreover, every member of $\mathcal A^\infty$ must change sign.
\end{proposition}
\begin{proof}
 By \eqref{v4:seed-mean-and-Q}, $a=3\operatorname{avg}(Lh)>0$.
 Both $Q,S$ are homogeneous quadratic and $L$ is linear, so division by $a$ preserves balance and strict positivity.
 If normalized $h$ were everywhere nonnegative, its mean $1$ and the condition $c>0$
 would force $\operatorname{avg}(hc)>0$, contradicting $Q=0$.
 Nor can $h$ be everywhere nonpositive, since its mean is $1$.
 It must therefore take both positive and negative values.
\end{proof}

\begin{proposition}[The actions of translation, reflection, and scaling]\label{v4:seed-prop-transformations}
 Set $h_\tau(s)=h(s+\tau)$ and $h^{\mathrm{rev}}(s)=h(-s)$.
 Then
\begin{align}
 Q(h_\tau)&=Q(h),&S(h_\tau)&=e^{2i\tau}S(h),&Lh_\tau&=(Lh)(\cdot+\tau),
 \label{v4:seed-phase-action}\\
 Q(h^{\mathrm{rev}})&=Q(h),&S(h^{\mathrm{rev}})&=\overline{S(h)},
 &Lh^{\mathrm{rev}}&=(Lh)(-\cdot).
 \label{v4:seed-reflection-action}
\end{align}
 For any real constant $a$, one has $Q(ah)=a^2Q(h)$, $S(ah)=a^2S(h)$, and $L(ah)=aLh$.
\end{proposition}
\begin{proof}
 Time translation multiplies $\widehat h_k$ by $e^{ik\tau}$.
 The phases cancel in every term of $Q$, while the two indices in $S$ sum to $2$, leaving $e^{2i\tau}$.
 Time reflection replaces $\widehat h_k$ by $\widehat h_{-k}=\overline{\widehat h_k}$.
 The remaining assertions follow directly from the definitions.
\end{proof}

The slice of mean $1$ is not itself preserved by general scaling.
In the unnormalized family, positive scaling preserves $c>0$, while negative scaling sends it to the other branch, with $c<0$.
The coordinate relation at the oscillator level is also explicit: for $a\ne0$, set
\[
 H_a(q,s)=a^{-1}H(aq,s),\qquad
 \Psi_a(x,s)=a^{-1}\Psi(ax,s).
\]
Then $\Psi_a$ satisfies the oscillator equation of the same form with control $H_a$,
the quadratic seed becomes $ah$, and the return map satisfies $P_a(x)=a^{-1}P(ax)$.
Thus mean normalization is compatible with scaling of the oscillator coordinates.

\subsection{Construction on general nonadjacent even-frequency supports}
\label{v4:seed-sec-sparse-support}

A choice of Fourier support can make the complex constraint automatic.
The remaining balance condition is then one ellipsoid equation, and a single estimate proves positivity for all its members.
Set
\begin{equation}\label{v4:seed-alpha-definition}
 \alpha_k=-\frac12\mathfrak m(k)
 =\frac{k^2-3}{2(k^2-1)(k^2-9)}>0,
 \qquad k=4,6,8,\ldots.
\end{equation}
Let $K\subset\{4,6,8,\ldots\}$, and assume that no two distinct elements differ by $2$.
Consider
\begin{equation}\label{v4:seed-sparse-series}
 h(s)=1+\sum_{k\in K}(a_k\cos ks+b_k\sin ks),
 \qquad
 \sum_{k\in K}\alpha_k(a_k^2+b_k^2)=\frac13.
\end{equation}
Frequency $2$ is excluded to avoid its pairing with the constant frequency.
Prohibiting positive frequencies that differ by $2$ excludes all remaining pairings in $S$.

\begin{lemma}[A summation identity for the weights]\label{v4:seed-lemma-sum-four}
 The weights satisfy
\begin{equation}\label{v4:seed-sum-four}
 \sum_{j\ge1}\alpha_{4j}=\frac1{12}.
\end{equation}
\end{lemma}
\begin{proof}
 Use the partial-fraction expansion of the cotangent,
\[
 \sum_{j\in\mathbb Z}\frac1{n^2j^2-a^2}
 =-\frac{\pi}{na}\cot\frac{\pi a}{n},\qquad a\notin n\mathbb Z,
\]
 which follows directly from $\pi\cot\pi z=1/z+\sum_{j\ge1}2z/(z^2-j^2)$
 by setting $z=a/n$; all the paired series converge absolutely.
 Combining this with \eqref{v4:seed-partial-fractions} gives, for every even $n\ge4$,
\begin{equation}\label{v4:seed-sum-m-n}
 \sum_{j\in\mathbb Z}\mathfrak m(nj)
 =\frac{\pi}{4n}\left(\cot\frac\pi n+\cot\frac{3\pi}n\right).
\end{equation}
 In particular, at $n=4$ the two cotangents are $1,-1$, so the right-hand side vanishes.
 The left-hand side is $1/3-4\sum_{j\ge1}\alpha_{4j}$, proving the claim.
\end{proof}

\begin{theorem}[The full ellipsoid on a nonadjacent support]\label{v4:seed-theorem-sparse}
 If \eqref{v4:seed-sparse-series} defines an $L^2$ function, then $Q=S=0$ and $Lh>0$.
 Setting $A_K=\sum_{k\in K}\alpha_k$, we have $A_K\le1/12$.
 If $K\ne\{4,8,12,\ldots\}$, then $A_K<1/12$ strictly, and the entire ellipsoid on this support satisfies
\begin{equation}\label{v4:seed-sparse-lower-bound}
 Lh\ge\frac13-2\sqrt{\frac{A_K}{3}}>0.
\end{equation}
\end{theorem}
\begin{proof}
 The condition $Q=0$ is exactly the displayed ellipsoid equation; the support assumptions ensure $S=0$.
 To estimate the sum of the weights, first observe that $\alpha_k$ is strictly decreasing for $k\ge4$.
 Set $t=k^2$; then
\[
 \frac{\dd}{\dd t}\frac{t-3}{2(t-1)(t-9)}
 =-\frac{(t-3)^2+12}{2(t-1)^2(t-9)^2}<0.
\]
 Partition all positive even frequencies $\ge4$ into the disjoint pairs
 $\{4j,4j+2\}$,$j=1,2,\ldots$.
 The set $K$ selects at most one element from each pair, whose weight is at most $\alpha_{4j}$.
 Summing gives $A_K\le\sum_j\alpha_{4j}=1/12$.
 Equality requires selecting the larger weight in every pair, hence selecting $4j$ rather than $4j+2$,
 and omitting no pair.
 Thus equality occurs only for $K=\{4,8,12,\ldots\}$.

 Since $L$ has order $-2$, it maps $L^2$ into $H^2$, so $c=Lh$ has a continuous representative.
 Its Fourier series is
\[
 c(s)=\frac13-2\sum_{k\in K}\alpha_k(a_k\cos ks+b_k\sin ks).
\]
 Cauchy--Schwarz in $\ell^2(K;\R^2)$ gives
\begin{align*}
 c(s)&\ge\frac13
 -2\left(\sum_{k\in K}\alpha_k(a_k^2+b_k^2)\right)^{1/2}
    \left(\sum_{k\in K}\alpha_k\right)^{1/2}\\
 &=\frac13-2\sqrt{\frac{A_K}{3}}.
\end{align*}
 The same bound gives absolute uniform convergence of the series for $c$.
 If $A_K<1/12$, the lower bound is strictly positive.
 If $A_K=1/12$ and $c(s_0)=0$, equality in Cauchy--Schwarz forces
\[
 (a_k,b_k)=2(\cos ks_0,\sin ks_0)\qquad(k\in K).
\]
 Indeed, the two weighted sequences have squared norms $1/3$ and $1/12$,
 so their positive proportionality factor is $2$.
 Equality of the weight sum requires $K=\{4,8,12,\ldots\}$, an infinite set.
 The relation $a_k^2+b_k^2=4$ on that set contradicts Parseval's condition for $h\in L^2$.
 Hence $c>0$ also in the equality case, and compactness gives a positive minimum for each fixed seed.
\end{proof}

\subsection{Examples: fourfold frequencies and mixed supports}
\label{v4:seed-sec-fourfold}

For the fourfold-frequency example, take $K=4\mathbb N$ in the general construction.
The resulting real Fourier series has the form
\begin{equation}\label{v4:seed-fourfold-series}
 h(s)=1+\sum_{j\ge1}\bigl(a_j\cos4js+b_j\sin4js\bigr).
\end{equation}
Since all nonzero frequencies are multiples of $4$, no pair of indices in $S$ can sum to $2$,
so $S(h)=0$ holds automatically.
Also,
\begin{equation}\label{v4:seed-fourfold-ellipsoid}
 Q(h)=0
 \quad\Longleftrightarrow\quad
 \sum_{j\ge1}\alpha_{4j}(a_j^2+b_j^2)=\frac13.
\end{equation}
The sine coefficients do not alter this diagonal form.

\begin{corollary}[Positivity for the fourfold-frequency example]\label{v4:seed-theorem-four-positive}
If \eqref{v4:seed-fourfold-series} defines an $L^2$ function satisfying
\eqref{v4:seed-fourfold-ellipsoid}, then $c=Lh$ is continuous and strictly positive.
Every smooth member belongs to $\mathcal A^\infty$.
\end{corollary}
\begin{proof}
Apply Theorem \ref{v4:seed-theorem-sparse} to $K=4\mathbb N$.
Here $A_K=1/12$ by \eqref{v4:seed-sum-four}; the equality argument in that theorem excludes a zero of $c$ for an $L^2$ seed.
\end{proof}

\begin{example}[Single-frequency seeds]\label{v4:seed-example-single-mode}
For each even integer $n\ge4$, take $K=\{n\}$ and choose zero phase.
The ellipsoid equation gives the sequence
\begin{equation}\label{v4:seed-single-mode}
 h(s)=1+\lambda_n\cos ns,
 \qquad
 \lambda_n=\sqrt{\frac{2(n^2-1)(n^2-9)}{3(n^2-3)}}.
\end{equation}
Then $Q=S=0$, and the corresponding drift is
\[
 c(s)=\frac13-\sqrt{\frac{2(n^2-3)}{3(n^2-1)(n^2-9)}}\cos ns.
\]
Its minimum is $1/3-2\sqrt{\alpha_n/3}>0$, since
$\alpha_n\le\alpha_4=13/210<1/12$.
Several members are listed below. All are real analytic, and their
minimal time periods follow from their single nonzero frequency.
\[
\begin{array}{c|c|c|c}
 n&\lambda_n&\min_s c(s)&\text{minimal time period}\\ \hline
 4&\sqrt{70/13}&\displaystyle\frac13-\sqrt{\frac{26}{315}}&\pi/2\\[6pt]
 6&\sqrt{210/11}&\displaystyle\frac13-\sqrt{\frac{22}{945}}&\pi/3\\[6pt]
 8&\sqrt{2310/61}&\displaystyle\frac13-\sqrt{\frac{122}{10395}}&\pi/4\\[6pt]
 10&\sqrt{6006/97}&\displaystyle\frac13-\sqrt{\frac{194}{27027}}&\pi/5
\end{array}
\]
The members with $4\mid n$ belong to the fourfold-frequency family.
The first member of the sequence is
\begin{equation}\label{v4:seed-first-single-mode}
 h(s)=1+\sqrt{\frac{70}{13}}\cos4s
                =1+\frac{\sqrt{910}}{13}\cos4s.
\end{equation}
The corresponding auxiliary function and drift coefficient are
\[
 f(s)=1+\frac6{\sqrt{910}}\cos4s,\qquad
 c(s)=\frac13-\frac{\sqrt{910}}{105}\cos4s>0.
\]
This first member is the only normalized quadratic seed explicitly used for the
construction in \cite[Theorem 3.1, equations (3.1)--(3.3); proof of Theorem 1.1]{Mao}, where it is
denoted by $a_1(s)$ and $a_2=0$.
The scalar parameter multiplying the seed there is accounted for by
Proposition \ref{v4:seed-prop-transformations}.
The characteristic condition for all the single-frequency seeds above is
verified in Section~\ref{v4:char-section}.
\end{example}

\begin{example}[Mixed frequencies $4$ and $10$]\label{v4:seed-example-four-ten}
Take $K=\{4,10\}$, so that
\[
 \alpha_4=\frac{13}{210},\qquad
 \alpha_{10}=\frac{97}{18018}.
\]
The full ellipsoid on this support can be parametrized by
\begin{equation}\label{v4:seed-four-ten-family}
\begin{split}
 h(s)={}&1+\sqrt{\frac{70t}{13}}\cos(4s+\theta_4)
          +\sqrt{\frac{6006(1-t)}{97}}\cos(10s+\theta_{10}),\\
 &0\le t\le1,\qquad \theta_4,\theta_{10}\in\R/(2\pi\mathbb Z).
\end{split}
\end{equation}
The two nonconstant modes contribute $t/3$ and $(1-t)/3$ to the
ellipsoid equation. Hence $Q=0$, and nonadjacency gives $S=0$.
Since $A_K=433/6435$, Theorem \ref{v4:seed-theorem-sparse} gives the uniform bound
\[
 c(s)\ge\frac13-2\sqrt{\frac{433}{19305}}>0
\]
for all parameters.
When $0<t<1$, both modes are present. Their greatest common divisor is $2$,
so the minimal time period is $\pi$, independently of the phases.
In particular, these seeds do not belong to the fourfold-frequency family.
At $t=1$ and $t=0$, the minimal periods are $\pi/2$ and $\pi/5$, respectively.
Corollary~\ref{v43:char-four-ten-endpoints} proves that neighborhoods of
both endpoints contain mixed members satisfying the characteristic
condition for geometric realization.

For a concrete member, take $t=1/2$ and $\theta_4=\theta_{10}=0$. Then
\[
 h(s)=1+\sqrt{\frac{35}{13}}\cos4s
        +\sqrt{\frac{3003}{97}}\cos10s,
\]
and its balance conditions can be checked directly:
\[
 Q=\frac13-\frac{13}{210}\frac{35}{13}
             -\frac{97}{18018}\frac{3003}{97}=0,
 \qquad S=0.
\]
For this member, the triangle inequality gives the sharper bound
\[
 c(s)\ge\frac13-\frac{\sqrt{455}}{105}
                    -\frac{\sqrt{291291}}{9009}
 >\frac{137}{2145}>0.
\]
The last strict inequality uses only $\sqrt{455}<22$ and $\sqrt{291291}<540$.
Thus positivity does not rely on numerical approximation.
\end{example}

\subsection{Finite-dimensional parameters, nonuniform lower bounds, and degenerate limits}
\label{v4:seed-sec-finite-dimensional}

For any finite nonadjacent support $K$ with $J$ elements, the coefficients in
\eqref{v4:seed-sparse-series} form an ellipsoid of dimension $2J-1$.
All its members are trigonometric polynomials in $\mathcal A^{\mathrm{trig}}$ by
Theorem~\ref{v4:seed-theorem-sparse}.
The choice $K=\{4\}$ gives the first single-frequency seed, while $K=\{4,10\}$
gives minimal period $\pi$ when both frequencies occur.
The dimension counts seed parameters; it is not the dimension of a moduli space of unmarked metrics.

Return to the fourfold-frequency family \eqref{v4:seed-fourfold-series}.
Retain only $1\le j\le J$, and set
\[
 A_J=\sum_{j=1}^J\alpha_{4j},\qquad
 \delta_J=\frac13-2\sqrt{\frac{A_J}{3}}.
\]
Since $A_J<1/12$, we have $\delta_J>0$.
The same proof gives the uniform estimate over the entire finite-dimensional ellipsoid,
\begin{equation}\label{v4:seed-finite-bound}
 c(s)\ge\delta_J\qquad\text{for all parameters and all }s.
\end{equation}
This lower bound is not merely sufficient; it is attained.
Take
\begin{equation}\label{v4:seed-extremizing-truncation}
 a_1=\cdots=a_J=(3A_J)^{-1/2},\qquad b_1=\cdots=b_J=0,
\end{equation}
Then the ellipsoid constraint holds and $c(0)=\delta_J$.
Since $A_J\uparrow1/12$, we have $\delta_J\downarrow0$.

\begin{corollary}[There is no uniform positive margin over the full family]\label{v4:seed-cor-no-uniform-gap}
 Even among the fourfold-frequency trigonometric polynomial members, there is no uniform constant $\delta>0$
 such that all balanced seeds satisfy $Lh\ge\delta$.
 Each fixed-bandwidth ellipsoid $J$, however, has the uniform positive margin in \eqref{v4:seed-finite-bound}.
\end{corollary}

The finite-dimensional geometry is also explicit.
When both sines and cosines are allowed, the transformation
$(a_j,b_j)\mapsto\sqrt{3\alpha_{4j}}(a_j,b_j)$
turns the constraint into the unit sphere in $\R^{2J}$, of dimension $2J-1$.
Retaining only cosines gives an ellipsoid surface in $\R^J$, of dimension $J-1$.
The former has only the quarter-period translation symmetry, whereas the latter also retains time-reflection symmetry.
General phase parameters cannot all be removed by a single time translation:
time translation has only one real parameter, whereas the phases at distinct frequencies can vary independently.

\begin{remark}[A distributional boundary limit]
In the formal limit of \eqref{v4:seed-extremizing-truncation}, each fixed $a_j$ tends to $2$.
In the sense of distributions, this corresponds to
\[
 1+2\sum_{j\ge1}\cos4js
 =\frac\pi2\sum_{n\in\mathbb Z}\delta(s-n\pi/2).
\]
This distributional limit lies outside $L^2$ and describes the loss of a uniform positivity bound.
\end{remark}

\subsection{General balance does not automatically imply positivity}
\label{v4:seed-sec-counterexample}

The quadratic-form weight of frequency $2$ is positive, and this frequency couples to both the constant and frequency $4$.
Such coupling can achieve $Q=S=0$ while destroying $c>0$.
The following seed satisfies all three balance equations although its drift changes sign.

\begin{example}[Balance holds but $c$ changes sign]\label{v4:seed-example-sign-change}
 Take
\begin{equation}\label{v4:seed-bad-seed}
 h(s)=1+a\cos2s+\frac{10}{3}\cos4s,
 \qquad a^2=\frac{670}{63}.
\end{equation}
 For a general $1+a\cos2s+b\cos4s$, \eqref{v4:seed-Q} and \eqref{v4:seed-S} give directly
\[
 Q=\frac13+\frac{a^2}{30}-\frac{13b^2}{210},
 \qquad S=-\frac a3+\frac{ab}{10}.
\]
 The first term of $S$ comes from $(0,2),(2,0)$, and the second from $(-2,4),(4,-2)$.
 Substituting $b=10/3$ and $a^2=670/63$ makes both expressions zero.
 However,
\[
 c(s)=\frac13+\frac a{15}\cos2s-\frac{26}{63}\cos4s.
\]
 Choose $\cos2s=-\operatorname{sign}(a)$, so that $\cos4s=1$ and hence
\[
 c(s)=-\frac5{63}-\frac{|a|}{15}<0.
\]
 The mean of $c$ remains $1/3$, so it does change sign.
\end{example}

Thus cubic balance alone does not imply the positivity required in Definition \ref{v4:seed-class-A}.

\begin{example}[A second reference seed]\label{v4:seed-ten}
The single-frequency member
\[
 h_{10}(s)=1+\sqrt{\frac{6006}{97}}\cos10s
\]
is balanced and has $c=Lh_{10}>0$. Its minimal time period is $\pi/5$, so it is
not a member of the fourfold-frequency subfamily. Section~\ref{v4:char-section}
verifies the nonzero-point characteristic condition for this seed and, more generally,
for every even single frequency $n\ge4$. The mixed $4,10$ family above illustrates
that the algebraic seed construction also allows independent frequencies and phases;
only the members satisfying the additional characteristic condition enter the
exact-closure theorem.
\end{example}

\begin{remark}[What the seed calculation establishes]\label{v4:seed-examples}
The positive families constructed above satisfy the cubic balance equations and $c>0$; the last counterexample separates these two requirements.
They supply the finite-order preparation proved next. The geometric realization
further uses the characteristic condition of Section~\ref{v4:char-section};
for the coefficient ellipsoids, this condition is established on the relative
neighborhoods specified there.
\end{remark}
\section{Formal recursion and finite-order initialization}
\label{v4:formal-chapter}

Let $h$ be a real smooth seed satisfying
\begin{equation}\label{v4:formal-seed-assumptions}
 h(s+\pi)=h(s),\qquad
 \operatorname{avg}h=1,\qquad
 Q(h)=S(h)=0,\qquad c=Lh>0.
\end{equation}
Here the average is taken over $[0,2\pi]$; $Q,S,L$ are the quadratic forms and multiplier defined in the section on seeds.

We seek a formal control of the form
\begin{equation}\label{v4:formal-control-series}
 \widehat H(q,s)=h(s)q^2+\sum_{j\ge4}\gamma_j(s)q^j.
\end{equation}
In the finite Taylor expansions below, $\gamma_2=h$ and $\gamma_3=0$
for this initial control.
We shall prove that each $\gamma_j$ can be chosen to be a real $\pi$-periodic trigonometric polynomial,
with the formal return map equal to the identity.
No convergence of the formal series is asserted here.

\subsection{Rotating coordinates and the recursion for finite jets}
\label{v4:formal-finite-jets}

Set $T=2\pi$ and write
\[
 C_s=\begin{pmatrix}\cos s&\sin s\\-\sin s&\cos s\end{pmatrix},
 \qquad m_s=\binom{-\sin s}{\cos s},
 \qquad \ell(x,s)=m_s\cdot x.
\]
For the oscillator
\[
 p_s=q,\qquad q_s=-p+H(q,s),\qquad (p,q)|_{s=0}=x,
\]
let $\Psi_s(x)$ denote the flow in the original coordinates and let
$\Phi_s(x)=C_s^{-1}\Psi_s(x)$.
Direct differentiation gives
\begin{equation}\label{v4:formal-rotating-ode}
 \partial_s\Phi_s(x)=m_sH(m_s\cdot\Phi_s(x),s),
 \qquad \Phi_0(x)=x.
\end{equation}
Since $C_T=\Id$, the full returns in the original and rotating coordinates agree.
We shall also repeatedly use
\[
 \ell(x,s+\pi)=-\ell(x,s),\quad
 \partial_s^2\ell=-\ell,\quad
 m_s\cdot m_t=\cos(s-t),\quad
 \det(m_t,m_s)=\sin(s-t).
\]

\begin{lemma}[The return jet depends only on the corresponding finite control jet]
\label{v4:formal-jet-dependence}
Suppose that two smooth controls $H,\widetilde H$ vanish to second order at $q=0$
and satisfy, for $0\le j\le N$,
$\partial_q^jH(0,s)=\partial_q^j\widetilde H(0,s)$.
Then at $x=0$ their flows and full returns have the same jet of order $N$.
The conclusion holds uniformly in time on any common interval of existence of the flows.
\end{lemma}
\begin{proof}
The integral equation gives an actual flow on a common compact time interval near the equilibrium: a sufficiently small initial ball stays inside a fixed control tube by Gronwall's inequality. Differentiating that integral equation produces linear Volterra equations for each successive initial-value derivative. Their unique solutions are continuous, and induction gives joint smoothness. This argument does not require time periodicity.
Taylor expansion on that interval therefore applies to these flows.
Write the finite Taylor expansion of the control as
$H(q,s)=\sum_{j=2}^N \gamma_j(s)q^j+O(q^{N+1})$,
and the homogeneous expansion of the flow as
$\Phi_s(x)=\sum_{n=1}^NU_n(x,s)+O(|x|^{N+1})$.
Here $U_1=x$ and $U_n(x,0)=0$ for $n\ge2$.
Taking successive homogeneous parts of \eqref{v4:formal-rotating-ode} gives
\begin{equation}\label{v4:formal-jet-recursion}
 \partial_sU_n
 =m_s\sum_{j=2}^n\gamma_j(s)
   \sum_{\substack{r_1+\cdots+r_j=n\\r_i\ge1}}
   \prod_{i=1}^j(m_s\cdot U_{r_i}),\qquad n\ge2.
\end{equation}
Since $j\ge2$, each $r_i$ in the inner sum is at most $n-1$.
Thus the $n$th equation depends only on $\gamma_2,\ldots,\gamma_n$ and the lower-order flow terms already determined.
Successive integration with the same zero initial values gives the same $U_2,\ldots,U_N$ for both controls.
The same conclusion follows by repeatedly differentiating the actual smooth ordinary differential equation with respect to its initial value,
so no convergence of the formal flow is presupposed.
Evaluation at $s=T$ proves the assertion about return jets.
\end{proof}

\begin{lemma}[The direct return produced by an added term]
\label{v4:formal-direct-return}
Let the background satisfy $H=O(q^2)$, and add $k(s)q^n$, where $n\ge2$.
Then the return terms of degree below $n$ are unchanged, and the change in degree $n$ is
\begin{equation}\label{v4:formal-direct-moment}
 (\mathcal L_nk)(x)=\int_0^T k(s)m_s\ell(x,s)^n\,\dd s
 =\frac1{n+1}\nabla_x\int_0^Tk(s)\ell(x,s)^{n+1}\,\dd s.
\end{equation}
If $k$ is $\pi$-periodic and $n$ is even, then $\mathcal L_nk=0$.
\end{lemma}
\begin{proof}
In \eqref{v4:formal-jet-recursion}, the new coefficient first appears in degree $n$.
It is accompanied by exactly $n$ linear flow terms, so its contribution is $k(s)m_s\ell^n$.
Integration gives the first equality; the second follows from $\nabla_x\ell=m_s$.
When $n$ is even, $m_{s+\pi}=-m_s$, whereas
$k(s+\pi)\ell(x,s+\pi)^n=k(s)\ell(x,s)^n$.
The integrals over the two half-periods therefore cancel exactly.
\end{proof}

\subsection{The parity step in the formal induction}
\label{v4:formal-leading-parity}

Lemma~\ref{v4:seed-lemma-leading-odd} shows that the first nonzero return error of a half-periodic control has odd degree.
We use it only after all lower-degree errors have been eliminated.
Once the error of odd degree $m+1$ is removed, a term of degree $m+2$ would be a leading even-degree term
and must therefore vanish. This permits the induction to advance by two degrees at each step.

For the initial quadratic control, the balance assumptions \eqref{v4:formal-seed-assumptions}
and Corollary~\ref{v4:seed-cor-order-five} give
\begin{equation}\label{v4:formal-initial-flat-order}
 P_{hq^2}(x)=x+O(|x|^5).
\end{equation}
Thus the first step treats the possible error of degree five.

\subsection{Direct gradient moments cover all gradient targets}
\label{v4:formal-gradient-section}

Let $\mathscr P_N$ denote the space on $\R^2$ of real homogeneous scalar polynomials of degree $N$,
which has dimension $N+1$.
If $N$ is even, the restriction of any $\psi\in\mathscr P_N$ to the unit circle
belongs to
\[
 \operatorname{span}\{1,\cos2\theta,\sin2\theta,
                  \ldots,\cos N\theta,\sin N\theta\}.
\]
Conversely, each angular expression in this space, multiplied by $r^N$, is a homogeneous polynomial.

\begin{lemma}[Circular moments of even powers form an isomorphism]
\label{v4:formal-ridge-isomorphism}
Let $N\ge2$ be even and set
\[
 \mathscr T_N=\operatorname{span}\{1,\cos2s,\sin2s,
                       \ldots,\cos Ns,\sin Ns\}.
\]
The map $a\mapsto\int_0^T a(s)\ell(x,s)^N\dd s$
is a linear isomorphism from $\mathscr T_N$ to $\mathscr P_N$.
\end{lemma}
\begin{proof}
Write $x=r(\cos\theta,\sin\theta)$, so that $\ell=r\sin(\theta-s)$.
For $0\le j\le N/2$, expanding $(e^{it}-e^{-it})^N/(2i)^N$ gives
\[
 \kappa_{N,j}:=\int_0^{2\pi}\cos(2jt)\sin^Nt\,\dd t
 =2\pi(-1)^j2^{-N}\binom{N}{N/2-j}\ne0.
\]
For $j=0$, this formula gives the integral of the constant mode.
Changing variables by $t=\theta-s$ and using the vanishing of the odd integrands yields
\[
 \int_0^T\cos(2js)\ell^N\dd s
   =\kappa_{N,j}r^N\cos(2j\theta),
\]
\[
 \int_0^T\sin(2js)\ell^N\dd s
   =\kappa_{N,j}r^N\sin(2j\theta),\qquad j\ge1.
\]
Thus the matrix in the corresponding real Fourier bases is diagonal, with every diagonal entry nonzero.
\end{proof}

\begin{corollary}[Odd-degree controls eliminate arbitrary gradient errors]
\label{v4:formal-gradient-surjective}
If $n\ge3$ is odd, then $\mathcal L_n$ maps the real $\pi$-periodic trigonometric polynomials
surjectively onto $\nabla\mathscr P_{n+1}$.
Only the frequencies $0,2,\ldots,n+1$ are needed in the control.
\end{corollary}
\begin{proof}
Given $\nabla\psi$, use Lemma \ref{v4:formal-ridge-isomorphism}
to choose $k$ such that $\int k\ell^{n+1}=(n+1)\psi$.
Substitution into \eqref{v4:formal-direct-moment} proves the claim.
\end{proof}

\subsection{The double integral for the mixed return and the Green kernel}
\label{v4:formal-mixed-section}

Fix an even integer $m\ge4$, assume that the background satisfies $H=hq^2+O(q^4)$,
and add $k(s)q^m$.
Write
\[
 X_2(x,s)=h(s)m_s\ell(x,s)^2,
 \qquad Y_m(x,s)=k(s)m_s\ell(x,s)^m.
\]
The quadratic flow term of the background is $U_2(x,s)=\int_0^sX_2(x,t)\dd t$.
The added flow term of degree $m$ is $\Delta U_m=\int_0^sY_m(x,t)\dd t$.
Its endpoint value vanishes by Lemma \ref{v4:formal-direct-return}.

\begin{lemma}[The mixed return of an even-degree control and the seed]
\label{v4:formal-mixed-return}
The added control does not change the return coefficients through degree $m$, and its change in degree $m+1$ is
\begin{equation}\label{v4:formal-mixed-double}
 \Delta B_{m+1}(x)=\int_0^T\int_0^s
 \left\{DX_2(x,s)[Y_m(x,t)]
       +DY_m(x,s)[X_2(x,t)]\right\}\dd t\dd s.
\end{equation}
This expression is linear in $k$.
\end{lemma}
\begin{proof}
Compare the terms of degree $m+1$ in \eqref{v4:formal-jet-recursion}.
The added term produces $DX_2[\Delta U_m]$ in the background quadratic vector field
and $DY_m[U_2]$ in the added vector field.
Interactions of background terms of degree at least four with the added displacement have degree at least $m+3$.
The quadratic self-interaction of the added term has degree $2m-1>m+1$.
Thus precisely the two displayed terms occur in the current degree, and integration gives the formula.
\end{proof}

We use the convention $\Curl V=\partial_{x^1}V^2-\partial_{x^2}V^1$.
For any $x$-independent vector $v$, one has $\Curl(vf)=\det(\nabla f,v)$.
Termwise differentiation therefore gives
\begin{align*}
 \Curl(DX_2(s)[Y_m(t)])
 &=2m h(s)k(t)\cos(s-t)\sin(s-t)\ell_s\ell_t^{m-1},\\
 \Curl(DY_m(s)[X_2(t)])
 &=2m k(s)h(t)\cos(s-t)\sin(s-t)\ell_s^{m-1}\ell_t.
\end{align*}
In this paragraph, $\ell_s,\ell_t$ denote only the values $\ell(x,s),\ell(x,t)$;
time derivatives are always written as $\partial_s\ell$.
The determinant is ordered as $(m_t,m_s)$, so both terms have the same positive sign.

Fix $x$ and set $a(s)=h(s)\ell(x,s)$ and $g(s)=k(s)\ell(x,s)^{m-1}$.
Using $2\sin u\cos u=\sin2u$ and interchanging the variables on one triangular integration region gives
\begin{equation}\label{v4:formal-green-square}
 \Curl\Delta B_{m+1}(x)
 =m\int_0^Tg(t)(\mathscr Ka)(t)\dd t,
 \qquad
 (\mathscr Ka)(t)=\int_0^T\sin(2|s-t|)a(s)\dd s.
\end{equation}
On the two sides, the kernel is respectively $\sin2(s-t)$ and $\sin2(t-s)$;
the absolute value cannot be removed to turn it into an antisymmetric kernel.

Splitting the integral at $s=t$ and differentiating twice directly verifies
\begin{equation}\label{v4:formal-green-ode}
 (\partial_t^2+4)\mathscr Ka=4a.
\end{equation}
In the first derivative the sine vanishes at the splitting point; in the second derivative each side contributes $2a(t)$.
The remaining integral is exactly $-4\mathscr Ka$.

Since $a(s+\pi)=-a(s)$, it contains only odd time frequencies.
The operator $\partial_s^2+4$ is invertible on this antiperiodic space:
its multiplier at frequency $r\in2\mathbb Z+1$ is $4-r^2\ne0$.
Let $W=(\partial_s^2+4)^{-1}a$ be the unique antiperiodic solution.
By \eqref{v4:formal-green-ode}, on $[0,T]$ we have
\[
 \mathscr Ka=4W+A(x)\cos2s+B(x)\sin2s.
\]
Since $g$ is also antiperiodic, the integrals of the products of the last two terms with $g$ over a full period vanish.
Consequently,
\begin{equation}\label{v4:formal-resolvent-curl}
 \Curl\Delta B_{m+1}(x)=4m\int_0^Tk(s)\ell(x,s)^{m-1}W(x,s)\dd s.
\end{equation}

\subsection{Coefficient equations for general seeds and the one-dimensional control operator}
\label{v4:formal-coefficient-section}

Use the auxiliary functions $f,c,d$ defined in
\eqref{v4:seed-f-equation}--\eqref{v4:seed-cd-definition}; their multipliers are given by \eqref{v4:seed-cd-fourier}.
In particular, $c=Lh$ is the positive coefficient prescribed in the seed condition.
The calculation in the proof of Lemma~\ref{v4:seed-lemma-cd-identity} gives
\begin{equation}\label{v4:formal-cd-system}
 c''+3c-2d'=h,\qquad 2c'+d''+3d=0.
\end{equation}
Since $\partial_s\ell=-r$, equation \eqref{v4:seed-cd-identity} says that
$c\ell+d\partial_s\ell$ solves the resolvent equation
$(\partial_s^2+4)W=h\ell$.
This expression is antiperiodic, so the uniqueness used above identifies it with $W$:
\begin{equation}\label{v4:formal-W-explicit}
 W=c\ell+d\partial_s\ell.
\end{equation}

The mixed identity in the following statement is the one in
\cite[equations (3.22)--(3.24)]{Mao}, written with the coefficients
$c,d$ of \eqref{v4:seed-cd-definition}.
\begin{theorem}[The mixed curl formula for general half-periodic seeds]
\label{v4:formal-mixed-curl}
For even $m\ge4$ and any real $\pi$-periodic $k$, we have
\begin{equation}\label{v4:formal-T-and-curl}
 T_mk=ck-\frac1m(dk)',\qquad
 \Curl\Delta B_{m+1}(x)=4m\int_0^T(T_mk)(s)\ell(x,s)^m\dd s.
\end{equation}
\end{theorem}
\begin{proof}
Substitute \eqref{v4:formal-W-explicit} into \eqref{v4:formal-resolvent-curl}.
Using $\ell^{m-1}\partial_s\ell=\partial_s(\ell^m)/m$ gives
\[
 \Curl\Delta B_{m+1}
 =4m\int_0^T\left(ck\ell^m+\frac{dk}{m}\partial_s(\ell^m)\right)\dd s.
\]
Integrate the second term by parts.
The boundary term $[dk\ell^m]_0^T$ vanishes by full periodicity.
The remaining expression is \eqref{v4:formal-T-and-curl}.
In particular, in $(dk)'=d'k+dk'$, the term $d'k$ must be retained.
\end{proof}

\subsection{Surjectivity of the even-degree curl moments}
\label{v4:formal-maximum-principle-section}

\begin{theorem}[Surjectivity of mixed curl moments]
\label{v4:formal-curl-surjectivity}
Let $m\ge2$ be even and let $c,d$ be real smooth $\pi$-periodic functions.
Suppose either that $c>0$ or that $c,d$ are given by
\eqref{v4:seed-cd-definition} for a real trigonometric polynomial seed
with highest positive Fourier frequency $K\ge4$.
Then the map
\[
 \mathcal A_m:k\longmapsto
  \int_0^T\left(ck-\frac1m(dk)'\right)\ell(x,s)^m\dd s
\]
maps real smooth $\pi$-periodic functions surjectively onto $\mathscr P_m$.
It remains surjective when the input is restricted to real $\pi$-periodic trigonometric polynomials.
\end{theorem}
\begin{proof}
Suppose a linear functional $\Lambda\in\mathscr P_m^*$ on the target annihilates the entire image.
Define
\[
 u(s)=\Lambda\bigl(\ell(\,\cdot\,,s)^m\bigr).
\]
This is a real trigonometric polynomial in $\mathscr T_m$.
For every smooth periodic $k$, integration by parts gives
\[
 0=\Lambda(\mathcal A_mk)
  =\int_0^Tk(s)\left(c(s)u(s)+\frac{d(s)}m u'(s)\right)\dd s.
\]
The coefficient in the integrand is smooth and periodic, so
\begin{equation}\label{v4:formal-adjoint-ode}
 mc(s)u(s)+d(s)u'(s)=0\quad\text{for every }s.
\end{equation}
Suppose $u\ne0$. If $c>0$, take a maximum point $s_0$ of $|u|$.
Then $u(s_0)\ne0$, $u'(s_0)=0$, and $c(s_0)>0$, contrary to
\eqref{v4:formal-adjoint-ode}.

For a polynomial seed as in the second hypothesis, the Fourier multiplier
in \eqref{v4:seed-f-equation} preserves its highest frequency $K$.
Let $\ell$ be the highest Fourier frequency of $u$, taking $\ell=0$
for constant $u$. Since $u$ is real, $0\le\ell\le m$ and
$\widehat u_\ell\ne0$. At frequency $K+\ell$ the sole contribution
in \eqref{v4:formal-adjoint-ode} is
\[
 \frac{\widehat f_K\widehat u_\ell}{9}
      \bigl[m(3-K^2)+2K\ell\bigr].
\]
This follows from $\widehat c_K=(3-K^2)\widehat f_K/9$ and
$\widehat d_K=-2iK\widehat f_K/9$. The bracket is at most
$-m(K-3)(K+1)<0$, again a contradiction. Thus $u=0$ under either
hypothesis.
Lemma \ref{v4:formal-ridge-isomorphism} shows that
$\{\ell(\,\cdot\,,s)^m:s\in[0,T]\}$ spans $\mathscr P_m$, and hence $\Lambda=0$.
There is no nonzero annihilator in the finite-dimensional target, proving surjectivity.

Finally, trigonometric polynomials are dense among smooth periodic functions in the $C^1$ norm,
and $\mathcal A_m:C^1\to\mathscr P_m$ is continuous.
The image of the trigonometric polynomials is a linear subspace of the finite-dimensional target, hence closed.
This image is also dense in the entire target, so it equals the entire target.
\end{proof}

\begin{remark}[The mean of a higher-order control]
The control space includes the constant time mode. If one imposes $\operatorname{avg}k=0$, the dual equation implies only that
$mcu+du'$ is constant; the preceding maximum principle argument then no longer applies as written.
The normalization $\operatorname{avg}h=1$ concerns the quadratic seed
and places no constraint on the mean of the higher-degree control coefficients $\gamma_m$.
\end{remark}

For each fixed $m$, finitely many Fourier modes can be selected
whose images form a basis of $\mathscr P_m$, and the resulting finite linear system can be solved.
The selected frequencies and the size of the solution may grow with $m$.
Only finitely many steps will be used for each closure construction. No uniform bound on the entire sequence of formal coefficients is needed.

\subsection{Transversality of the cubic constraints: a separate check of the quadratic variation}
\label{v4:formal-transversality-section}

The preceding mixed return lemma for a finite addition $kq^m$ requires $m\ge4$.
For $m=2$, a finite addition also creates a cubic self-interaction of $k$,
so that conclusion cannot be applied unchanged.
However, under the first-order parameter variation $h_\varepsilon=h+\varepsilon k$ of the seed,
these additional self-interactions are $O(\varepsilon^2)$.
Starting from
\[
 B_3[h](x)=\int_0^T\int_0^sDX_h(x,s)[X_h(x,t)]\dd t\dd s,
 \qquad X_h=h(s)m_s\ell^2,
\]
and differentiating directly with respect to $\varepsilon$ gives exactly two mixed terms.
The preceding Green kernel calculation still applies, and hence
\begin{equation}\label{v4:formal-m2-variation}
 D_h(\Curl B_3)[k]
 =8\int_0^T\left(ck-\frac12(dk)'\right)\ell(x,s)^2\dd s.
\end{equation}

\begin{proposition}[The constraints are transverse throughout the balanced family with positive drift]
\label{v4:formal-transversality}
In the mean-$1$ hyperplane of real smooth $\pi$-periodic functions,
the set $\{h:Q(h)=0,S(h)=0,Lh>0\}$ is an everywhere smooth submanifold of codimension $3$.
More precisely, the derivative of the constraint map
$\mathcal G(h)=(Q(h),\operatorname{Re}S(h),\operatorname{Im}S(h))$
at every point of this set remains surjective when restricted to zero-mean directions.
\end{proposition}
\begin{proof}
The cubic formula in the section on seeds reads
\[
 B_3[h](x)=\pi\{Q(h)|x|^2+\operatorname{Re}(S(h)z^2)\}Jx,
 \quad z=x^1+ix^2,\quad Jx=(-x^2,x^1).
\]
For a homogeneous quadratic scalar $a$, $\Curl(aJx)=2a+x\cdot\nabla a=4a$.
Thus the linear correspondence from the three constraints to $\Curl B_3\in\mathscr P_2$ is an isomorphism.
Equation \eqref{v4:formal-m2-variation} and Theorem \ref{v4:formal-curl-surjectivity}
show that $D\mathcal G(h)$ is surjective on the full space of periodic directions.

At a balanced point, quadratic homogeneity gives
$D\mathcal G(h)[h]=2\mathcal G(h)=0$.
Any direction $k$ can be replaced by $k_0=k-(\operatorname{avg}k)h$.
Since $\operatorname{avg}h=1$, the function $k_0$ has zero mean, and
$D\mathcal G(h)[k_0]=D\mathcal G(h)[k]$.
Surjectivity therefore persists upon restriction to the normalized hyperplane.

To make the function-space statement precise, one may first work in any real periodic Sobolev space of sufficiently high order.
The functions $Q,S$ are continuous quadratic functions, and $h\mapsto Lh$ maps continuously into the continuous functions,
so $Lh>0$ is an open condition.
The kernel of the surjective derivative has a finite-dimensional complement, and the Banach implicit function theorem gives a local graph of codimension three.
Alternatively, in the smooth space one may choose three fixed smooth control directions as a complement
and apply the implicit function theorem only to these three real parameters, obtaining the same local representation compatible with all smooth seminorms.
\end{proof}

\subsection{Induction eliminating the full vector error}
\label{v4:formal-induction-section}

\begin{lemma}[An explicit potential for a homogeneous curl-free vector field]
\label{v4:formal-gradient-potential}
Let $V$ be a homogeneous vector polynomial of degree $n$, with $\Curl V=0$.
Then
\[
 \psi(x)=\frac{x\cdot V(x)}{n+1}\in\mathscr P_{n+1},
 \qquad \nabla\psi=V.
\]
\end{lemma}
\begin{proof}
The curl-free condition gives $\partial_iV_j=\partial_jV_i$.
Thus
$\partial_i(x_jV_j)=V_i+x_j\partial_iV_j
=V_i+x_j\partial_jV_i=(n+1)V_i$,
where the last step uses Euler's homogeneous identity.
\end{proof}

\begin{theorem}[Formal identity return for balanced seeds]
\label{v4:formal-formal-identity}
For every smooth seed $h$ satisfying \eqref{v4:formal-seed-assumptions},
or every real mean-one $\pi$-periodic trigonometric polynomial seed
with $Q(h)=S(h)=0$ and highest frequency at least four,
there exist real $\pi$-periodic trigonometric polynomials $\gamma_j$, $j\ge4$,
such that the formal return of \eqref{v4:formal-control-series} is the identity.
More specifically, at each completed odd truncation order $M\ge5$,
\[
 H_M=hq^2+\sum_{j=4}^M\gamma_jq^j
 \quad\text{satisfies}\quad
 P_{H_M}(x)=x+O(|x|^{M+2}).
\]
\end{theorem}
\begin{proof}
By \eqref{v4:formal-initial-flat-order}, the quadratic seed already satisfies $P=x+O(|x|^5)$.
Let $m\ge4$ be even, suppose the control coefficients of degree below $m$ have been chosen, and assume
\[
 P_H(x)=x+B_{m+1}(x)+O(|x|^{m+2}).
\]
Here $B_{m+1}$ is a general homogeneous vector polynomial of odd degree, and its curl belongs to $\mathscr P_m$.

By Theorem \ref{v4:formal-curl-surjectivity}, choose a real periodic trigonometric polynomial $k_m$ such that
\[
 4m\mathcal A_mk_m=-\Curl B_{m+1}.
\]
Adding $k_mq^m$ does not change the lower-order return terms already eliminated,
and Theorem \ref{v4:formal-mixed-curl} shows that the new term of degree $m+1$, denoted by $\widetilde B_{m+1}$, is curl-free.
By Lemma \ref{v4:formal-gradient-potential}, it is the gradient of some
$\psi\in\mathscr P_{m+2}$.

Next, by Corollary \ref{v4:formal-gradient-surjective}, choose $k_{m+1}$ such that
$\mathcal L_{m+1}k_{m+1}=-\nabla\psi$.
Adding $k_{m+1}q^{m+1}$ eliminates the entire vector error of degree $m+1$, giving
$P=x+O(|x|^{m+2})$.
The new control is still $\pi$-periodic, so Lemma~\ref{v4:seed-lemma-leading-odd}
also eliminates the possible leading term of degree $m+2$, and in fact gives
\[
 P=x+O(|x|^{m+3}).
\]
This completes the induction step from $m$ to $m+2$.

Proceed successively with $m=4,6,8,\ldots$, and set
$\gamma_m=k_m$, $\gamma_{m+1}=k_{m+1}$.
Each step adds only controls of higher degree, so the return jets already determined remain unchanged.
At any fixed degree, the corresponding return coefficient vanishes after finitely many steps,
which is precisely the meaning of formal identity return.
All added coefficients are obtained from finite Fourier controls.
\end{proof}

In particular the sign-changing drift in Example~\ref{v4:seed-example-sign-change}
allows formal identity return at all orders. The geometric realization argument
also requires the characteristic hypothesis of Section~\ref{v4:char-section}.

\subsection{Finite truncations on a fixed nonzero initial-value region}

The formal calculation at the equilibrium is used only to prepare a small-parameter
approximation. The region on which the return will be corrected need not contain
that equilibrium.

\begin{theorem}[Finite-order approximation on fixed regions]\label{v4:formal-approximation}
Fix $h\in\mathcal A^\infty$, an odd integer $M\ge5$, and a finite truncation
$H_M$ from Theorem~\ref{v4:formal-formal-identity}. Put
\begin{equation}\label{v4:formal-scaled-control}
 H^{[M]}_{\varepsilon,h}(q,s)
 =\varepsilon^{-1}H_M(\varepsilon q,s)
 =\varepsilon h(s)q^2+\sum_{j=4}^M\varepsilon^{j-1}\gamma_j(s)q^j.
\end{equation}
The polynomial on the right defines the value at $\varepsilon=0$.
For every fixed compact initial-value set $K\subset\R^2$, the actual flows exist
on $[0,2\pi]$ for all sufficiently small $|\varepsilon|$. For every finite $k$,
\begin{equation}\label{v4:formal-residual-estimate}
 \|P_{H^{[M]}_{\varepsilon,h}}-\Id\|_{C^k(K)}
 \le C_{M,k,K,h}|\varepsilon|^{M+1}.
\end{equation}
The initial controls satisfy
$H^{[M]}_{\varepsilon,h}=\varepsilon hq^2+O(\varepsilon^3)$ in every fixed
smooth norm on a compact control tube.
\end{theorem}
\begin{proof}
Let $P_M$ denote the actual return map of $H_M$ near zero. The finite jet theorem
and the formal induction give a smooth vector-valued function $R_M$ with
\[
 P_M(\zeta)-\zeta
   =\sum_{|\alpha|=M+2}\zeta^\alpha R_{M,\alpha}(\zeta),
 \qquad \zeta\in\R^2.
\]
This follows from the integral form of the multivariable Taylor remainder;
each $R_{M,\alpha}$ is smooth on a fixed small ball. The rescaled actual flow is
\[
 \Psi_{\varepsilon}(x,s)=\varepsilon^{-1}\Psi_M(\varepsilon x,s)
 \qquad(\varepsilon\ne0).
\]
Differentiating shows that it solves the oscillator with control
\eqref{v4:formal-scaled-control} and initial value $x$.
For $x$ in a fixed neighborhood of $K$, its trajectories exist throughout the
period when $|\varepsilon|$ is small, because the corresponding initial values
$\varepsilon x$ lie in the common small ball for $\Psi_M$.
Consequently
\[
 P_{H^{[M]}_{\varepsilon,h}}(x)-x
 =\varepsilon^{M+1}\sum_{|\alpha|=M+2}x^\alpha
                      R_{M,\alpha}(\varepsilon x).
\]
Differentiation in $x$ to any fixed finite order proves
\eqref{v4:formal-residual-estimate}; the constants are uniform on $K$.
The polynomial expression for the control gives the last assertion directly.
It also gives joint smoothness of the rescaled flow through $\varepsilon=0$
by the integral equation, where the limiting flow is the harmonic rotation.
\end{proof}

A fixed smooth velocity cutoff can be inserted in
\eqref{v4:formal-scaled-control}. Choose its plateau to contain all velocity
values on the harmonic trajectories from a neighborhood of $K$, with a positive
margin. The rescaled flows remain inside that plateau for small parameters by
continuous dependence and a finite-time Gronwall estimate. Uniqueness then
identifies the cutoff and uncut flows there. Thus the same approximation can be
used in a fixed compact control cylinder without changing its return estimate.

\subsection{Finite selections near a reference seed}

The finite preparation respects local seed parameters. We state its
dependence explicitly, since it allows the same initialization procedure
to be used on the seed neighborhoods constructed in the next section.

\begin{proposition}[Smooth selection of a finite preparation]
\label{v43:formal-finite-selection}
Fix $h^\circ\in\mathcal A^\infty$ and an odd integer $M\ge5$.
There are a neighborhood $\mathcal V$ of $h^\circ$ in a finite $C^{r_0}$
norm, with $r_0\ge2$, and finite-dimensional spaces $E_j$ of real
$\pi$-periodic trigonometric polynomials, $4\le j\le M$, such that the
coefficients
\[
 \gamma_j:\mathcal V\cap\mathcal A^\infty\longrightarrow E_j
\]
can be chosen smoothly and
\[
 H_M[h](q,s)=h(s)q^2+\sum_{j=4}^M \gamma_j(h)(s)q^j
 \quad\hbox{satisfies}\quad P_{H_M[h]}(x)=x+O(|x|^{M+2}).
\]
Smoothness is understood in the local seed charts of
Proposition~\ref{v4:formal-transversality}; in particular it holds along
every smooth finite-dimensional family of seeds in this neighborhood.

For each finite $k$ and compact initial-value set $K$, there is a finite
integer $r\ge r_0$ such that, on
$\{h\in\mathcal V\cap\mathcal A^\infty:\|h\|_{C^r}\le R\}$, the
coefficients have uniform $C^k$ bounds. The scaled controls
\eqref{v4:formal-scaled-control} have a common small-amplitude range
and satisfy
\begin{equation}\label{v43:formal-uniform-residual}
 \|P_{H^{[M]}_{\varepsilon,h}}-\Id\|_{C^k(K)}
 \le C_{M,k,K,R}|\varepsilon|^{M+1}.
\end{equation}
The neighborhood and the spaces $E_j$ may depend on $M$.
\end{proposition}
\begin{proof}
Choose a basis of $\mathscr P_m$ for each even $m=4,6,\ldots,M-1$.
Theorem~\ref{v4:formal-curl-surjectivity} supplies $m+1$ real
trigonometric polynomials $k_{m,1},\ldots,k_{m,m+1}$ whose images under
$k\mapsto4m\mathcal A_m k$ at $h^\circ$ form a basis.
Let $A_m(h)$ be their coordinate
matrix in the chosen target basis. Its entries depend linearly and
continuously on $c_h,d_h$, and hence smoothly on $h$.
After shrinking a finite $C^{r_0}$ neighborhood, all the finitely many
determinants satisfy
\[
 |\det A_m(h)|\ge\tfrac12|\det A_m(h^\circ)|>0.
\]
Their inverses are smooth functions of their entries and are uniformly
bounded in this neighborhood.

Suppose that the coefficients below degree $m$ have been selected.
The finite jet recursion \eqref{v4:formal-jet-recursion} expresses the
current return coefficient $B_{m+1}$ by finitely many products and
time integrals of $h$ and those selected coefficients. It therefore
depends smoothly on the seed. If $b_m(h)$ is the coordinate vector
of $-\Curl B_{m+1}$, set
\[
 \gamma_m(h)=\sum_{\nu=1}^{m+1}
       \bigl(A_m(h)^{-1}b_m(h)\bigr)_\nu k_{m,\nu}.
\]
The remaining vector coefficient is curl-free by
Theorem~\ref{v4:formal-mixed-curl}. Let $\psi\in\mathscr P_{m+2}$
be its potential from Lemma~\ref{v4:formal-gradient-potential}.
The inverse moment map of Lemma~\ref{v4:formal-ridge-isomorphism}
selects $\gamma_{m+1}(h)$ smoothly in a fixed finite Fourier space so that
\[
 \int_0^T \gamma_{m+1}(h)(s)\ell(x,s)^{m+2}\,\dd s=-(m+2)\psi(x).
\]
By \eqref{v4:formal-direct-moment}, its direct return contribution
is $-\nabla\psi$.
The parity step gives the next vanishing degree as in
Theorem~\ref{v4:formal-formal-identity}. Induction through the
finitely many values of $m$ proves the assertion, with
$E_m=\operatorname{span}\{k_{m,\nu}\}$ and with the indicated fixed
Fourier spaces for $E_{m+1}$.

At each prescribed finite level, this construction uses only finitely
many products, time integrals, and inverses of matrices whose
determinants have the displayed lower bounds. Boundedness of the
corresponding finitely many seed derivatives therefore bounds every
selected coefficient in $C^k$. On a compact velocity cylinder the
scaled controls and the finitely many derivatives needed in the flow
equations have uniform bounds. Gronwall's inequality supplies a common
full-period tube for small $|\varepsilon|$, as well as uniform bounds
for the required initial-value and amplitude $\varepsilon$ derivatives of the flow.
The return derivatives in $\varepsilon$ through order $M$ vanish by
the finite jet construction. Taylor's formula with integral remainder,
also after $k$ initial-value derivatives, gives
\eqref{v43:formal-uniform-residual}. A velocity cutoff may be chosen
with a common plateau containing this tube, exactly as above.
\end{proof}

This proposition concerns one finite preparation. Joint smooth selection
of completed exact-return controls is a separate question.

The role of the integer $M$ is limited and explicit. The local inverse and the
nonlinear estimates first determine a finite derivative budget and finite negative
powers of $|\varepsilon|$. One then chooses a single odd $M$ large enough to make
\eqref{v4:formal-residual-estimate} pay those costs. Finally $|\varepsilon|$ is
reduced to absorb the finitely many coefficient constants. This finite selection
provides the initial control for the local exact-return argument.
\section{A characteristic condition at a nonzero initial value}
\label{v4:char-section}

The preceding construction supplies arbitrarily accurate finite approximations for
all seeds in $\mathcal A^\infty$. We now select the seeds for which the local linear
problem has a single simple characteristic direction. This selection is made at a
nonzero initial value, and the resulting inverse will be used on a neighborhood
of that point.

\subsection{The general seed symbol}

For $h\in\mathcal A^\infty$, use the functions $f_h,c,d$ of
\eqref{v4:seed-f-equation}--\eqref{v4:seed-cd-definition} and put
\begin{equation}\label{v4:char-vector}
 v_h(s)=c(s)m_s-d(s)e_s,\qquad
 g_{h,x}(s)=x\cdot v_h(s).
\end{equation}
The vector $v_h$ never vanishes, since $c>0$ and $e_s,m_s$ are orthonormal.
It is $\pi$-antiperiodic. In particular, $g_{h,x}(s+\pi)=-g_{h,x}(s)$;
its zero set is invariant under the half-period translation.

The first variation of the two branches is computed in
\cite[equations (5.3)--(5.9)]{Mao}. We express it in terms of $v_h$.
At the quadratic reference $H=\varepsilon h(s)q^2$, the flow and its
Jacobian have first variations
\begin{equation}\label{v4:char-first-variation}
 \Psi_s(x)=C_s\bigl(x+\varepsilon U_s(x)+O(\varepsilon^2)\bigr),\qquad
 \partial_sU_s=h(s)m_s(m_s\cdot x)^2,
\end{equation}
\[
 M_s=D_xU_s,
 \qquad
 \partial_s M_s=2h(s)(m_s\cdot x)m_s\otimes m_s,
 \qquad M_0=0.
\]
The forcing is antiperiodic. Its unique antiperiodic primitive $\mathsf K_s$ is a symmetric
matrix and satisfies $M_s=\mathsf K_s-\mathsf K_0$. Define
\[
 t_s=e_s^{\mathsf T}\mathsf K_sm_s,
 \qquad a_s=e_s^{\mathsf T}\mathsf K_se_s,
 \qquad b_s=m_s^{\mathsf T}\mathsf K_sm_s.
\]
Differentiation using $e_s'=m_s$, $m_s'=-e_s$ gives
\[
 a_s'=2t_s,\qquad t_s'=b_s-a_s,\qquad
 b_s'=-2t_s+2h(s)(m_s\cdot x).
\]
Consequently $t_s''+4t_s=2h(s)(m_s\cdot x)$.
By \eqref{v4:seed-cd-identity}, the unique antiperiodic solution is
\begin{equation}\label{v4:char-T}
 t_s=2\bigl(c(s)x\cdot m_s-d(s)x\cdot e_s\bigr)=2g_{h,x}(s).
\end{equation}
Uniqueness holds because the homogeneous solutions $\cos2s,\sin2s$ are periodic
under translation by $\pi$, whereas $t_s$ is antiperiodic.

The two incidence branches must be compared at matching conormal directions,
not merely at times differing by $\pi$. This distinction contributes to the
coefficient. Write the first-order covectors as
\[
 \alpha_s=e_s+\varepsilon(\mathsf K_s-\mathsf K_0)e_s+O(\varepsilon^2),\qquad
 \beta_s=m_s+\varepsilon(\mathsf K_s-\mathsf K_0)m_s+O(\varepsilon^2).
\]
At the opposite branch write $s_-=s+\pi+\varepsilon\delta+O(\varepsilon^2)$.
The expansion
\[
 \beta_{s_-}=-m_s+\varepsilon\{\delta e_s+(\mathsf K_s+\mathsf K_0)m_s\}
                 +O(\varepsilon^2)
\]
and matching its direction with $-\beta_s$ give $\delta=-2t_s$ by comparison
of the $e_s$ components. Similarly,
\[
 \alpha_{s_-}=-e_s+\varepsilon\{-\delta m_s+(\mathsf K_s+\mathsf K_0)e_s\}
                 +O(\varepsilon^2).
\]
As $\mathsf K_s$ is symmetric, $\det(e_s,m_s)=1$, and the $\mathsf K_0$ terms cancel,
\begin{equation}\label{v4:char-determinant}
 \det(\alpha_s,\alpha_{s_-})
 =\varepsilon(-\delta+2t_s)+O(\varepsilon^2)
 =8\varepsilon g_{h,x}(s)+O(\varepsilon^2).
\end{equation}
The first parameter derivative of the return is zero by the antiperiodicity in
\eqref{v4:char-first-variation}. Thus endpoint periodization does not change
this first variation at the reference.
The determinant calculation defines the homogeneous seed model
\begin{equation}\label{v4:char-principal-model}
 p_h^0(x,\xi)=8|\xi|g_{h,x}(s),\qquad \xi=|\xi|m_s.
\end{equation}
We study its characteristic directions here. The complete operator
construction identifies this model with the principal symbol after column
normalization and multiplication of the scalar row by $-\Lambda$;
the precise identification is
\eqref{v4:lin-schur}--\eqref{v41:lin-schur-seed}.
A nonzero real elliptic scalar multiple has the same characteristic
condition.

\subsection{Admissible nonzero initial values}

\begin{definition}[Seeds with one simple characteristic line]
\label{v4:seed-class-B}
A pair $(h,x_*)$, where $h\in\mathcal A^\infty$ and $x_*\ne0$, is admissible
if there is $s_*$ such that
\begin{equation}\label{v4:char-simple-zero}
 g_{h,x_*}(s)=0\ \Longleftrightarrow\ s-s_*\in\pi\mathbb Z,
 \qquad g'_{h,x_*}(s_*)\ne0.
\end{equation}
Let $\mathcal B^\infty$ be the set of seeds admitting at least one such $x_*$.
The zero condition is taken modulo $\pi$; the function $g_{h,x_*}$ itself is
antiperiodic, not a function on $\R/\pi\mathbb Z$.
\end{definition}

This is a condition on the full circle of cotangent directions. It prevents an
additional characteristic direction from being left outside a chosen local cone.
It is sufficient for the single-transport construction used here; its failure is
not asserted to be an obstruction to every possible closure argument.

\begin{proposition}[Openness of the characteristic condition]
\label{v4:char-openness}
For each admissible pair $(h^\circ,x_*^\circ)$, the same condition holds for all
$(h,x_*)$ sufficiently close in a fixed finite $C^r$ norm for $h$ and the usual
norm for $x_*$, with $h\in\mathcal A^\infty$. One may take $r\ge2$.
The characteristic roots, their transverse derivative lower bound, and the
noncharacteristic bound away from them may be chosen on one common neighborhood.
\end{proposition}
\begin{proof}
The Fourier multipliers in \eqref{v4:seed-cd-fourier} show that
$h\mapsto(c,d)$ is continuous from $C^r$ to $C^1$ for $r\ge2$.
For example, twice integrating the Fourier coefficient of $h$ by parts gives
$|\widehat h_k|\le C\|h\|_{C^2}(1+|k|)^{-2}$; after one derivative the
$c$ and $d$ series remain absolutely summable. The zero modes are bounded directly.
Thus $(h,x)\mapsto g_{h,x}$ is continuous into $C^1$.

Choose two disjoint small intervals about the two antipodal roots in a full
period. At the reference, the derivative has a fixed sign and magnitude bounded
below on each interval. The values at the two ends of an interval have opposite
signs. On the compact complement, $|g|$ has a positive minimum.
Small $C^1$ perturbations preserve all these strict inequalities. The intermediate
value theorem and strict monotonicity give one and only one zero in each interval;
there are none elsewhere. Antiperiodicity pairs the two roots, and the derivative
lower bound proves simplicity. The implicit function theorem gives their smooth
dependence on any smooth finite-dimensional set of parameters.
\end{proof}

\begin{theorem}[Factorization on the entire cotangent circle]
\label{v4:characteristic-factorization}
Fix an admissible pair $(h^\circ,x_*^\circ)$. After a fixed rotation of base
coordinates, there are a neighborhood $U$ of $x_*^\circ$, with
$0\notin\overline U$, and a neighborhood of $h^\circ$ in a finite $C^r$ norm
such that
\begin{equation}\label{v4:char-factorization}
 p_h^0(x,\xi)=a_h(x,\xi)\bigl(\xi_1-\beta_h(x)\xi_2\bigr),
 \qquad x\in U,\quad \xi\ne0.
\end{equation}
Here $a_h$ is real, smooth away from $\xi=0$, homogeneous of degree zero, and
bounded away from zero on the entire unit cotangent circle; $\beta_h$ is real
and smooth. The conclusion also holds for every smooth real, odd,
degree-one homogeneous perturbation sufficiently close to $p_h^0$ in
$C^1(\overline U\times S^1)$, after taking $U$ with compact closure
in the same base neighborhood.
\end{theorem}
\begin{proof}
Rotate once so that the reference characteristic line is $\xi_1=0$.
Proposition~\ref{v4:char-openness}, applied also to the base variable $x$, gives
exactly two simple antipodal roots on each cotangent circle. In the chart
$\xi_2>0$ the implicit function theorem writes the root as
$\xi_1=\beta_h(x)\xi_2$; oddness of $p_h^0$ under $\xi\mapsto-\xi$
gives the same equation in the opposite chart.

Near this line put $t=\xi_1/\xi_2$ and divide using
\[
 p_h^0(x,(t,1))
 =(t-\beta_h(x))\int_0^1
 \partial_{\xi_1}p_h^0\bigl(x,(\beta_h(x)+\tau(t-\beta_h(x)),1)\bigr)
 \,\dd\tau.
\]
The integral is smooth and bounded away from zero on a fixed narrow chart,
by the transverse derivative bound. Homogeneity gives a degree-zero quotient.
In the complementary directions both the linear factor and $p_h^0$ are nonzero,
so ordinary division gives the same quotient. Compactness gives a positive
lower bound for its absolute value on the whole circle. Oddness makes the
quotient even, so the two opposite charts agree.

The same argument applies to a real odd perturbation close enough in the finite
$C^1$ symbol seminorms used above. Differentiating the displayed integral formula
and the implicit root equation shows that each higher seminorm requires only
finitely many higher coefficient derivatives, with the same fixed derivative
offset. Product, reciprocal, and composition estimates retain at most one high
coefficient norm. Thus the factorization can be used on the fixed low-norm
neighborhood of the later operator estimates.
\end{proof}

The vector field $\partial_{x^1}-\beta_h(x)\partial_{x^2}$ never vanishes.
A flow box straightens it on a smaller neighborhood. The real elliptic factor
can then be removed by a finite elliptic reduction, leaving a transport operator
with a complete order-zero remainder. This reduction concerns the full symbol
near the characteristic line; an exact inverse still has to account for the
lower-order and finite-frequency terms.

\subsection{An angular criterion for general seeds}

The characteristic condition can be tested through the direction of $v_h$.
Since $c>0$, a smooth real lift of its angle is
\begin{equation}\label{v41:angular-lift}
 \vartheta_h(s)=s+\frac\pi2+\arctan\frac{d(s)}{c(s)},\qquad
 v_h(s)=|v_h(s)|(\cos\vartheta_h(s),\sin\vartheta_h(s)).
\end{equation}
The functions $c,d$ are $\pi$-periodic, so
$\vartheta_h(s+\pi)=\vartheta_h(s)+\pi$.
Differentiating $v_h=c m_s-d e_s$ gives
\begin{equation}\label{v41:angular-derivative}
 v_h'=(c'-d)m_s-(c+d')e_s,\qquad
 \vartheta_h'=
 \frac{\det(v_h,v_h')}{|v_h|^2}
 =\frac{c^2+d^2+c d'-d c'}{c^2+d^2}.
\end{equation}
Thus positivity of the numerator makes the direction map a diffeomorphism
of the projective circle. In that case every nonzero $x_*$ is admissible.
The following criterion also permits intervals where the angle decreases.

\begin{proposition}[Total angular variation]\label{v41:angular-criterion}
For $h\in\mathcal A^\infty$, define
\begin{equation}\label{v41:angular-variation}
 V(h)=\int_0^{2\pi}|\vartheta_h'(s)|\,\dd s.
\end{equation}
If $V(h)<6\pi$, then $h\in\mathcal B^\infty$.
More precisely, the set $E_h\subset\R/\pi\mathbb Z$ of directions $a$
for which $(h,(\cos a,\sin a))$ is admissible satisfies
\begin{equation}\label{v41:good-direction-measure}
 |E_h|\ge\frac{6\pi-V(h)}4.
\end{equation}
Here $|E_h|$ denotes ordinary angular measure on the circle of length $\pi$.
\end{proposition}
\begin{proof}
Put $x_a=(\cos a,\sin a)$. The zeros of $x_a\cdot v_h(s)$ are
exactly the inverse images of $a+\pi/2$ under the map
$\vartheta_h:\R/\pi\mathbb Z\to\R/\pi\mathbb Z$.
By \eqref{v41:angular-lift}, this map has degree one.
For a regular value its inverse images are finite, and their oriented
count is one. Their unoriented count on a half-period is therefore
$1,3,5,\ldots$. On a full period the number of zeros consequently obeys
\[
 Z_h(a):=\#\{s\bmod2\pi:x_a\cdot v_h(s)=0\}
       \in\{2,6,10,\ldots\}.
\]
At these zeros the derivative has absolute value
$|v_h(s)|\,|\vartheta_h'(s)|>0$, so all of them are simple.

The critical values have measure zero. Indeed, on an open neighborhood
of the compact critical set one can make $|\vartheta_h'|<\delta$;
the sum of the lengths of the images of its component intervals is
at most $\delta\pi$, using local lifts. Letting $\delta$ decrease to zero
proves the assertion. On every component of
$\{\vartheta_h'\ne0\}$ the angle is strictly monotone.
Change variables on these intervals, splitting at multiples of $\pi$
when necessary, and sum the nonnegative integrals. This gives the
one-dimensional counting formula
\[
 V(h)=\int_0^\pi Z_h(a)\,\dd a.
\]
Regular directions in $E_h$ have two zeros, while the other regular
directions have at least six. The exceptional critical values do not
contribute to the integral. Hence
\[
 V(h)\ge2|E_h|+6(\pi-|E_h|),
\]
which proves \eqref{v41:good-direction-measure} and the conclusion.
\end{proof}

Since $\int_0^{2\pi}\vartheta_h'\,\dd s=2\pi$, the sufficient condition
is equivalently
$\int_0^{2\pi}(\vartheta_h')_-\,\dd s<2\pi$, where
$r_-:=\max\{-r,0\}$.
This criterion and the determinant test in \eqref{v41:angular-derivative}
provide separate ways to check the additional condition on a balanced seed.
The relation $\mathcal A^\infty=\mathcal B^\infty$ is left undetermined
here, and $V(h)<6\pi$ is used only as a sufficient criterion.

\subsection{Fourfold frequencies and all higher even single frequencies}

\begin{proposition}[The fourfold reference]\label{v4:char-four}
The seed $h_4(s)=1+\sqrt{70/13}\cos4s$ belongs to $\mathcal B^\infty$,
with admissible initial value $x_*=(1,0)$.
\end{proposition}
\begin{proof}
Put $t=6/\sqrt{910}$. The coefficient formulas give
\[
 c(s)=\frac{3-13t\cos4s}{9},\qquad
 d(s)=\frac{8t\sin4s}{9}.
\]
At $x_*=(1,0)$ the characteristic function is
$g(s)=-c(s)\sin s-d(s)\cos s$.
The fourfold polynomial of \cite[equation (5.10)]{Mao} is obtained by
writing $z=\sin^2s$ and using
$\sin4s=4\sin s\cos s(1-2\sin^2s)$ yields
\begin{equation}\label{v4:char-four-factor}
 g(s)=-\frac{\sin s}{9}
       \bigl(3+19t+8tz-40tz^2\bigr).
\end{equation}
The polynomial in parentheses is concave on $0\le z\le1$.
Its endpoint values are $3+19t$ and $3-13t$, both strictly positive:
$13t=78/\sqrt{910}<3$ follows by squaring positive quantities.
Thus the only roots are $s\in\pi\mathbb Z$, and
$g'(0)=-(3+19t)/9\ne0$.
\end{proof}

\begin{remark}[Equality in the angular measure bound]
\label{v42:fourfold-measure-remark}
For $h_4$, the measure bound \eqref{v41:good-direction-measure} is attained.
Here is an exact computation. Put $t=6/\sqrt{910}$, $L=\pi/2$, and
$\phi(s)=\vartheta_{h_4}(s)-\pi/2=s+\arctan(d(s)/c(s))$.
The quarter-period symmetry gives
\[
 \phi(s+L)=\phi(s)+L,\qquad \phi(L-s)=L-\phi(s).
\]
For $u=\cos4s$, direct expansion yields
\[
 \det(v_{h_4},v_{h_4}')
 =\frac{70u^2+2\sqrt{910}\,u-83}{1365}.
\]
The numerator is negative at $u=-1,0$ and positive at $u=1$.
Its only root in $[-1,1]$ is the number $\rho$ below. Define
\begin{equation}\label{v42:fourfold-critical-values}
 \begin{aligned}
 \rho&=\frac{8\sqrt{105}-\sqrt{910}}{70},
 &s_c&=\frac14\arccos\rho,\\
 b&=\phi(s_c)
   =s_c+\arctan\frac{8t\sqrt{1-\rho^2}}{3-13t\rho}.
 \end{aligned}
\end{equation}
Thus $0<s_c<L/2$, and $\phi$ is strictly increasing, decreasing, and
increasing on $(0,s_c)$, $(s_c,L-s_c)$, and $(L-s_c,L)$, respectively.
The two critical values are $b$ and $L-b$.

To count the inverse images, we also check that these branches remain in
the same angular interval. Write
$\mathfrak p(z)=3+19t+8tz-40tz^2$, which is positive on $[0,1]$ by the
preceding proof. Equation~\eqref{v4:char-four-factor} and the
quarter-period symmetry give
\[
 v_{h_4}^1(s)=-\frac{\sin s}{9}\mathfrak p(\sin^2s),\qquad
 v_{h_4}^2(s)=\frac{\cos s}{9}\mathfrak p(\cos^2s).
\]
Hence $v_{h_4}^1<0<v_{h_4}^2$ for $0<s<L$, and $0<\phi(s)<L$ there.
Together with the strict decrease between the critical points, this gives
$0<L-b<L/2<b<L$. The three monotone branches therefore have images
$(0,b)$, $(L-b,b)$, and $(L-b,L)$. A direction in
$(0,L-b)\cup(b,L)$ has one inverse image in $(0,L)$, whereas a direction
in $(L-b,b)$ has three. Translation by $L$ gives the same count in the
next interval, with no overlapping images from adjacent intervals.
Since $x_a\cdot v_{h_4}(s)=0$ is equivalent to $\phi(s)=a$ modulo $\pi$,
the full-period zero count is consequently two or six for almost every
direction, and all these zeros are simple.

Each angular interval of length $L$ contributes $2(L-b)$ to the good
directions, and each parameter interval of length $L$ contributes
$b+(2b-L)+b$ to the variation. Thus
\begin{equation}\label{v42:fourfold-measure-equality}
 V(h_4)=16b-2\pi,\qquad
 |E_{h_4}|=2\pi-4b=\frac{6\pi-V(h_4)}4.
\end{equation}
These exact expressions give $|E_{h_4}|\approx139.469664^\circ$ and
$V(h_4)/(2\pi)\approx1.450337$. The equality concerns the measure
estimate; optimality of the sufficient threshold $V<6\pi$ is a separate
question.
\end{remark}

\begin{proposition}[Even single frequencies at every nonzero initial value]
\label{v4:char-high-frequency}
For every even integer $n\ge6$, let
\[
 h_n(s)=1+\lambda_n\cos ns,\qquad
 \lambda_n=\sqrt{\frac{2(n^2-1)(n^2-9)}{3(n^2-3)}}.
\]
Then $(h_n,x_*)$ is admissible for every $x_*\ne0$.
In particular, together with Proposition~\ref{v4:char-four}, every even
single-frequency seed $h_n$, $n\ge4$, belongs to $\mathcal B^\infty$.
\end{proposition}
\begin{proof}
Membership in $\mathcal A^\infty$ was established in
Example~\ref{v4:seed-example-single-mode}. Put $t=n^2$ and write
\[
 c(s)=\frac13-A_n\cos ns,\qquad d(s)=B_n\sin ns,
 \qquad
 A_n^2=\frac{2(t-3)}{3(t-1)(t-9)},\quad
 B_n=\frac{2nA_n}{t-3}.
\]
For $u=\cos ns$, expansion of \eqref{v41:angular-derivative} gives
\begin{align}
 \det(v_{h_n},v_{h_n}')
 &=(A_n^2-B_n^2)u^2+\frac{nB_n-2A_n}{3}u
       +\frac19+B_n^2-nA_nB_n \notag\\
 &=(A_n^2-B_n^2)u^2+\frac{2A_n}{t-3}u
       +\frac19+B_n^2-2A_n^2-\frac{6A_n^2}{t-3}.
 \label{v41:single-angular-polynomial}
\end{align}
For $n\ge6$, one has $B_n\le A_n$, since $2n\le n^2-3$.
Moreover, the derivative of $(t-3)/((t-1)(t-9))$ has numerator
$-(t-3)^2-12$, so $A_n\le A_6$ for $t\ge36$.
The exact values and a rational upper bound are
\[
 A_6^2=\frac{22}{945},\qquad A_6<\frac2{13};
\]
the latter follows from $22\cdot169<4\cdot945$.
Discarding the nonnegative quadratic term and $B_n^2$ in
\eqref{v41:single-angular-polynomial}, and using $|u|\le1$, yields
\begin{align}
 \det(v_{h_n},v_{h_n}')
 &\ge\frac19-2A_6^2-\frac{6A_6^2+2A_6}{33}\notag\\
 &>\frac19-\frac{44}{945}-\frac{132/945+4/13}{33}
   =\frac{2297}{45045}>\frac1{20}.
 \label{v4:char-high-margin}
\end{align}
Thus $\vartheta_{h_n}'>0$ everywhere. Its half-period increment is $\pi$
by \eqref{v41:angular-lift}, so its map on the projective circle is a
diffeomorphism. For any $x_*\ne0$, $x_*\cdot v_{h_n}=0$ has exactly one
root modulo $\pi$, and at that root
\[
 |(x_*\cdot v_{h_n})'|
       =|x_*|\,|v_{h_n}|\,\vartheta_{h_n}'>0.
\]
This proves the required condition.
\end{proof}

The point $x_*$ in the proposition may be chosen anywhere away from zero.
The local inverse and the completed control will be constructed separately
near each chosen point. For $h_4$ the preceding polynomial argument supplies
an admissible direction without requiring its angle to be monotone everywhere.

\subsection{Local families of arbitrarily many seed parameters}

\begin{corollary}[Relative neighborhoods and finite-dimensional seed sets]
\label{v4:char-families}
Every seed in Propositions~\ref{v4:char-four} and
\ref{v4:char-high-frequency} has a relative neighborhood in
$\mathcal A^\infty$, in a fixed finite $C^r$ topology, contained in
$\mathcal B^\infty$. For each finite $J\ge1$, the fourfold coefficient ellipsoid
on frequencies $4,8,\ldots,4J$ has a relatively open subset of dimension $2J-1$
contained in $\mathcal B^\infty$.
\end{corollary}
\begin{proof}
The first assertion is Proposition~\ref{v4:char-openness}.
For the second, the coefficient vector of $h_4$ is
$(a_1,b_1,a_2,b_2,\ldots)=(\sqrt{70/13},0,0,0,\ldots)$.
The derivative of the ellipsoid equation with respect to $a_1$ is nonzero there.
Thus it is a smooth $(2J-1)$-dimensional hypersurface near that vector.
All its members are in $\mathcal A^\infty$ by
Theorem~\ref{v4:seed-theorem-sparse}; sufficiently small coefficient
neighborhoods are small in every fixed finite $C^r$ norm.
Intersecting with the relative neighborhood from the first assertion gives the
claimed open subset.
\end{proof}

\begin{corollary}[Mixed four--ten seeds near the single-frequency endpoints]
\label{v43:char-four-ten-endpoints}
For the family \eqref{v4:seed-four-ten-family}, there are
$t_0,t_1\in(0,1/2)$ such that every member with
$0\le t<t_0$ or $1-t_1<t\le1$ belongs to $\mathcal B^\infty$,
for all phases $\theta_4,\theta_{10}$.
For $0\le t<t_0$, every nonzero initial value is admissible.
In particular, both endpoint neighborhoods contain genuinely mixed seeds.
\end{corollary}
\begin{proof}
All these seeds lie in $\mathcal A^\infty$ by
Example~\ref{v4:seed-example-four-ten}. At $t=0$ they are the
translates $1+\lambda_{10}\cos(10s+\theta_{10})$ of $h_{10}$.
The determinant bound \eqref{v4:char-high-margin} gives
$\det(v,v')>1/20$, uniformly in the angle and the phase.
Write $c_t,d_t$ for the coefficients of the mixed seed and
$c_0,d_0$ for those of its ten-frequency endpoint. Their finite
Fourier formulas give
\[
 \|c_t-c_0\|_{C^1}+\|d_t-d_0\|_{C^1}\le C\sqrt t,
\]
with $C$ independent of the two phases. These functions and their
first derivatives are uniformly bounded over the parameter set.
Since the determinant is $c^2+d^2+cd'-dc'$, its change is also at
most $C'\sqrt t$. Choose $t_0>0$ so that $C'\sqrt{t_0}<1/40$.
Then the determinant is greater than $1/40$ throughout this
neighborhood. The angle map is a diffeomorphism of the projective
circle, proving admissibility for every nonzero initial value.

At the other endpoint, translate the variable by $-\theta_4/4$.
The translated seed tends to $h_4$ in every finite $C^r$ norm as
$t\uparrow1$, uniformly in both phases: the fourth coefficient
changes by $O(1-t)$ and the tenth coefficient is $O(\sqrt{1-t})$.
Proposition~\ref{v4:char-openness}, applied at $(h_4,(1,0))$,
therefore supplies one $t_1>0$ for all the translated seeds.
Translation of a seed translates the angular parameter and applies
a constant planar rotation to $v_h$, as follows from its definition in
\eqref{v4:char-vector}. It preserves the existence of an admissible
nonzero initial value. Translating back proves the assertion near
$t=1$.
\end{proof}

This corollary gives open sets of mixed examples from the full
coefficient ellipsoid of Section~\ref{v4:seed-sec-fourfold}.
The intervening parameters remain candidates for the general
characteristic tests; no conclusion about the whole ellipsoid is
needed for these local families.

\begin{proposition}[A coupled family on frequencies two, four, and six]
\label{v44:coupled-realization}
Write $b_0=\sqrt{70/13}$ and
\[
 h_{a,b,t}(s)=1+a\cos 2s+b\cos4s+t\cos6s.
\]
For all $t$ in a neighborhood of $0$ there are smooth real functions
$a(t),b(t)$ with $a(0)=0$, $b(0)=b_0$, for which
$(h_{a(t),b(t),t},(1,0))$ is admissible. All three nonconstant coefficients
are nonzero for sufficiently small $t\ne0$.
There is also an admissible member with $t=1$ and $229/100<b<23/10$.
\end{proposition}
\begin{proof}
Substitution in \eqref{v4:seed-Q} and \eqref{v4:seed-S} gives
\[
 Q=\frac13+\frac{a^2}{30}-\frac{13b^2}{210}-\frac{11t^2}{630},
 \qquad S=-\frac a3+\frac{ab}{10}+\frac{bt}{42}.
\]
The three summands in $S$ come respectively from the pairs of frequencies
$(0,2)$, $(-2,4)$ and $(-4,6)$ and their reversals. At $(0,b_0,0)$
the two diagonal entries of the Jacobian of $(S,Q)$ with respect to
$(a,b)$ are $-1/3+b_0/10$ and $-13b_0/105$, both nonzero,
and the two off-diagonal entries vanish. The implicit function theorem
gives $a(t),b(t)$ and
\[
 a'(0)=-\frac{b_0}{42(-1/3+b_0/10)}\ne0.
\]
Positivity of $Lh$ and the single simple zero for $x_*=(1,0)$ persist by
Propositions~\ref{v4:char-four} and~\ref{v4:char-openness}.

For $t=1$, put $a=5b/(70-21b)$ to make $S=0$. Then $630Q$ equals
$199+21a^2-39b^2$. This expression is positive at $b=229/100$ and
negative at $b=23/10$. At any intervening zero, $(a,b)$ lies in the rectangle
\[
 [52/100,53/100]\times[229/100,23/10].
\]
Put $z=\cos2s$. The Fourier multipliers in \eqref{v4:seed-cd-fourier}
give $c(s)=C(z)$ and $g_{h,(1,0)}(s)=-\sin s\,P(z)$, where
\[
 \begin{aligned}
 C(z)&=\tfrac13+\tfrac{a}{15}z-\tfrac{13b}{105}(2z^2-1)
                         -\tfrac{11}{315}(4z^3-3z),\\
 P(z)&=C(z)+(1+z)
       \left[-\tfrac{4a}{15}+\tfrac{16b}{105}z+
                                \tfrac4{315}(4z^2-1)\right].
 \end{aligned}
\]
To check positivity on $[-1,1]$, write each cubic $R(z)$ at $z=2u-1$ as
$\sum_{j=0}^3 B_j\binom3j u^j(1-u)^{3-j}$. Each $B_j$ is affine in
$(a,b)$, so its minimum on the rectangle is attained at a vertex.
The minima over the four vertices, multiplied by $94500$, are
\[
 \begin{array}{c|rrrr}
 R&B_0&B_1&B_2&B_3\\ \hline
 C&4551&58542&93747&4566\\
 P&4551&30054&67960&51186
 \end{array}
\]
For instance $B_0=R(-1)$ and $B_1=R(-1)+2R'(-1)/3$.
The other two coefficients are $B_2=R(1)-2R'(1)/3$ and $B_3=R(1)$.
The four Bernstein basis functions are nonnegative and sum to one.
It follows that $C,P>0$ throughout $[-1,1]$. Thus $g_{h,(1,0)}$
has exactly one zero modulo $\pi$, and its derivative at zero is
$-P(1)\ne0$.
\end{proof}

The dimensions in Corollary~\ref{v4:char-families} count input seed
parameters; a geometric separation statement
for the completed metrics is proved using the invariant in
Section~\ref{v4:geom-tensors}. The exact-closure construction realizes each
member separately, with local bounds near a fixed admissible pair. Smooth
selection of all completed metrics, and a corresponding moduli-space dimension,
are separate questions.
\part{A local inverse near a nonzero initial value}\label{part:analysis}
\section{The actual return differential and a periodic comparison}
\label{v4:endpoint-chapter}

The equation to be solved is an equation for the original oscillator flow.
We first calculate its differential. We then construct a periodic comparison
family whose linearized integral admits the two-branch analysis in the next
section. The discrepancy between the two integrals is retained explicitly.

\subsection{Fixed flow tubes and actual control variations}
Write $T=2\pi$. Fix a bounded initial-value neighborhood $O$, a compact
subset $K\Subset O$, and a control cylinder
$\mathcal W=(-V,V)\times\R/(T\mathbb Z)$, with the arguments of a control
written as $H(q,s)$. Choose $V$ larger than every harmonic velocity sampled
from a neighborhood of $\overline O$. A sufficiently small $C^1$ neighborhood
of zero, with any further finite low norms bounded as specified below, has a common flow tube
containing $[0,T]\times\overline O$. To see this directly, write the actual
integral equation relative to the harmonic rotation $C_s$:
\[
 \Psi_H(x,s)=C_sx+
   \int_0^s C_{s-t}e_2 H(q_H(x,t),t)\,\dd t.
\]
The bound $|\Psi_H-C_sx|\le T\|H\|_{C^0}$, applied up to a first possible
exit time, prevents exit from a fixed larger tube when the low ball is small.
Local uniqueness then continues the same solution throughout the period.
The variational integral equations give joint smoothness in $(x,s)$ and in
each finite collection of smooth control parameters on this same tube.
Only a fixed interior velocity collar is sampled in all subsequent estimates.

Put $A_H=D_x\Psi_H$, $P_H=\Psi_H(\cdot,T)$, and
$N(H)=P_H-\Id$. Differentiating the oscillator gives
\begin{equation}\label{v4:endpoint-variational}
 \partial_s A_H=
 \begin{pmatrix}0&1\\-1&H_q(q_H,s)\end{pmatrix}A_H,
 \qquad A_H(x,0)=\Id.
\end{equation}
In particular,
\[
 \det A_H(x,s)=\exp\left(\int_0^s H_q(q_H(x,t),t)\,\dd t\right)>0.
\]
The bounds for this determinant and its reciprocal are uniform on the
chosen low ball.

\begin{lemma}[Differential of the actual return]\label{v4:actual-return-differential}
For a real smooth periodic variation $\eta$, define
\begin{equation}\label{v4:true-source-operator}
 \mathcal M_H^{\rm tr}\eta(x)=
 \int_0^T A_H(x,s)^{-1}e_2\,
                 \eta(q_H(x,s),s)\,\dd s.
\end{equation}
Then
\begin{equation}\label{v4:actual-derivative}
 DN_H[\eta](x)=DP_H(x)\,\mathcal M_H^{\rm tr}\eta(x).
\end{equation}
\end{lemma}
\begin{proof}
Let $Z=\partial_\tau\Psi_{H+\tau\eta}|_{\tau=0}$. Its initial value is zero
and its equation is
\[
 Z_s=\begin{pmatrix}0&1\\-1&H_q(q_H,s)\end{pmatrix}Z
             +e_2\eta(q_H,s).
\]
Variation of constants gives $Z(T)=A_H(T)\mathcal M_H^{\rm tr}\eta$.
Since $A_H(T)=DP_H$, this proves the identity. In particular the endpoint
Jacobian is part of the actual differential, even when the return is close
to the identity.
\end{proof}

For later use the second variation in the same direction is equally concrete.
Writing $W=\partial_\tau^2\Psi_{H+\tau\eta}|_{0}$ and denoting the second
component of $Z$ by $Z_q$, it satisfies
\begin{equation}\label{v4:second-variation}
 W_s=\begin{pmatrix}0&1\\-1&H_q(q_H,s)\end{pmatrix}W
   +e_2\{H_{qq}(q_H,s)Z_q^2+2\eta_q(q_H,s)Z_q\},\qquad W(0)=0.
\end{equation}
Both variations are obtained from the same perturbed flow.

\subsection{Bounds with a fixed derivative offset}
All norms in this subsection are taken on the chosen flow and sampling
collars. A low output index and a low background index play different
roles; we keep this distinction in the following statement.

\begin{lemma}[Uniform flow and variation estimates]\label{v43:endpoint-tame}
Fix an integer $b\ge2$ and a bound $M$ for
$\|H\|_{C^{b+4}(\mathcal W)}$, in addition to the $C^1$ smallness
ensuring the common flow tube. The return map is twice differentiable
between the finite control and target spaces appearing below. There is
an integer $a\ge4$, independent of $k$, such that, for every $k\ge b$,
\begin{align}
 \|DN_H\eta\|_{C^k}
 &\le C_k\{\|\eta\|_{C^{k+a}}
       +(1+\|H\|_{C^{k+a}})\|\eta\|_{C^{b+a}}\},
       \label{v4:endpoint-first-bound}\\
 \|D^2N_H[\eta,\zeta]\|_{C^k}
 &\le C_k\{\|\eta\|_{C^{k+a}}\|\zeta\|_{C^{b+a}}
       +\|\zeta\|_{C^{k+a}}\|\eta\|_{C^{b+a}}\notag\\
 &\hspace{38mm}+(1+\|H\|_{C^{k+a}})
                  \|\eta\|_{C^{b+a}}\|\zeta\|_{C^{b+a}}\}.
       \label{v4:endpoint-second-bound}
\end{align}
The constants depend on $b,M,k$ and the collars, but are uniform over
these backgrounds and independent of the variation directions.
One may take $a=4$ in the two displayed estimates. The flow,
its inverse Jacobian and their first two control variations have the
same form of estimate, with a further fixed derivative offset when
needed. None of these offsets depends on $k$.
\end{lemma}
\begin{proof}
At target order $k$, consider the integral equation as an equation in
$C([0,T],C^k(\overline O;\R^2))$, after slightly reducing $O$ if
necessary. Its derivative with respect to the trajectory is a
Volterra operator. If its coefficient has multiplication norm $M_k$,
its $j$th iterate has norm at most $(TM_k)^j/j!$.
The resulting convergent inverse series and the ordinary Taylor
formula for $H(q,s)$ prove twice differentiability when the control
has $k+2$ derivatives. They also justify the first and second
variational equations, including the mixed forcing
\[
 H_{qq}Z_{\eta,q}Z_{\zeta,q}
       +\eta_qZ_{\zeta,q}+\zeta_qZ_{\eta,q}.
\]
Here each coefficient is evaluated along the unperturbed flow.

We next obtain bounds with only one high norm. Put $B=b+4$.
For nonzero directions divide $\eta,\zeta$ by their $C^B$ norms,
and denote the normalized functions by $\nu_1,\nu_2$. Adjoin constant
parameters $a_1,a_2$ to the state, with vector field
\[
 \mathcal V(p,q,a,s)
 =\bigl(q,-p+H(q,s)+a_1\nu_1(q,s)+a_2\nu_2(q,s),0,0\bigr).
\]
All its flow derivatives below are evaluated at $a=0$. For an
integer $K\ge B$, set
$Q_K=1+\|H\|_{C^K}+\sum_{j=1}^2\|\nu_j\|_{C^{K-1}}$.
The only nonzero derivatives of a term $a_j\nu_j$ at $a=0$
include one $a_j$ derivative. Interpolation between the bounded
low norms and these high norms therefore gives
\[
 \|D^r\mathcal V\|\le C_KQ_K^{(r-1)/(K-1)}
       \quad(1\le r\le K).
\]
In particular $D\mathcal V$ has a uniform low bound.
If $J_j$ is a total derivative of degree $j$ in the initial state
and the two parameters, its highest term is $D\mathcal V\,J_j$.
Every other term has the form
$D^r\mathcal V[J_{i_1},\ldots,J_{i_r}]$, with
$r\ge2$, $i_\nu\ge1$ and $\sum i_\nu=j$. Its interpolation exponent is
\[
 \frac{r-1}{K-1}+\sum_{\nu=1}^r\frac{i_\nu-1}{K-1}
       =\frac{j-1}{K-1}.
\]
Induction and variation of constants give
$\|J_j\|\le C_{K,j}Q_K^{(j-1)/(K-1)}$ for $j\le K$;
the propagator in this induction depends only on the low bound for
$D\mathcal V$. Take $K=\max\{B,k+2\}$ and undo the direction
normalizations. Since $K\le k+4$, this gives
\eqref{v4:endpoint-first-bound}--\eqref{v4:endpoint-second-bound}
with $a=4$, after evaluation at $s=T$. A zero direction contributes
zero by linearity or bilinearity.

Spatial derivatives of $A_H$ use one further spatial derivative
of the flow. Differentiating $A_H^{-1}A_H=\Id$ isolates the highest
inverse derivative; the remaining terms partition the total
derivative degree, and the undifferentiated inverse matrix is
uniformly bounded by~\eqref{v4:endpoint-variational}. The same
interpolation and induction therefore give a fixed extra offset
for the inverse Jacobian. Time derivatives follow by differentiating
the ODE, with one fewer derivative on its right side. These arguments
also give the stated joint flow estimates on the whole time interval.
\end{proof}

In particular, to use a pure low-output estimate at order $b$ we
bound the background in $C^{b+4}$ once. Later the common Sobolev
index $b_0$ is chosen to dominate this finite requirement and the
other low indices; it is not increased with a high target order.

\subsection{Periodization on domains with collars}
Choose $\chi\in C^\infty([0,T])$ with $0\le\chi\le1$, equal to zero
near $0$ and to one near $T$.
For $0\le t\le1$ set
\begin{equation}\label{v4:periodization-definition}
 E_{s,t}=\Id+t\chi(s)N(H),\qquad Q_{s,t}=E_{s,t}^{-1},
 \qquad Q_s=Q_{s,1},\qquad
 \widetilde\Psi_H(x,s)=\Psi_H(Q_s(x),s).
\end{equation}
These are local inverse maps on nested convex initial-value domains.
There is no assumption that they are globally invertible on $\R^2$.

\begin{lemma}[Protection of the inverse domain]\label{v4:inverse-domain}
Let $U'\Subset U$ be concentric balls with boundary gap $d>0$.
If $\|N(H)\|_{C^0(U)}\le d/4$ and
$\|D_xN(H)\|_{C^0(U)}\le1/2$, then every $Q_{s,t}$ in
\eqref{v4:periodization-definition} is defined on $U'$ and maps $U'$ into $U$.
\end{lemma}
\begin{proof}
For $x\in U'$, the map $y\mapsto x-t\chi(s)N(H)(y)$ is a contraction
of the closed ball $\overline B(x,d/2)$ into itself. Its fixed point solves
$E_{s,t}(y)=x$ and belongs to $U$. The first derivative of $E_{s,t}$ is
invertible there. The inverse-function theorem identifies the fixed points
with a smooth inverse on $U'$. This establishes the domain inclusion before
any estimate is evaluated at $Q_{s,t}(x)$.
\end{proof}

\begin{lemma}[A smooth periodic comparison family]\label{v4:periodic-comparison}
On fixed nested domains satisfying Lemma~\ref{v4:inverse-domain},
$\widetilde\Psi_H$ has matching derivatives of every order at $0,T$.
It therefore defines a smooth family parametrized by $\R/T\mathbb Z$.
\end{lemma}
\begin{proof}
Time periodicity of the original equation and uniqueness give
$\Psi_H(x,T+t)=\Psi_H(P_H(x),t)$ wherever both sides are defined.
Near $s=0$ the comparison is $\Psi_H(x,s)$, while near $s=T$ it is
$\Psi_H(P_H^{-1}(x),s)=\Psi_H(x,s-T)$.
These equalities hold on neighborhoods, so all time and spatial derivatives
match. For the quantitative estimates put
$F_s(y)=y+t\chi(s)N(H)(y)$. The first inverse derivative is
$DQ_{s,t}=(DF_s\circ Q_{s,t})^{-1}$, with norm at most two
under Lemma~\ref{v4:inverse-domain}. At spatial derivative degree $j\ge2$,
differentiating $F_s\circ Q_{s,t}=\Id$ gives
\[
 (DF_s\circ Q_{s,t})D^jQ_{s,t}
 =-\sum_{\substack{r\ge2,\ i_\nu\ge1\\i_1+\cdots+i_r=j}}
       c_{r,i}\,(D^rF_s\circ Q_{s,t})
           [D^{i_1}Q_{s,t},\ldots,D^{i_r}Q_{s,t}],
\]
with fixed combinatorial coefficients. For time, interpolation, or
control parameters, apply this formula to the augmented inverse identity
$(y,a)\mapsto(F(y,a),a)$, where $a$ lists the finitely many parameters
being differentiated. Its inverse Jacobian is triangular and uniformly
bounded at the required low indices. This includes the explicit
parameter derivatives of $F$, as well as those of the inverse map.
The degree excess in each product is
$(r-1)+\sum(i_\nu-1)=j-1$. Lemma~\ref{v43:endpoint-tame},
the inverse margin, and interpolation bound the right side with one
high control norm and fixed offsets. Composition with $\Psi_H$
uses the same partition formula. Thus the inverse maps and the
periodic comparison have these estimates on the common collars.
\end{proof}

Write $\widetilde q_H$ for the second component of this family and put
\begin{equation}\label{v4:comparison-operator}
 \widetilde w_H=(D_x\widetilde\Psi_H)^{-1}e_2,
 \qquad
 \mathcal M_H\eta(x)=\int_0^T\widetilde w_H(x,s)
                       \eta(\widetilde q_H(x,s),s)\,\dd s.
\end{equation}
The family $\widetilde\Psi_H$ generally does not solve the oscillator.
In particular its periodicity is not the exact return that the paper seeks.

\subsection{The error between the two integral operators}
For a vector field $Z$ write $Q^*Z=(DQ)^{-1}Z\circ Q$. If
$Z_s=A_H^{-1}e_2\eta(q_H,s)$, the chain rule gives
$\mathcal M_H\eta=\int_0^TQ_s^*Z_s\,\dd s$.
The inverse identity yields
\[
 Q_s-\Id=-\chi N(H)\circ Q_s,
 \qquad (DQ_s)^{-1}=\Id+\chi D_xN(H)\circ Q_s.
\]
Consequently the exact pointwise difference is
\begin{align}
 Q_s^*Z_s-Z_s={}&\chi D_xN(H)(Q_s)Z_s(Q_s)\notag\\
 &+\int_0^1 DZ_s\bigl(x+t(Q_s(x)-x)\bigr)
                         (Q_s(x)-x)\,\dd t.
 \label{v4:pullback-difference}
\end{align}
Every summand contains the actual residual or one of its derivatives.

\begin{proposition}[Residual-controlled comparison]\label{v4:comparison-defect}
For convex $U'\Subset U$ on which the preceding inverses and segments remain
inside $U$, assume the $C^{b+4}$ background bound of
Lemma~\ref{v43:endpoint-tame}. Then the fixed-level estimate is
\begin{equation}\label{v4:comparison-low}
 \|DP_H(\mathcal M_H^{\rm tr}-\mathcal M_H)\eta\|_{C^b(U')}
 \le C_b\|N(H)\|_{C^{b+1}(U)}\|\eta\|_{C^{b+4}(\mathcal W)}.
\end{equation}
On a fixed low neighborhood, its high-level version has a fixed offset $a$:
\begin{align}
 \|DP_H(\mathcal M_H^{\rm tr}-\mathcal M_H)\eta\|_{C^k(U')}
 \le C_k\{&\|N(H)\|_{C^{k+a}(U)}\|\eta\|_{C^{b+a}}
       +\|N(H)\|_{C^{b+a}(U)}\|\eta\|_{C^{k+a}}\notag\\
       &+(1+\|H\|_{C^{k+a}})
                    \|N(H)\|_{C^{b+a}(U)}\|\eta\|_{C^{b+a}}\}.
 \label{v4:comparison-high}
\end{align}
\end{proposition}
\begin{proof}
Let $Z_s=A_H^{-1}e_2\eta(q_H,s)$. Lemma~\ref{v43:endpoint-tame}
at the required finite indices and the inverse-map recursion give
\[
 \sup_s\|Z_s\|_{C^{b+1}}\le C_b\|\eta\|_{C^{b+4}},\qquad
 \|Q_s-\Id\|_{C^b}\le C_b\|N(H)\|_{C^b},\qquad
 \|DP_H\|_{C^b}+\|Q_s\|_{C^{b+1}}\le C_b.
\]
For example the $b+1$ derivatives of $A_H^{-1}$ use at most
$b+2$ spatial derivatives of the flow, and are covered by the
assumed background bound. Substitute $Q_s-\Id=-\chi N(H)\circ Q_s$
in~\eqref{v4:pullback-difference}. The product and composition
rules at order $b$ bound the first summand by
$C_b\|N(H)\|_{C^{b+1}}\|\eta\|_{C^{b+4}}$.
The integral summand has the same bound, using the displayed
$C^{b+1}$ estimate for $Z_s$. Integration in $s$ and multiplication
by $DP_H$ prove~\eqref{v4:comparison-low}.

For the high estimate retain the explicit factor $N(H)$ or $DN(H)$
in each summand. Distribute the high differentiations among this
factor, the source $\eta$, and the flow or inverse-map coefficients.
The chain-rule partition and the product interpolation from
Lemma~\ref{v43:endpoint-tame} give one high factor. If it is the
residual or the source, this gives the first or second term in
\eqref{v4:comparison-high}; if it is a coefficient, the explicit
residual and source factors remain at the low level, giving the
third term. The one spatial derivative in $DZ_s$ and the inverse
Jacobian costs are absorbed in a single finite $a$. Convexity
and Lemma~\ref{v4:inverse-domain} justify all evaluations. At a
fixed order, shrinking collars contribute only polynomial costs
in their inverse gaps, as used in the later iteration.
\end{proof}

\subsection{The incidence geometry used by the complete kernels}
For the actual flow one has the identity
\[
 \det D_x(q_H,\partial_s q_H)
 =\det D_x(q_H,-p_H+H(q_H,s))=\det A_H>0.
\]
The corresponding assertion for the comparison family includes a residual
term. Indeed
\[
 \partial_sQ_s=-\chi'(s)DQ_s\,N(H)(Q_s),
 \quad
 \partial_s\widetilde q_H=(\partial_sq_H)\circ Q_s
                   +Dq_H(Q_s)\partial_sQ_s.
\]
Differentiating these formulas once in $x$, and using the bounded inverse
matrices, bounds the determinant error by $C\|N(H)\|_{C^2(U)}$ on a
smaller domain. It is therefore small on the common low ball.
The complete kernels require the corresponding statement for all
endpoint pairs, for which the following quantitative form is useful.

\begin{lemma}[The two branches for every endpoint pair]\label{v43:all-pairs-twist}
Let $O_c$ be a convex domain with compact closure inside the comparison
chart. Suppose
\[
 \sup_{\overline O_c\times[0,T]}
 \left|D_x\widetilde q_H-m_s\right|\le\delta,\qquad
 \sup_{\overline O_c\times[0,T]}
 \left|\partial_sD_x\widetilde q_H-m_s'\right|\le\delta,
 \qquad \delta\le\tfrac14.
\]
For all $x,y\in\overline O_c$, define
\begin{equation}\label{v4:averaged-normal}
 \overline n_s(x,y)=\int_0^1D_x\widetilde q_H(y+t(x-y),s)\,\dd t.
\end{equation}
Then
\[
 \tfrac34\le|\overline n_s|\le\tfrac54,\qquad
 \det(\overline n_s,\partial_s\overline n_s)\ge\tfrac7{16}.
\]
For each endpoint pair and each sign, the map
$(s,\sigma)\mapsto\sigma\overline n_s(x,y)$ is a diffeomorphism from
$\R/(2\pi\mathbb Z)\times(0,\infty)$, or respectively from
$\R/(2\pi\mathbb Z)\times(-\infty,0)$, onto $\R^2\setminus\{0\}$.
Its absolute Jacobian is
$|\sigma|\det(\overline n_s,\partial_s\overline n_s)$.
\end{lemma}
\begin{proof}
Convexity keeps every segment in $O_c$, so averaging preserves both
$\delta$ bounds. The norm bounds follow by the triangle inequality.
Since $|m_s|=|m_s'|=1$ and $\det(m_s,m_s')=1$,
\[
 \det(\overline n_s,\partial_s\overline n_s)
       \ge1-2\delta-\delta^2\ge\tfrac7{16}.
\]
The angle of $\overline n_s$ has derivative
$\det(\overline n_s,\partial_s\overline n_s)/|\overline n_s|^2
\ge7/25$. The straight homotopy
$m_s+t(\overline n_s-m_s)$ never vanishes, so its direction has
degree one. A strictly increasing angle with degree one is a
diffeomorphism of the circle. Thus each direction determines
exactly one $s$ on either branch, and the radius determines
$\sigma$ by $|\sigma|=|\xi|/|\overline n_s|$. The inverse function
theorem and the displayed Jacobian give the smooth inverse,
including smooth dependence on $(x,y)$.
\end{proof}

The hypotheses hold on a common neighborhood of $H=0$. Indeed,
the integral equation and its spatial derivative give
$D_x\Psi_H-C_s=O(\|H\|_{C^1})$ uniformly; the endpoint inverse
identity gives $DQ_s-\Id=O(\|H\|_{C^1})$. Differentiating the
comparison once in time and once in space uses the displayed
formula for $\partial_sQ_s$ and at most two spatial derivatives
of $N(H)$. The corresponding finite flow estimates therefore give
\[
 \|D_x\widetilde q_H-m_s\|_{C^0}
 +\|\partial_sD_x\widetilde q_H-m_s'\|_{C^0}
       \le C\|H\|_{C^3(\mathcal W)}
\]
on the same collars. Restricting this one finite low norm makes
$\delta\le1/4$. The fundamental theorem of calculus also gives
$\widetilde q_H(x,s)-\widetilde q_H(y,s)
=\overline n_s(x,y)\cdot(x-y)$, which is the phase identity used below.

\subsection{The identity used by the nonlinear iteration}
If a comparison right inverse $\mathcal B_H$ satisfies
$\mathcal M_H\mathcal B_Hf=f$ on an inner domain, put
$\mathcal Q_H=\mathcal B_H(DP_H)^{-1}$. Then, on that domain,
\begin{equation}\label{v4:approximate-newton-identity}
 DN_H\mathcal Q_H-\Id
 =DP_H(\mathcal M_H^{\rm tr}-\mathcal M_H)
                         \mathcal B_H(DP_H)^{-1}.
\end{equation}
This follows directly from Lemma~\ref{v4:actual-return-differential}.
The right side retains a residual factor. Substituting the residual itself
as data makes it quadratic in residual size, with a loss of derivatives.
Sections~\ref{v4:newton-chapter} and~\ref{v4:convergence-chapter} handle
that loss by smoothing and prove convergence for the actual map $N$.
\section{Complete branch operators and their normalization}\label{v4:lin-normalization-section}

We now prepare the linear equation used to close the oscillator.  Fix a
seed $h^\circ$ satisfying Definition~\ref{v4:seed-class-B}, and a nonzero
initial point $x_*$ witnessing its simple-characteristic condition.  All
domains and estimates in this section are chosen near this pair.
We allow $h$ to vary in a sufficiently small relative neighborhood of
$h^\circ$ among the balanced seeds.  The control is written
\begin{equation}\label{v4:lin-control-chart}
 H_{\varepsilon,h,z}(q,s)=\varepsilon\bigl(h(s)q^2+z(q,s)\bigr).
\end{equation}
Here the control correction $z$ ranges over real smooth $2\pi$-periodic
functions in a fixed low-norm ball. At $\varepsilon=0$ the flow is the harmonic rotation for
every $h,z$.

The formal preparation and the characteristic criterion are given by
Theorems~\ref{v4:formal-approximation}
and~\ref{v4:characteristic-factorization}. We now construct the complete
operators on a common neighborhood of the admissible pair.

\subsection{Domains and the two complete normal kernels}
Write the periodic comparison family of
\eqref{v4:periodization-definition} and Lemma~\ref{v4:periodic-comparison}
as $\widetilde\Psi_s=(\widetilde p_s,\widetilde q_s)$ and set
\[
 w_s=(D_x\widetilde\Psi_s)^{-1}e_2,\qquad
 (\mathcal M_H a)(x)=\int_0^{2\pi}w_s(x)a(\widetilde q_s(x),s)\,\dd s.
\]
This is the comparison operator, rather than the differential formed
from the true flow.  We use Lebesgue densities in the initial-value
plane and the product density $\dd q\,\dd s$ in the control cylinder.

Choose a bounded convex outer initial-value domain on which the true
flow, the periodization maps and their inverses have common margins.
It contains all the endpoint supports used in this section and their
line segments.  Choose real compactly supported base cutoffs
$\chi_A,\chi_D$, with $\chi_A=1$ on a neighborhood of
$\supp\chi_D$, and a real even cutoff $\theta(q)$ equal to one on
all velocity values over a smaller reduction domain.  Its entire support
is inside the fixed open interval $I=(-V,V)$. These choices are made with
strict margins at the harmonic reference.  Finite-time flow estimates
preserve the margins for the backgrounds under consideration.  In
particular, every evaluation point occurring in a branch formula or in
the periodization lies in one fixed compact subset of
$W=I\times\R/(2\pi\mathbb Z)$.

A collar of a compact set is an open neighborhood containing its
closed $\delta$-neighborhood for some $\delta>0$. We take successively
smaller such neighborhoods of $x_*$. On the outermost one used in an
elliptic reduction, the base cutoffs and the sampled velocity cutoff
are identically one and the specified principal symbols admit the
divisions used there. We call this the reduction domain (or platform).
The coefficient domain is chosen inside it so that the factorization
of the scalar principal symbol holds on its closure. A local order
bound on an open domain means that multiplication on both sides by
arbitrary smooth cutoffs compactly supported there gives an operator
of that order. A global order bound applies on the entire plane.
The operators defined below retain their bounded-frequency and smooth
kernel terms. A word is a finite ordered product of the indicated
operators and cutoffs. A primitive amplitude is the kernel amplitude
of an individual factor before any products are taken.

The operators $\mathcal A$ and $\mathcal A^*$ pass between the
control cylinder and the initial-value plane. The normal operators
$K_\pm$, their normalized columns $T_\pm$ and the matrix $M$ act
on the plane. The matrix $\mathcal I$ maps the two scalar inputs of
$M$ back to controls, so that $\mathcal A\mathcal I=M$. The subscripts
$\pm$ specify the signs of the velocity frequency in the two branches.

Let $\mathcal A=\chi_A\mathcal M_H$ be the ambient extension.
Choose $\gamma_+\in C^\infty(\R;[0,1])$ equal to zero for
$\sigma\le1$ and equal to one for $\sigma\ge2$, and put
$\gamma_-(\sigma)=\gamma_+(-\sigma)$.  Define
\begin{equation}\label{v4:lin-normals}
 C_\pm=\theta\gamma_\pm(D_q)\theta,\qquad
 K_\pm=\mathcal A C_\pm\mathcal A^*,\qquad D_q=-i\partial_q.
\end{equation}
The filters are smooth; their transition frequencies remain in the
operators.  Intermediate branch operators act on complex Sobolev
spaces: a single frequency sign need not preserve real functions.
The final physical inverse is made real only after the exact vector
equation has been solved.  The adjoint and the resulting complete matrix kernels are
\begin{align}
 (\mathcal A^*v)(q,s)
 &=\int\chi_A(x)w_s(x)\cdot v(x)
                 \delta(q-\widetilde q_s(x))\,\dd x,\label{v4:lin-adjoint}\\
 K_{\pm,ij}(x,y)
 &=\frac{\chi_A(x)\chi_A(y)}{2\pi}
    \int_0^{2\pi}\int_\R
       e^{i\sigma(\widetilde q_s(x)-\widetilde q_s(y))}
       w_{s,i}(x)w_{s,j}(y)\nonumber\\
 &\hspace{18mm}\cdot\theta(\widetilde q_s(x))
       \theta(\widetilde q_s(y))\gamma_\pm(\sigma)
                         \dd\sigma\,\dd s.\label{v4:lin-complete-kernel}
\end{align}
These identities follow by inserting the Fourier kernel of the
one-dimensional multiplier between the incidence kernels. This is
the normal-operator representation for a generalized Radon transform
with a directional weight \cite{Quinto1980,HomanZhou2017}. Here the
change of variables can be made for every pair of endpoints.

For every pair in the convex outer domain, set
\[
 n_s(x,y)=\int_0^1d_x\widetilde q_s(y+t(x-y))\,\dd t.
\]
Then $\widetilde q_s(x)-\widetilde q_s(y)=(x-y)\cdot n_s(x,y)$.
At the harmonic reference $n_s=m_s=(-\sin s,\cos s)$.
Lemma~\ref{v43:all-pairs-twist}, applied on this convex outer domain,
gives
\begin{equation}\label{v4:lin-all-pairs-twist}
 |n_s(x,y)|\ge c_0,\qquad
 \det(n_s(x,y),\partial_sn_s(x,y))\ge c_0>0
\end{equation}
for all endpoint pairs, with the common choice $c_0=7/16$.
The same lemma proves that the normalized direction has degree one
and is strictly increasing. For either sign of $\sigma$, the change
of variables
\[
 (s,\sigma)\longmapsto\xi=\sigma n_s(x,y)
\]
therefore covers $\R^2\setminus0$ once and has absolute Jacobian
$|\sigma|\det(n_s,\partial_sn_s)$.  Let
$(s_\pm(x,y,\xi),\sigma_\pm(x,y,\xi))$ be its inverses.
Formula~\eqref{v4:lin-complete-kernel} becomes exactly
\begin{equation}\label{v4:lin-fixed-diagonal}
 K_{\pm,ij}(x,y)=(2\pi)^{-2}\int e^{i(x-y)\cdot\xi}
                                   k_{\pm,ij}(x,y,\xi)\,\dd\xi,
\end{equation}
where
\begin{equation}\label{v4:lin-full-amplitude}
 k_{\pm,ij}=\left.
 \frac{2\pi\chi_A(x)\chi_A(y)w_{s,i}(x)w_{s,j}(y)
       \theta(\widetilde q_s(x))\theta(\widetilde q_s(y))
       \gamma_\pm(\sigma)}
 {|\sigma|\det(n_s,\partial_sn_s)}
                  \right|_{(s,\sigma)=(s_\pm,\sigma_\pm)}.
\end{equation}
The formula is set to zero near $\xi=0$, where the frequency filter
already vanishes.  It includes the finite-frequency transition.
Thus $K_\pm$ have order $-1$ with the fixed phase $(x-y)\cdot\xi$.

The same incidence relation gives
\[
 q=\widetilde q_s(x),\quad
 \xi=\sigma d_x\widetilde q_s(x),\quad
 \eta_s=-\sigma\partial_s\widetilde q_s(x),\quad
 \eta_q=\sigma.
\]
The map $x\mapsto(\widetilde q_s(x),\partial_s\widetilde q_s(x))$
has derivative uniformly close to the invertible harmonic matrix
with rows $m_s,m_s'$. Integrating its derivative along segments in
the convex chart gives injectivity and a uniform lower Lipschitz
bound, as well as a nonzero Jacobian. Each frequency
sign is consequently a graph branch. This graph observation explains
the two branches; the exact change of variables above is what controls
their complete kernels.

\subsection{Divergence normalization and complete reference identities}
Use the full multiplier $\Lambda=(1-\Delta)^{1/2}$ and the rows
\begin{equation}\label{v4:lin-rows}
 \mathcal D=\Lambda^{-1}\chi_D\Div,\qquad
 \mathcal C=\Lambda^{-1}\chi_D\Curl.
\end{equation}
The derivative acts before the cutoff.  This ordering matters.  At
$\varepsilon=0$, the local comparison has the exact expression
\[
 \mathcal M_0 a=\nabla_x\int_0^{2\pi}
                  \left(\int_0^{m_s\cdot x}a(q,s)\,\dd q\right)\dd s.
\]
Since $\chi_A=1$ near $\supp\chi_D$, it follows that
\begin{equation}\label{v4:lin-exact-zero-curl}
 \mathcal C K_{\pm,0}=0
\end{equation}
as an identity of complete operators.  A cutoff inserted before curl
would introduce an additional derivative and would not have this property.

Put $N_\pm=\mathcal D K_\pm\mathcal D^*$.  These scalar operators
have order $-1$ and are elliptic on the reduction platform.  Indeed, if
$n_s=d_x\widetilde q_s$ and $\Delta_s=\det(n_s,\partial_sn_s)$,
then $n_s\cdot w_s=1$.  On the branch $\xi=\sigma n_s$ write
$e=\xi/|\xi|$ and $e^\perp=Je$, where $J(a,b)=(-b,a)$.
Formula~\eqref{v4:lin-full-amplitude} gives
\[
 \sigma_{-1}(K_\pm)
  =\frac{\kappa_\pm}{|\xi|}w_{s_\pm}\otimes w_{s_\pm},
 \qquad \kappa_\pm=\frac{2\pi|n_{s_\pm}|}{\Delta_{s_\pm}},
\]
and therefore
\[
 \sigma_{-1}(N_\pm)
 =\frac{\kappa_\pm(e\cdot w_{s_\pm})^2}{|\xi|}
 =\frac{2\pi}{|\xi|\,|n_{s_\pm}|\Delta_{s_\pm}}.
\]
Choose one real frequency excision, zero near zero and one at high
frequency.  On the common platform, invert this principal symbol,
multiply by the excision, extend with fixed spatial cutoffs, and use
left quantization.  This specifies an actual order-one operator
$Q_\pm$ whose principal symbol is
\begin{equation}\label{v4:lin-Q}
 q_{1,\pm}=\frac{|\xi|\,|n_{s_\pm}|\Delta_{s_\pm}}{2\pi}.
\end{equation}
Both signs use the same extension, endpoint cutoffs, frequency profile
and quantization rule.
Define
\begin{equation}\label{v4:lin-branch-injection}
 B_\pm=C_\pm\mathcal A^*\mathcal D^*Q_\pm,\qquad
 T_\pm=\mathcal A B_\pm=K_\pm\mathcal D^*Q_\pm.
\end{equation}
The injections have compact velocity support.  The complete columns
$T_\pm$ have order zero, and their principal symbols are
\begin{equation}\label{v4:lin-column-principal}
 t_\pm=-i\frac{w_{s_\pm}}{e\cdot w_{s_\pm}},\qquad
 \sigma_0(\mathcal D T_\pm)=1,\qquad
 \sigma_0(\mathcal C T_\pm)=
 \Gamma_\pm=\frac{e^\perp\cdot w_{s_\pm}}{e\cdot w_{s_\pm}}.
\end{equation}
The signs follow from $D_j=-i\partial_{x^j}$ and
$\mathcal D^*=-\nabla\chi_D\Lambda^{-1}$.  All identities of
symbols are read on the common platform; the operators themselves
retain their endpoint cutoffs.

\begin{lemma}[Complete matching at the harmonic reference]\label{v4:lin-reference}
For every $h,z$, at $\varepsilon=0$ one has
\[
 K_{-,0}=K_{+,0},\qquad Q_{-,0}=Q_{+,0},\qquad T_{-,0}=T_{+,0}.
\]
The equality includes all finite-frequency terms.
\end{lemma}
\begin{proof}
At zero parameter, $\widetilde q_{s+\pi}(x)=-\widetilde q_s(x)$
and $w_{s+\pi}=-w_s$.  In~\eqref{v4:lin-complete-kernel} make
the substitution $(s,\sigma)\mapsto(s+\pi,-\sigma)$.
The phase is unchanged, the two vector signs cancel, the even velocity
cutoff is unchanged, and the filters are interchanged.  The integration
is over a full time circle.  Hence the two kernels agree exactly.
At the same reference $|n_s|=\Delta_s=1$, so the identical complete
prescription~\eqref{v4:lin-Q} gives identical $Q_\pm$.
Formula~\eqref{v4:lin-branch-injection} proves the remaining equality.
\end{proof}

\subsection{The normalized physical matrix}
For nonzero parameter set
\begin{equation}\label{v4:lin-normalized-columns}
 \mathcal I=\left[B_+,\frac{B_--B_+}{\varepsilon}\right],\qquad
 M=\left[T_+,\frac{T_--T_+}{\varepsilon}\right].
\end{equation}
Then $\mathcal A\mathcal I=M$ is an exact equality in the original
two-component target.  Lemma~\ref{v4:lin-reference} and Hadamard's
formula give
\begin{equation}\label{v4:lin-hadamard}
 \frac{T_-(\varepsilon,h,z)-T_+(\varepsilon,h,z)}{\varepsilon}
 =\int_0^1\partial_\varepsilon(T_--T_+)
                    (t\varepsilon,h,z)\,\dd t.
\end{equation}
Because the complete kernels have already been put on the fixed
diagonal, the derivative on the right preserves their order.  Thus
$M$ extends smoothly through zero as an order-zero matrix.  There is
no corresponding assertion that $B_{-,0}=B_{+,0}$: the injections
use different source-frequency branches.  The factor
$|\varepsilon|^{-1}$ in $\mathcal I$ remains in every control estimate.

Applying the rows in~\eqref{v4:lin-rows}, write
\begin{equation}\label{v4:lin-row-matrix}
 \begin{pmatrix}\mathcal D\\\mathcal C\end{pmatrix}M
   =\begin{pmatrix}\Pi&A_1\\B_1&F_1\end{pmatrix}.
\end{equation}
The principal divergence symbols of both columns in
\eqref{v4:lin-branch-injection} are one at every background.
Consequently both $\Pi-\Id$ and $A_1$ have local order at most $-1$.
The scalar Schur principal symbol is consequently
\begin{equation}\label{v4:lin-schur}
 \varsigma_{\varepsilon,h,z}(x,\xi)
       =\frac{\Gamma_-(x,\xi)-\Gamma_+(x,\xi)}{\varepsilon}.
\end{equation}
At a diagonal endpoint write
$A_s=D_x\widetilde\Psi_s$, with row covectors $\alpha_s,\beta_s$,
and put $w_s=A_s^{-1}e_2$.  On either branch
$e=\pm\beta_s/|\beta_s|$ and $\beta_s\cdot w_s=1$.  Hence
\begin{equation}\label{v41:lin-gamma-rows}
 \Gamma=\frac{(Je)\cdot w_s}{e\cdot w_s}
       =\det(\beta_s,w_s)
       =\frac{\alpha_s\cdot\beta_s}{\det A_s}.
\end{equation}
For the last equality, if the rows of $A_s$ are $(a,b)$ and $(c,d)$,
then $w_s=(-b,a)/\det A_s$.  The factors $|\beta_s|$ cancel in
the quotient, for both signs of the branch.

\begin{proposition}[The comparison symbol at a general background]
\label{v43:schur-background}
Let $f(q,s)=h(s)q^2+z(q,s)$ be real, smooth and $T$-periodic in
$s$, with $T=2\pi$. No half-periodicity is assumed for $z$.
Put
\begin{align}
 \mathsf U_s^f(x)&=\int_0^s m_t f(m_t\cdot x,t)\,\dd t,\notag\\
 \mathsf L_s^f(x)&=D_x\mathsf U_s^f(x)
       =\int_0^s f_q(m_t\cdot x,t)m_t\otimes m_t\,\dd t,\notag\\
 \mathsf B_s^f(x)&=\mathsf L_s^f(x)-\chi(s)\mathsf L_T^f(x).
 \label{v43:comparison-first-matrices}
\end{align}
Extend these integrals to real $s$ and set $\chi(s+T)=\chi(s)+1$.
Then $\mathsf B_s^f$ is symmetric and $T$-periodic. The normalized
Schur symbol at zero parameter is
\begin{equation}\label{v43:schur-general-formula}
 \varsigma_{0,f}(x,m_s)
   =2e_s^{\mathsf T}
        \bigl(\mathsf B_{s+\pi}^f(x)-\mathsf B_s^f(x)\bigr)m_s,
 \qquad e_s=(\cos s,\sin s).
\end{equation}
It is linear in $f$, real, and odd in the cotangent direction.
For each finite integer $r\ge0$ and fixed compact initial-value
domain $\overline U$ in the comparison chart,
\begin{align}
 \|\varsigma_{0,f}\|_{C^r(\overline U\times S^1)}
       &\le C_r\|f\|_{C^{r+1}(W)},\notag\\
 \|\varsigma_{0,h,z}+8g_h\|_{C^r(\overline U\times S^1)}
       &\le C_r\|z\|_{C^{r+1}(W)}.
 \label{v43:schur-zero-continuity}
\end{align}
Here $g_h(x,s)=g_{h,x}(s)$.
For every allowed real background the finite-parameter symbol is
also real and odd.
\end{proposition}
\begin{proof}
The first variation of the integral equation in Section~\ref{v4:endpoint-chapter}
gives, on the fixed flow tube,
\[
 \Psi_{\varepsilon f}(x,s)
     =C_s\bigl(x+\varepsilon\mathsf U_s^f(x)\bigr)
            +O(\varepsilon^2),\qquad
 N(\varepsilon f)=\varepsilon\mathsf U_T^f+O(\varepsilon^2).
\]
Differentiating the inverse identity defining $Q_s$ therefore gives
\[
 Q_s(x)=x-\varepsilon\chi(s)\mathsf U_T^f(x)
               +O(\varepsilon^2).
\]
Consequently
\begin{align}
 \widetilde\Psi_{\varepsilon f}(x,s)
 &=C_s\bigl(x+\varepsilon[
       \mathsf U_s^f(x)-\chi(s)\mathsf U_T^f(x)]\bigr)
                +O(\varepsilon^2),\notag\\
 D_x\widetilde\Psi_{\varepsilon f}(x,s)
 &=C_s\bigl(\Id+\varepsilon\mathsf B_s^f(x)\bigr)
                +O(\varepsilon^2).
 \label{v43:comparison-general-first-variation}
\end{align}
The remainders hold in every fixed finite norm when enough
derivatives of $f$ are bounded. Periodicity of the integrands gives
$\mathsf L_{s+T}^f=\mathsf L_s^f+\mathsf L_T^f$, proving the
asserted periodicity of $\mathsf B^f$. Its symmetry follows directly
from~\eqref{v43:comparison-first-matrices}.

The two rows in~\eqref{v43:comparison-general-first-variation}
and~\eqref{v41:lin-gamma-rows} now give
\[
 \Gamma_{\varepsilon f}(x,s)
       =2\varepsilon e_s^{\mathsf T}\mathsf B_s^f(x)m_s
                      +O(\varepsilon^2).
\]
At the fixed direction $m_s$, the reference times of the positive
and negative branches are $s$ and $s+\pi$. Their parameter
derivatives do not contribute to the first coefficient, since
$\Gamma_{0,f}$ vanishes identically in time. Subtracting the
two coefficients proves~\eqref{v43:schur-general-formula}.
In particular the matrix between $e_s^{\mathsf T}$ and $m_s$ is
\begin{equation}\label{v43:schur-endpoint-term}
 \int_s^{s+\pi} f_q(m_t\cdot x,t)m_t\otimes m_t\,\dd t
       -[\chi(s+\pi)-\chi(s)]\mathsf L_T^f.
\end{equation}
The endpoint term is retained for a general correction $z$.

The formula is real and linear in $f$. Replacing $s$ by $s+\pi$
changes its sign because $\mathsf B^f$ is $T$-periodic.
Each of its $r$ spatial and angular derivatives uses at most
$r+1$ derivatives of $f$ on the fixed sampling collar. This proves
the first estimate in~\eqref{v43:schur-zero-continuity}.
For the seed $f=hq^2$, the integrand defining $\mathsf U^f$
changes sign under $t\mapsto t+\pi$, so
$\mathsf U_T^f=\mathsf L_T^f=0$. The first variations preceding
\eqref{v4:char-determinant} yield
$\mathsf L_s^f=\mathsf K_s-\mathsf K_0$ and
$\mathsf K_{s+\pi}=-\mathsf K_s$. Thus
$\varsigma_{0,h,0}=-4e_s^{\mathsf T}\mathsf K_sm_s=-8g_{h,x}(s)$.
Linearity proves the second estimate.

Finally, at any finite parameter the all-pair change of variables
gives $s_+(x,-\xi)=s_-(x,\xi)$ and the reverse identity.
Moreover $e\cdot w_{s_\pm}=\pm|\beta_{s_\pm}|^{-1}\ne0$.
Changing $e$ to $-e$ changes both signs in the quotient
\eqref{v41:lin-gamma-rows}; hence
$\Gamma_+(x,-\xi)=\Gamma_-(x,\xi)$. These quotients are real.
This proves realness and oddness for nonzero parameter, and
continuity gives the same conclusion at zero.
\end{proof}

In particular, with $t_s=e_s^{\mathsf T}\mathsf K_sm_s$
as in~\eqref{v4:char-T}, the seed calculation gives
\begin{equation}\label{v41:lin-schur-seed}
 \varsigma_{\varepsilon,h,0}(x,\xi)
     =-4t_s+O(\varepsilon)
     =-8g_{h,x}(s)+O(\varepsilon).
\end{equation}
Here $s$ is defined by $\xi=|\xi|m_s$ at zero parameter.
The true flow and its periodic comparison have the same first
variation in this seed case because $\mathsf U_T^{hq^2}=0$.
Multiplication by $-\Lambda$ therefore gives the reference principal symbol
$8|\xi|g_{h,x}(s)$ of~\eqref{v4:char-principal-model}.

\begin{lemma}[A finite norm controlling the characteristic]
\label{v43:schur-continuity}
On a bounded $C^6(W)$ neighborhood of the reference forcing
$f^\circ=h^\circ q^2$, with the fixed flow and inverse margins,
there are $\varepsilon_0>0$ and $C$ such that
\begin{equation}\label{v43:schur-epsilon-continuity}
 \|\varsigma_{\varepsilon,f}-\varsigma_{0,f}\|
             _{C^1(\overline U\times S^1)}
       \le C|\varepsilon|,\qquad |\varepsilon|\le\varepsilon_0.
\end{equation}
Consequently, on this neighborhood,
\begin{equation}\label{v43:schur-reference-continuity}
 \|\varsigma_{\varepsilon,h,z}+8g_{h^\circ}\|_{C^1}
 \le C\bigl(|\varepsilon|+
          \|h-h^\circ\|_{C^2}+\|z\|_{C^2}\bigr).
\end{equation}
\end{lemma}
\begin{proof}
Adjoin $\varepsilon$ as a constant parameter to the oscillator.
The finite variational equations used in
Lemma~\ref{v43:endpoint-tame} give uniform joint derivatives
of the flow through degree four in $(x,s,\varepsilon)$ on the
bounded $C^6$ ball. Differentiating
$(\Id+\chi N)\circ Q_s=\Id$ isolates each highest derivative
of $Q_s$ with the uniformly invertible coefficient
$\Id+\chi DN$. Thus the comparison flow has the same finite
joint bounds.

The direction equation for a branch is the direction of
$d_x\widetilde q_s$. Its time derivative is bounded below by
Lemma~\ref{v43:all-pairs-twist}. Differentiating this scalar
equation isolates each highest branch derivative with that
same nonzero denominator. For example, write the direction equation as
$G(\varepsilon,s)=0$ with $x$ and the direction fixed. On either branch,
\[
 s_\varepsilon=-G_\varepsilon/G_s,\qquad
 s_{\varepsilon\varepsilon}
 =-(G_{\varepsilon\varepsilon}
       +2G_{s\varepsilon}s_\varepsilon+G_{ss}s_\varepsilon^2)/G_s .
\]
One further derivative in $x$ or in the angular variable uses at most
three derivatives of $G$ and another power of the same bounded reciprocal.
The function $G$ already contains one derivative of $\widetilde q_s$,
so these terms use the comparison flow through degree four.
Formula~\eqref{v41:lin-gamma-rows}, whose denominator is uniformly
nonzero, therefore bounds
$\partial_\varepsilon^2(\Gamma_--\Gamma_+)$ in $C^1$.
Only these finite derivatives are required; the choice $C^6$
leaves a fixed margin in this count.

Since $\Gamma_--\Gamma_+=0$ at zero parameter, Hadamard's formula
gives
\[
 \varsigma_{\varepsilon,f}
    =\int_0^1\partial_\varepsilon(\Gamma_--\Gamma_+)
                                      (t\varepsilon,f)\,\dd t.
\]
Subtract its value at zero and integrate once more in the
parameter. The preceding bound proves
\eqref{v43:schur-epsilon-continuity}. Combining this with
\eqref{v43:schur-zero-continuity} and the linear dependence
of $g_h$ on $h$ proves~\eqref{v43:schur-reference-continuity}.
\end{proof}

We apply Theorem~\ref{v4:characteristic-factorization} in this
$C^1(\overline U\times S^1)$ topology. On a small fixed
neighborhood of the reference roots, their transverse derivatives
remain bounded away from zero; on the complementary compact set
the symbol remains nonzero. Thus there is the same simple
characteristic line and a real elliptic quotient. The finite
Sobolev index $b_0$ used below is chosen once to include
$b_0\ge8$, which controls $C^6$ on the two-dimensional cylinder,
and the other finite low indices of the construction. A small
$H^{b_0}$ background ball then gives both the bounded $C^6$
neighborhood and the required small $C^2$ perturbation.
After a fixed
rotation of the target chart, the first-order principal symbol has the form
\begin{equation}\label{v4:lin-factor}
 a_0(x,\xi)\bigl(\xi_1-\beta(x)\xi_2\bigr),
\end{equation}
with $a_0$ real and elliptic of degree zero and $\beta$ real.
The opposite root has the same slope, by the branch exchange identity.
This controls the whole cosphere, rather than only one conic sector.

\subsection{Finite reductions and uniform dependence on the seed}
All elliptic inverses used next are finite parametrices.  For example,
put $R_\Pi=\Pi-\Id$ and
\[
 Q=\sum_{j=0}^{d-1}(-R_\Pi)^j,\qquad
 E_\Pi=Q\Pi-\Id=\Pi Q-\Id=-(-R_\Pi)^d.
\]
These are actual global operator identities.  In addition,
$QE_\Pi=E_\Pi Q$, because both are polynomials in $R_\Pi$.
Only their order
assertions are local: $E_\Pi$ has order $-d$ on the reduction
platform, with separated input pieces retained.  Matrix multiplication
gives the exact formula
\begin{equation}\label{v4:lin-triangular}
 \begin{pmatrix}Q&0\\-B_1Q&1\end{pmatrix}
 \begin{pmatrix}\Pi&A_1\\B_1&F_1\end{pmatrix}
 \begin{pmatrix}1&-QA_1\\0&1\end{pmatrix}
 =\begin{pmatrix}1&0\\0&F_1-B_1QA_1\end{pmatrix}
 +\begin{pmatrix}
 E_\Pi&-QE_\Pi A_1\\-B_1E_\Pi&B_1QE_\Pi A_1
 \end{pmatrix}.
\end{equation}
Multiplication of the second row by $-\Lambda$ raises error orders
by at most one.  Taking $d=N+1$ makes the displayed matrix error of
order at most $-N$ on a smaller platform.  Its complete exterior
terms are not discarded.

For completeness, the dependence required here can be read directly
from the formulas.  Finite-time flow variations solve linear Volterra
equations; derivatives of the inverse comparison maps are found by
differentiating their inverse identities.  The branch times are
obtained from a scalar implicit equation whose derivative is bounded
below by~\eqref{v4:lin-all-pairs-twist}.  Formula
\eqref{v4:lin-full-amplitude} then uses products, these branch
evaluations, and reciprocals of uniformly nonzero quantities.
The root in~\eqref{v4:lin-factor} is another nondegenerate implicit
function.  Near the root its elliptic quotient is the integral of
the transverse derivative; away from the root it is ordinary division.
All these operations have fixed derivative costs.  Increasing the
Sobolev index increases differentiation order, but does not repeat the
construction or increase its finite depth.

\begin{lemma}[The finite list of target operators]\label{v43:finite-target-words}
The target reductions are finite sums of words generated by the
complete normal amplitudes $K_\pm$, the principal prescriptions
$Q_\pm$, the rows $\mathcal D,\mathcal C$ and their adjoints,
fixed cutoffs, and the full multipliers $\Lambda^{\pm1}$.
The finite scalar reduction in Section~\ref{v4:lin-inverse-section}
also uses the prescribed principal factors $A_0,E_0,J_0$ and the
first-order vector field $V_0$. The normalized second column is
formed by one parameter derivative of the primitive products
and the fixed integral in~\eqref{v4:lin-hadamard}.
The list and its maximal word length are fixed before the Sobolev
index is varied.

The global orders and the improved local orders used later are
as follows. A local order is asserted only after restriction to
the common reduction platform.
\begin{equation}\label{v43:raw-word-orders}
\begin{array}{c|c|c}
 \text{operator}&\text{global order at most}&\text{local order at most}\\ \hline
 K_\pm,\ Q_\pm&-1,\ 1&-1,\ 1\\
 \mathcal D,\mathcal C,\mathcal D^*,\mathcal C^*&0&0\\
 T_\pm,\ M,\Pi,B_1,F_1&0&0\\
 A_1,\ \Pi-\Id&0&-1\\
 Q,\ E_\Pi&0,\ 0&0,\ -d\\
 P&1&1\\
 A_0,E_0,E_N,V,L&0&0\\
 C_0&1&0\\
 U,\ S_N,\ T_N&1,\ 1,\ 0&1,\ -N,\ -N
\end{array}
\end{equation}
For each such actual finite word of global order $m$ there are
fixed finite indices $b,q$, independent of $s$, for which
\begin{equation}\label{v4:lin-word-bound}
 \|Tf\|_{H^s}\le C_s\bigl[\|f\|_{H^{s+m}}
  +(1+\|h\|_{H^{s+q}}+\|z\|_{H^{s+q}})\|f\|_{H^b}\bigr],
\end{equation}
on a common low background neighborhood. The same estimate with
a lower order applies to a localized defect only after its local
symbol cancellation and its separated input pieces have been
estimated. In particular a defect containing a global identity
term is not assigned a negative global order.
\end{lemma}
\begin{proof}
The full amplitudes~\eqref{v4:lin-full-amplitude} have order $-1$;
the prescribed reciprocals~\eqref{v4:lin-Q} have order $1$.
The four divergence and curl rows have order zero. Thus
$T_\pm=K_\pm\mathcal D^*Q_\pm$ has order zero. To exhibit the
derivative budget, commute a multi-index $\alpha$ with this
actual product:
\begin{equation}\label{v43:column-commutator}
\begin{aligned}
 \operatorname{ad}_{\partial}^{\alpha}
       (K_\pm\mathcal D^*Q_\pm)
 &=\sum_{\alpha_1+\alpha_2+\alpha_3=\alpha}
   \binom{\alpha}{\alpha_1,\alpha_2,\alpha_3}
   (\operatorname{ad}_{\partial}^{\alpha_1}K_\pm)\\
 &\qquad\cdot
   (\operatorname{ad}_{\partial}^{\alpha_2}\mathcal D^*)
   (\operatorname{ad}_{\partial}^{\alpha_3}Q_\pm).
\end{aligned}
\end{equation}
Here $\operatorname{ad}_{\partial_j}T=[\partial_j,T]$.
Use the fixed mapping indices
$L^2\to H^{-1}\to H^{-1}\to L^2$, from right to left.
The fixed-index estimate~\eqref{v4:kernel-fixed-index} then
requires only a fixed finite frequency list for each factor.
The primitive triangular estimate~\eqref{v4:kernel-triangular}
uses at most $|\alpha_\nu|+q_\nu$ background derivatives in
factor $\nu$. After the bounded low orders have been absorbed,
interpolation bounds their product by one background norm at
$|\alpha|+q$, with $q$ fixed. The identity
\[
 \partial^\alpha(Tf)=
   \sum_{\gamma\le\alpha}\binom{\alpha}{\gamma}
     (\operatorname{ad}_{\partial}^{\gamma}T)
                                    \partial^{\alpha-\gamma}f
\]
and the same interpolation between the low source norm and
$H^s$ prove~\eqref{v4:lin-word-bound} for this order-zero column.
One parameter derivative first applied to its primitive factors
adds a fixed cost to $q$; integration over $0\le t\le1$ preserves
the bound. This proves the estimate for the second column
without dividing a target estimate by $\varepsilon$.

For a word $T=A_1\cdots A_k$ of factor orders $m_1,\ldots,m_k$,
put $r_j=\sum_{i=1}^j m_i$ and $r_0=0$. Every commutator on its
$j$th factor uses the fixed mapping $H^{r_j}\to H^{r_{j-1}}$.
With $\sum_j|\gamma_j|=a$, its product therefore maps
$H^{r_k}$ to $L^2$ with total coefficient excess at most $a$
over a common low index $b_c$. For $s>b$ large enough that
$s+r_k>b$ and $|\alpha|\le s$, the source in a term of
$\partial^\alpha(Tf)$ has order at most $s-a+r_k$. Interpolate
between the low and high source norms $H^b,H^{s+r_k}$, and between
the background norms $H^{b_c},H^{s+q}$. Choose $b\ge r_k$ and
$q\ge b_c$. The sum of the high weights is at most
\[
 \frac{a}{s+q-b_c}
   +\frac{(s-a+r_k-b)_+}{s+r_k-b}\le1,\qquad 0\le a\le s.
\]
If the second numerator vanishes, only the first fraction occurs.
The fixed mapping estimates absorb the remaining bounded low
indices. Weighted arithmetic--geometric mean now gives
\eqref{v4:lin-word-bound} with $m=r_k$ at the integer levels used
in the iteration. For a fixed fractional shift $0\le\sigma<1$ take
the same commutator expansion at $s_0=\lfloor s\rfloor$, with
$H^{r_j+\sigma}\to H^{r_{j-1}+\sigma}$ at factor $j$.
This gives the preceding source order $s-a+r_k$ for $s=s_0+\sigma$.
The finitely many words and their adjoints
use a common maximum of the fixed offsets.

The row matrix and the triangular factors are globally of
order zero. The principal identities
\eqref{v4:lin-column-principal} give the two local order drops
for $A_1$ and $\Pi-\Id$. Their finite polynomial then gives
the local order $-d$ for $E_\Pi$. Multiplication by
$-\Lambda$ gives $P$ global order at most one.
The definitions~\eqref{v4:lin-scalar-finite} and
\eqref{v4:lin-matrix-identities} give the remaining global
orders in the table. In particular $E_N$ is a finite sum of
order-zero words; the negative local orders belong to its
defect and to $S_N,T_N$, after cancellation on the platform.
The complete difference $C_0=P-A_0V_0$ has its order-one
principal part cancelled there.

For these local improvements, use the finite Taylor formula
of Appendix~\ref{v4:kernel-appendix}, retaining its integral
remainder. Each cancellation lowers the local diagonal order
by the stated number. Insert the fixed collar cutoffs between
factors; every complementary piece has a separated primitive
kernel, or a complete multiplier tail, controlled by
Lemmas~\ref{v4:kernel-separated} and
\ref{v4:kernel-multiplier-tails}. This proves the localized
estimates with the same finite-list convention.
The estimates on completed composition remainders use fixed
frequency seminorm lists; these remainders are not regarded
as new primitive amplitudes with an unrestricted triangular
bound. For a coordinate change, the primitive factors are
conjugated first and then composed. Proposition
\ref{v43:kernel-conjugated-finite-word} proves that this
procedure gives the scalar grading~\eqref{v43:kernel-B-grading}
with an offset independent of $s$.
\end{proof}

In~\eqref{v4:lin-word-bound}, norms of $h$ are on its time circle and
norms of $z$ are on $W$; equivalently $h$ can be extended constantly
in $q$ with the chosen cutoff.  The one parameter derivative used in
\eqref{v4:lin-hadamard} is included in the fixed offset $q$.
None of these target operators has a negative power of $\varepsilon$.

These estimates apply around every admissible pair $(h^\circ,x_*)$.
Indeed, write $\widetilde z=(h-h^\circ)q^2+z$ on the fixed velocity
interval.  After a fixed outer cutoff, this is a bounded linear map
on every same-order Sobolev space, and
$H_{\varepsilon,h,z}=\varepsilon(h^\circ q^2+\widetilde z)$.
Thus seed variation uses the same complete control-background chart;
no separate high-order uniformity assumption is required.
The domains, low indices and constants may depend on that pair.
The construction applies to smooth seeds.

The following table records the passage from the branch kernels to the
physical inverse.  Here an order refers to pseudodifferential order unless
a Sobolev loss or gain is stated.  The constructions of $U,V,L,C,R,G$
are given in the next section; the scope column distinguishes exact
identities from local order statements.
\begin{center}
\small
\begin{tabular}{@{}>{\raggedright\arraybackslash}p{.22\textwidth}>{\raggedright\arraybackslash}p{.17\textwidth}>{\raggedright\arraybackslash}p{.14\textwidth}>{\raggedright\arraybackslash}p{.35\textwidth}@{}}
\hline
Stage & Order or loss & Parameter cost & Scope\\
\hline
$K_\pm$, $Q_\pm$ & $-1$, $1$ & bounded & Complete kernels; $Q_\pm$ invert the principal symbol of $N_\pm$ on the platform.\\
$M$, $\mathcal I$ & $0$; loss $1$ & bounded; $|\varepsilon|^{-1}$ & $\mathcal A\mathcal I=M$ globally.\\
$U$, $V$, $L$ & at most $1,0,0$ & bounded & Exact finite identities; defects locally of order $-4$.\\
$C$ & Sobolev loss $0$ & bounded & $DC=\Id$ on a common inner domain.\\
$R$, $G$ & gain $3/2$; loss $0$ & bounded & Localized target error; actual $L^2$ inverse.\\
$\mathcal B$ & Sobolev loss $2$ & $|\varepsilon|^{-1}$ & Exact comparison equation on $\Omega_-$.\\
\hline
\end{tabular}
\end{center}
\section{A local right inverse with complete error correction}\label{v4:lin-inverse-section}

The reduction in the preceding section has two distinct defects.
Its scalar transport equation contains a complete order-zero operator,
and its return to the original vector target leaves a smoothing error.
We remove them in this order.  The first correction uses a short
transport interval; the second uses a small output support.

\begin{theorem}[Local comparison inverse]\label{v4:local-inverse}
Fix a seed $h^\circ$ satisfying Definition~\ref{v4:seed-class-B},
and a nonzero point $x_*$ satisfying its characteristic condition.
There are fixed domains $\Omega_-\Subset\Omega_+$ containing $x_*$,
a compact velocity interval, finite integers $b_0,b,q$, and positive
constants $\delta,\varepsilon_0$ such that the following holds.  For
\[
 0<|\varepsilon|<\varepsilon_0,\qquad
 \|h-h^\circ\|_{H^{b_0}}+\|z\|_{H^{b_0}(W)}<\delta,
\]
there is a real linear map $\mathcal B_{\varepsilon,h,z}$ from smooth
vector fields on $\Omega_+$ to smooth real $2\pi$-periodic controls
with that fixed compact velocity support, such that
\begin{equation}\label{v4:lin-exact}
 \mathcal M_{H_{\varepsilon,h,z}}\mathcal B_{\varepsilon,h,z}f=f
                \quad\text{on }\Omega_-.
\end{equation}
For every integer $s\ge b$,
\begin{equation}\label{v4:lin-tame}
 \|\mathcal B_{\varepsilon,h,z}f\|_{H^s}
 \le\frac{C_s}{|\varepsilon|}
 \left[\|f\|_{H^{s+2}}
 +(1+\|h\|_{H^{s+q}}+\|z\|_{H^{s+q}})\|f\|_{H^{b+2}}\right].
\end{equation}
In this estimate the source norms are those of the zero extension
of $\zeta f$ to the plane, where $\zeta\in C_c^\infty(\Omega_+)$ equals one near
$\overline{\Omega_-}$. They are bounded by the corresponding norms
of $f$ on $\Omega_+$. The constants are uniform on the stated background
neighborhood. The same domains, support, low neighborhood, cutoff and
operators work for all $s$.
\end{theorem}

The theorem concerns the periodic comparison operator.  The difference
from the true return differential is the residual-dependent term of
Proposition~\ref{v4:comparison-defect}, which remains in the nonlinear
iteration.

\subsection{A complete scalar inverse on a short interval}
We first work on $\mathcal X=\R_t\times\mathbb T_y$, where the
transverse circumference is fixed.  Let $B$ be a complete order-zero
pseudodifferential operator with fixed compact endpoint localizations.
Its low symbol and kernel seminorms are bounded and
$\|B\|_{L^2\to L^2}\le M$.  The extension defining $B$ is chosen
before the interval below is shortened.

Choose fixed smooth functions $\eta_0,\chi_0$ supported in $(-2,2)$
and equal to one on $[-1,1]$.  For $0<\tau\le1$, let
\begin{equation}\label{v4:lin-integral}
 T_\tau f(t,y)=\eta_0(t/\tau)
       \int_{-\infty}^t\chi_0(a/\tau)f(a,y)\,\dd a.
\end{equation}
Schur's test in the $t$ variable gives
$\|T_\tau\|_{2\to2}\le C\tau$.  Fix $\tau$ so that
$C\tau M<1/2$ and define actual $L^2$ operators
\begin{equation}\label{v4:lin-scalar-inverse}
 G_B=(\Id+T_\tau B)^{-1},\qquad Q_B=G_BT_\tau.
\end{equation}
They have norms at most $2$ and $2C\tau$, respectively.

For the estimates in this subsection, represent $B$ in fixed charts by
compact-output left symbols $b_\ell$ and integral kernels $k_\ell$ with
compact endpoint supports. Write $\|k\|_{\rm Sch}$ for the larger of the
two Schur norms $\sup_u\int|k(u,v)|\,\dd v$ and
$\sup_v\int|k(u,v)|\,\dd u$. For this representation put
\[
 |B|_j=\sum_\ell\max_{|\alpha|\le j}
   \|\partial_u^\alpha b_\ell\|_{L^\infty_{u,\xi}}
  +\sum_\ell\max_{|\alpha|+|\gamma|\le j}
   \|\partial_u^\alpha\partial_v^\gamma k_\ell\|_{\rm Sch}.
\]
The representations used for the oscillator are supplied by
Proposition~\ref{v43:kernel-conjugated-finite-word}.

\begin{lemma}[The same scalar solution at all orders]\label{v4:lin-scalar-lemma}
The operator $Q_B$ preserves smoothness and $H^s$ for every real $s$.
Its output is supported in $|t|<2\tau$ and satisfies
\begin{equation}\label{v4:lin-scalar-equation}
 (\partial_t+B)Q_Bf=f\quad\text{on }|t|<\tau.
\end{equation}
For fixed finite $b,q_0$, and integer $s\ge b$,
\begin{equation}\label{v4:lin-scalar-tame}
 \|Q_Bf\|_{H^s}\le C_s
 \left[\|f\|_{H^s}+(1+|B|_{s+q_0})\|f\|_{H^b}\right].
\end{equation}
The constants may depend on the chosen $\tau$. The
adjoint satisfies analogous Sobolev and tame estimates, and continuity at fixed negative
indices uses only finitely many background seminorms.
\end{lemma}
\begin{proof}
For the $L^2$ solution $v=Q_Bf$ one has
\begin{equation}\label{v4:lin-scalar-fixedpoint}
 v=T_\tau(f-Bv).
\end{equation}
This proves its support.  To justify differentiating it, write
$I_\chi f(t)=\int_{-\infty}^t\chi(a)f(a)\,\dd a$.
Integration by parts gives
\[
 [\partial_t,T_\tau]=\eta_\tau'I_{\chi_\tau}
                   +\eta_\tau I_{\chi_\tau'},\qquad
 [\partial_y,T_\tau]=0.
\]
Every iterated coordinate commutator is a finite sum of the same
integrals; its $L^2$ norm is bounded by
$C_\alpha\tau^{1-|\alpha|}$. A coordinate commutator of $B$
differentiates its symbol or its kernel and is still $L^2$ bounded.

Put $K=T_\tau B$, $g=T_\tau f$.  If $\mathsf T_a$ is translation
in one coordinate and $K_a=\mathsf T_aK\mathsf T_{-a}$, then
\[
 (\Id+K_a)\frac{\mathsf T_av-v}{a}
 =\frac{\mathsf T_ag-g}{a}-\frac{K_a-K}{a}v.
\]
The inverse on the left has norm at most two, because translation is
unitary.  The commutator bounds control the second quotient, so
$f\in H^1$ implies $v\in H^1$.  Differentiate the equation weakly
through order $k-1$ and repeat the argument.  The commutators act on
at most $k-2$ derivatives of $v$, and their first weak derivatives
use only the already known $H^{k-1}$ regularity.  This proves
$Q_B:H^k\to H^k$ by induction with the same $L^2$ inverse.
Equation~\eqref{v4:lin-scalar-equation} now follows by differentiating
\eqref{v4:lin-scalar-fixedpoint} on the common plateau.

With the fixed chart extensions, the grading above has the interpolation
inequality for Banach-valued $C^j$ functions.  The Fourier--Schur
estimate of Appendix~\ref{v4:kernel-appendix} gives, for a fixed $q_0$,
\[
 \|\operatorname{ad}_{\partial}^{\alpha}B\|_{2\to2}
                  \le C_\alpha |B|_{|\alpha|+q_0}.
\]
In the kernel part the commutator differentiates by
$\partial_x+\partial_y$; its norm is bounded by the independent
endpoint grading just defined.

Fix $b\ge2$ and bound $|B|_{b+1+q_0}$ in the low neighborhood.
Write $Z_j=\|v\|_{H^j}$ and $X_s=1+|B|_{s+q_0}$.  For
$1\le j\le s$, coefficient interpolation and source interpolation
use the weights
\[
 \theta_j=\frac{j-1}{s-1},\qquad
 \psi_j=\max\left(0,\frac{s-j-b}{s-1-b}\right).
\]
For $s>b+1$, their sum is at most one: in the nontrivial range it is
$1-b(j-1)/((s-1)(s-1-b))$.  Weighted arithmetic--geometric mean
therefore gives
\[
 |B|_{j+q_0}Z_{s-j}\le C_s(Z_{s-1}+X_sZ_b).
\]
The differentiated equation and the same inverse $G_B$ now yield
\[
 Z_s\le C_s(\|f\|_{H^s}+Z_{s-1}+X_sZ_b).
\]
The fixed low levels give $Z_b\le C\|f\|_{H^b}$.  Interpolation
between $H^b$ and $H^s$ bounds the middle term by
$\lambda Z_s+C_{s,\lambda}Z_b$.  Choose $\lambda$ using only
the low coefficient constant, and absorb.  This proves
\eqref{v4:lin-scalar-tame}.  In particular, the high symbol norm
is never required to be small.

For the adjoint, use
\[
 Q_B^*=T_\tau^*(\Id+B^*T_\tau^*)^{-1},\qquad
 T_\tau^*f(a,y)=\chi_\tau(a)\int_a^\infty\eta_\tau(t)f(t,y)\,\dd t.
\]
Put $w=(\Id+B^*T_\tau^*)^{-1}f$. Translation is unitary and
$\|B^*T_\tau^*\|_{2\to2}=\|T_\tau B\|_{2\to2}<1/2$.
For every coordinate derivative $\partial$ on the cylinder,
$[\partial,B^*]=[\partial,B]^*$ since $\partial^*=-\partial$.
The identical identity holds for $T_\tau^*$. Differentiating the
equation for $w$ therefore requires precisely the commutator bounds
above, with the product order reversed. In particular, differentiating
$j$ times inside $B^*T_\tau^*$ costs at most $|B|_{j+q_0}$ at
fixed $\tau$. The preceding difference quotient and interpolation
estimates give the same $H^s$ bound for $w$ and then for $T_\tau^*w$.
At a fixed integer $k>0$, apply this bound at order $k$ (or its
finite-order version from the difference quotients) to obtain
$|\langle Q_B f,g\rangle|\le C_k\|f\|_{H^{-k}}\|g\|_{H^k}$.
Duality and interpolation now give the stated boundedness on every
fixed real Sobolev space, using only finitely many seminorms there.
The extensions agree on common domains by density.
\end{proof}

The extension defining $B$ and its $L^2$ bound are fixed before $\tau$.
Thus shortening the integration interval changes only $T_\tau$.

\subsection{The scalar factor and its retained localization errors}
Starting with~\eqref{v4:lin-triangular}, multiply its scalar row by
$-\Lambda$.  Call the resulting actual scalar operator $P$.
Its principal symbol is~\eqref{v4:lin-factor}.  Choose compactly
supported coefficient extensions of $\beta$ and the elliptic factor
$a_0$, using a smaller platform on which their margins are uniform.
Let
\[
 V_0=D_1-\beta D_2,\qquad v=iV_0=\partial_{x^1}-\beta\partial_{x^2}.
\]
Fix an order-zero quantization $A_0$ of $a_0$ with endpoint cutoffs,
and define the complete difference $C_0=P-A_0V_0$.  It has order zero
on the platform, although its raw global order is at most one.
This definition retains input-cutoff commutators and the derivative
of $\beta$ produced by integration by parts.

Let $E_0$ be the identically localized principal reciprocal of $A_0$,
put $R_0=E_0A_0-\Id$, and fix $d=N+1$.  Define
\[
 E_N=\sum_{j=0}^{d-1}(-R_0)^jE_0,\qquad
 R_N=-(-R_0)^d.
\]
The exact identities are
\begin{equation}\label{v4:lin-scalar-finite}
 E_NA_0=\Id+R_N,\qquad
 E_NP=V_0+E_NC_0+R_NV_0.
\end{equation}
Locally $R_N$ has order $-d$, so $R_NV_0$ has order at most $-N$.
The global identity term in $R_N$ will be retained in the localization
estimates below.

Choose fixed cutoffs $\theta_o,\theta_i$ with supports in the
coefficient platform, equal to one on a common smaller region, and put
\[
 \widetilde B=iE_NC_0,\qquad B_0=\theta_o\widetilde B\theta_i.
\]
The difference is retained exactly as
\begin{equation}\label{v4:lin-B-localization}
 \widetilde B-B_0=(1-\theta_o)\widetilde B
                       +\theta_o\widetilde B(1-\theta_i).
\end{equation}
The first term has exterior output.  On a smaller output region the
second has separated input and output.  Both will be included in the
vector error below.

\begin{lemma}[Uniform scalar bound]\label{v43:lin-B0-bound}
On a finite low background neighborhood there is a constant $M_0$ such
that $\|B_0\|_{L^2\to L^2}\le M_0$. The neighborhood and the extension
defining $B_0$ can be chosen independently of the transport interval.
\end{lemma}
\begin{proof}
The polynomial defining $E_N$ has the exact expansion
\begin{equation}\label{v43:lin-EN-polynomial}
 E_N=\sum_{k=0}^{d-1}(-1)^k\binom d{k+1}(E_0A_0)^kE_0.
\end{equation}
Thus $E_N$ is a finite sum of order-zero words with compact endpoint
cutoffs. Choose $\eta_B=1$ near $\supp\theta_o$, with support still
inside the platform where the principal symbols of $P$ and $A_0V_0$
agree. Then
\begin{equation}\label{v43:lin-B0-split}
 B_0=i\theta_oE_N\eta_BC_0\theta_i
       +i\theta_oE_N(1-\eta_B)C_0\theta_i.
\end{equation}
In the first term, localize the input of $\eta_BC_0\theta_i$ near its
output. Its diagonal part has order zero by the one-term finite
composition formula and principal cancellation. Its remaining kernels
have separated endpoints. The fixed-index and separated-kernel estimates
of Appendix~\ref{v4:kernel-appendix} therefore give an $L^2$ bound for
this factor. In the second term use the global order-one bound
$C_0\theta_i:L^2\to H^{-1}$. Expand each word in
\eqref{v43:lin-EN-polynomial} across a finite chain of collars between
$\supp\theta_o$ and $\supp(1-\eta_B)$. Each summand contains one
separated primitive factor, so Lemmas~\ref{v4:kernel-separated} and
\ref{v4:kernel-multiplier-tails} give
$\theta_oE_N(1-\eta_B):H^{-1}\to L^2$. Consequently
\begin{align*}
 \|B_0\|_{2\to2}\le{}&
 \|\theta_oE_N\|_{2\to2}\|\eta_BC_0\theta_i\|_{2\to2}\\
 &+\|\theta_oE_N(1-\eta_B)\|_{H^{-1}\to L^2}
                      \|C_0\theta_i\|_{L^2\to H^{-1}}.
\end{align*}
Each norm on the right uses finitely many fixed symbol and kernel
seminorms. Taking their maximum gives one finite background index and
the asserted bound $M_0$.
\end{proof}

Use centered Euclidean coordinates $(t,y)$ in the initial-value
plane after the fixed rotation in \eqref{v4:lin-factor}. Choose an
open rectangle $\mathcal O$ containing the supports of $\theta_i,\theta_o$
and embed it by $(t,y)\mapsto(t,y\bmod l)$ in a cylinder $\mathcal X$
of circumference $l$ greater than its transverse width. Denote this
isometric chart by $\iota$. Fix a real $\rho\in C_c^\infty(\mathcal O)$,
$0\le\rho\le1$, equal to one on a neighborhood of these supports and of
$x_*$. For functions on the plane and the cylinder, respectively, put
\begin{equation}\label{v44:scalar-chart-transfer}
 (\mathsf E f)(\iota x)=\rho(x)f(x),\qquad
 (\mathsf Jg)(x)=\rho(x)g(\iota x)\qquad(x\in\mathcal O),
\end{equation}
extending each output by zero. The chart preserves density, so
$\mathsf J=\mathsf E^*$ and $\mathsf J\mathsf E$ is multiplication
by $\rho^2$. These fixed cutoffs and the Euclidean chart give bounded
maps on each Sobolev space, by the product rule, duality and
interpolation. Transfer the complete scalar kernel by setting
\[
 B_{\mathcal X}=\mathsf E B_0\mathsf J.
\]
Since $\rho^2 B_0=B_0=B_0\rho^2$, we have on the whole physical plane
\begin{equation}\label{v44:scalar-chart-intertwining}
 \mathsf J B_{\mathcal X}=B_0\mathsf J .
\end{equation}
In the chart the transferred kernel is exactly the original compact
endpoint kernel of $B_0$. In particular, all contributions from
intermediate integrations of the full plane multipliers remain in it.
Each endpoint cutoff has support inside a single lift of the chart,
so the local amplitude and separated-kernel estimates apply with the
same fixed derivative offsets. Also $\|B_{\mathcal X}\|_{2\to2}\le M_0$.

Extend $\beta$ to a smooth function $\beta_{\mathcal X}$ on the
cylinder by a fixed chart cutoff, equal to $\beta$ near $\supp\rho$
and compactly supported in $t$. This cutoff can first be fixed on the
plane. Take the circumference $l$ large enough that, in a lift of
the chart, the horizontal displacement bound
$\int_{\R}\|\beta_{\mathcal X}(t,\cdot)\|_\infty\,\dd t$ leaves the
original and transported endpoint supports strictly inside one
transverse period. The time support and low norm of the fixed planar
extension give this bound uniformly before $l$ is chosen. Write
$v_{\mathcal X}=\partial_t-\beta_{\mathcal X}\partial_y$. Solve
\[
 \partial_tY(t,y)=-\beta_{\mathcal X}(t,Y(t,y)),\qquad Y(0,y)=y.
\]
The coefficient is bounded and smooth, so this defines a global
cylinder diffeomorphism $X(t,y)=(t,Y(t,y))$. Let
$\Lambda_{\rm tr}$ be a uniform bound for
$\int_{\R}\|\partial_y\beta_{\mathcal X}(t,\cdot)\|_\infty\,\dd t$.
Differentiating the flow equation gives
\[
 e^{-\Lambda_{\rm tr}}\le\partial_yY\le e^{\Lambda_{\rm tr}}.
\]
For the pullback $\mathsf P f=f\circ X$ and its genuine inverse this yields
$\|\mathsf P\|_{2\to2},\|\mathsf P^{-1}\|_{2\to2}
\le e^{\Lambda_{\rm tr}/2}$.
The complete conjugate is $B=\mathsf P B_{\mathcal X}\mathsf P^{-1}$, so
globally
\[
 \mathsf P(v_{\mathcal X}+B_{\mathcal X})\mathsf P^{-1}=\partial_t+B.
\]
In particular $\|B\|_{2\to2}\le e^{\Lambda_{\rm tr}}M_0$.
Proposition~\ref{v43:kernel-conjugated-finite-word} applies to the primitive
factors in \eqref{v43:lin-EN-polynomial} and to the actual difference
$C_0=P-A_0V_0$. It first conjugates these factors and then performs
their finite compositions; the principal cancellation in $C_0$ is
preserved. It gives, for the symbol-and-kernel grading in
Lemma~\ref{v4:lin-scalar-lemma}, a fixed $q_*$ such that
\begin{equation}\label{v43:lin-conjugated-grading}
 |B|_j\le C_j\bigl(1+\|h\|_{H^{j+q_*}}
                         +\|z\|_{H^{j+q_*}}\bigr),\qquad j\ge0.
\end{equation}
Its transported endpoint supports lie in common compact sets.
Include the finite index $b+1+q_0+q_*$ in the low background norm of
the scalar lemma. We can now choose $\tau$ with
$C\tau e^{\Lambda_{\rm tr}}M_0<1/2$ and apply that lemma to the
same $B$ at every order. Write $Q_{\mathcal X}=\mathsf P^{-1}Q_B\mathsf P$.
On the physical plane define
\begin{equation}\label{v44:physical-scalar-solver}
 C=\operatorname{diag}(\Id,\mathsf JQ_{\mathcal X}\mathsf E),\qquad
 D=\operatorname{diag}(\Id,v+B_0).
\end{equation}
On the open set where $\rho$ is identically one near the point,
$v\mathsf Jg=\mathsf Jv_{\mathcal X}g$. Combining this with
\eqref{v44:scalar-chart-intertwining} gives
$(v+B_0)\mathsf JQ_{\mathcal X}\mathsf E f
=\mathsf J(v_{\mathcal X}+B_{\mathcal X})Q_{\mathcal X}\mathsf E f$ there.
On the part of this set with $|t|<\tau$, the last expression equals
$\mathsf J\mathsf E f=f$. Hence for every physical source
\begin{equation}\label{v4:lin-block-inverse}
 DC=\Id
\end{equation}
on this common inner domain. Indeed, $Y(t,\cdot)$
is a circle diffeomorphism for each $t$, and therefore
\begin{equation}\label{v43:lin-common-slab}
 X(( -\tau,\tau)\times\mathbb T)
             =(-\tau,\tau)\times\mathbb T
\end{equation}
for every allowed background. In centered physical-chart coordinates,
let $d_{\rm chart},d_{\rm red}>0$ be the distances from $x_*$ to the
complements of the interior plateau of $\rho$ and of the reduction
platform, respectively. Choose
\[
 r_{\rm mid}<\tfrac14\min\{\tau,d_{\rm chart},d_{\rm red}\}.
\]
The ball $B(x_*,r_{\rm mid})$ is a common domain for
\eqref{v4:lin-block-inverse} and all physical identifications.
The scalar formula \eqref{v44:physical-scalar-solver} defines $C$
on the entire plane. We give its uniform estimates, since both
positive and negative orders enter the vector equation. At each
fixed integer $k\ge0$ the chain rule and change of variables bound
$\mathsf P$ and $\mathsf P^{-1}$ on $H^k$ using the $C^k$ norms of
$X,X^{-1}$ and the low Jacobian bounds. Their adjoints are pullbacks
by the respective inverse maps times the corresponding Jacobians,
so duality gives the $H^{-k}$ estimates with at most one additional
derivative. The low norms needed for each fixed $k$ thus remain
finite. For large $s$, in a derivative of $g\circ X$ of order $s$
the term with $k\ge1$ derivatives on $g$ has factors of derivatives
of $X$ of orders $a_1,\ldots,a_k\ge1$ with $\sum a_j=s$.
The Jacobian bound controls $(D^kg)\circ X$ in $L^2$. Interpolating
the factors of $X$ between $C^1$ and $C^s$ uses total weight
$(s-k)/(s-1)$. The source weight between $H^b$ and $H^s$ is
$(k-b)_+/(s-b)$, and the sum is at most one for $b\ge1$.
The inverse map obeys the same count. The flow estimates in
Appendix~\ref{v4:kernel-appendix} then control both pullbacks with
one high background norm at a fixed offset. Combining these bounds
with \eqref{v4:lin-scalar-tame}, \eqref{v43:lin-conjugated-grading}
and the fixed multiplications in \eqref{v44:scalar-chart-transfer}
gives, for a finite $q_C$ and sufficiently large $s$,
\begin{equation}\label{v44:physical-block-tame}
 \|Cf\|_{H^s}\le C_s\bigl\{\|f\|_{H^s}
   +(1+\|h\|_{H^{s+q_C}}+\|z\|_{H^{s+q_C}})\|f\|_{H^b}\bigr\}.
\end{equation}
In forming each product, the high-source term uses the preceding
factor at order $s$, and the high-coefficient term uses it at the
single low order $b$. Increasing the background low index once for
these fixed orders proves the displayed estimate.
The separated-support estimates below are applied to the
pseudodifferential factors outside $C$.

\subsection{Recombination in the original vector target}
All operators in the following identities act on the initial-value
plane, with $C$ given by \eqref{v44:physical-scalar-solver}. In
particular $M=\mathcal A\mathcal I$ uses the physical injection of
Section~\ref{v4:lin-normalization-section}. Put
\begin{equation}\label{v44:physical-rows-columns}
 \begin{aligned}
 U&=Z\begin{pmatrix}Q&0\\-B_1Q&1\end{pmatrix}
             \binom{\mathcal D}{\mathcal C},&
 V&=\begin{pmatrix}1&-QA_1\\0&1\end{pmatrix},\\
 Z&=\operatorname{diag}(1,iE_N)\operatorname{diag}(1,-\Lambda),&
 D_-&=\operatorname{diag}(1,\Lambda^{-1}).
 \end{aligned}
\end{equation}
Both multipliers in this display are the full plane multipliers of
\eqref{v4:lin-rows}. The finite reductions give the two identities
\begin{equation}\label{v4:lin-matrix-identities}
 UMV=D+S_N,\qquad LU=\Id+T_N.
\end{equation}
To see the first identity, let $\mathsf R_\Pi$ be the last matrix in
\eqref{v4:lin-triangular}, including each term containing $E_\Pi$.
Since $P=-\Lambda(F_1-B_1QA_1)$, \eqref{v4:lin-scalar-finite} yields
\begin{equation}\label{v44:physical-scalar-defect}
 S_N
   =Z\mathsf R_\Pi+
     \operatorname{diag}\bigl(0,\widetilde B-B_0+iR_NV_0\bigr).
\end{equation}
This formula is exact for the chosen global operators. On the common
inner cutoffs, $\mathsf R_\Pi$ has order at most $-d$, so the first
term has order at most $1-d$. The last term has the same order from
$R_N$. The difference $\widetilde B-B_0$ has zero kernel when both
endpoints are in the common interior of the two cutoff plateaux.
Its exterior input and output terms are given by
\eqref{v4:lin-B-localization}. Thus the local diagonal order of
$S_N$ is at most $-N$ for $d=N+1$, and its raw global order is at
most one.

For the second identity set $W_0=D_-U$. On the inner collar its
order-zero principal symbol is
\[
 \operatorname{diag}(1,-i/a_0)
 \begin{pmatrix}1&0\\-\sigma_0(B_1)&1\end{pmatrix}
 \begin{pmatrix}i e^{\mathsf T}\\ i(e^\perp)^{\mathsf T}\end{pmatrix},
 \qquad e=\xi/|\xi| .
\]
Indeed $\sigma_0(Q)=1$, $\sigma_0(E_N)=1/a_0$, and the last
matrix consists of the divergence and curl symbols. Its rows are
$i$ times an orthonormal pair. The ellipticity and low bounds of
$a_0$ therefore give an inverse symbol with uniform bounds on the
collar. Multiply this inverse symbol by a fixed high-frequency
excision and localize it with spatial cutoffs supported in that
elliptic domain and equal to one on a smaller collar. Left
quantization gives an order-zero $J_0$ on the plane. With
$R_J=J_0W_0-\Id$, set
\[
 J_N=\sum_{j=0}^{N-1}(-R_J)^jJ_0,\qquad
 L=J_ND_-,\qquad T_N=-(-R_J)^N.
\]
The finite geometric sum proves the second identity on the whole
plane. Its defect $T_N$ has raw global order zero and local diagonal
order at most $-N$, since $R_J$ has order at most $-1$ on the
elliptic collar. Appendix~\ref{v4:kernel-appendix} gives the estimates
also for the separated kernels of these products and the full
multiplier tails.

The approximate physical control is
$\mathcal B^{\rm app}=\mathcal IVCU$.  Multiplying the first
identity of~\eqref{v4:lin-matrix-identities} by $L$ and then by $CU$
gives the useful global algebraic identity
\begin{equation}\label{v4:lin-recombination}
 MVCU-\Id=L(DC-\Id)U+T_N+(LS_N-T_NMV)CU.
\end{equation}
Indeed $LUMV=(\Id+T_N)MV=L(D+S_N)$. Right multiplication by
$CU$ gives \eqref{v4:lin-recombination} after adding and subtracting
$LU$. Every product in this calculation is a composition on the
plane. All cutoff terms from the chart transfer of the scalar solve
are contained in the actual operator $DC-\Id$. Its restriction to
the inner domain vanishes by \eqref{v4:lin-block-inverse}.

We seek a target gain of $3/2$ and take $N=4$; the finite pivot and
scalar depths above are then $N+1=5$.  There are only finitely many primitive operator
words.  If their largest length is $p$, choose a nested chain of
$p+4$ compact collars strictly inside $B(x_*,r_{\rm mid})$.
The gaps
are fixed and positive.  Let $\chi$ be an inner output cutoff and
$\sigma$ a larger source cutoff, both equal to one on a common
domain $\Omega_0$.  They can be chosen so that $\sigma=1$ near
$\supp\chi$, and so that all later corrected source supports lie
inside $\Omega_0$.  Let $\kappa$ equal one on all these collars.
Its support is contained in the common interior of the plateaux
$\theta_o=\theta_i=1$ and every region where the finite symbol
cancellations hold.  Choose a further cutoff $\eta$, equal to one
near $\supp\kappa$, whose support is still in this common platform.
Let $\kappa_t$ equal one on the collars and have support in the
common domain where $DC=\Id$.  All these cutoffs take values in $[0,1]$.
Define the actual localized target error by
\begin{equation}\label{v4:lin-R-definition}
 R=\chi(MVCU-\Id)\sigma.
\end{equation}
On $\Omega_0$, for every source supported in its prescribed larger
source collar, the physical equation is
\begin{equation}\label{v4:lin-approximate-identity}
 \mathcal M_H\mathcal B^{\rm app}f=f+Rf.
\end{equation}
Indeed $M=\chi_A\mathcal M_H\mathcal I$ on the plane, $\chi_A=\chi=1$
on $\Omega_0$ and $\sigma f=f$ for these sources. In particular
\eqref{v4:lin-approximate-identity} involves the same physical
control $\mathcal IVCU f$ as the exact operator product.

\begin{lemma}[Complete vector error]\label{v4:lin-complete-error}
The actual error in~\eqref{v4:lin-R-definition} satisfies, uniformly
on one fixed finite low neighborhood,
\begin{equation}\label{v4:lin-low-gain}
 \|R\sigma_0f\|_{H^{3/2}}\le C\|f\|_{L^2},
\end{equation}
where $\sigma_0$ has support in the larger source collar and is one
near $\overline{\Omega_0}$.  For fixed finite $b,q_1$ it also satisfies
\begin{equation}\label{v4:lin-high-gain}
 \|Rf\|_{H^s}\le C_s\left[\|f\|_{H^{s-3/2}}
 +(1+\|h\|_{H^{s+q_1}}+\|z\|_{H^{s+q_1}})\|f\|_{H^b}\right].
\end{equation}
Neither $R$ nor these low constants has a negative power of
$\varepsilon$.
\end{lemma}
\begin{proof}
The identity~\eqref{v4:lin-recombination} splits $R$ exactly into
\begin{align*}
 R_1&=\chi L(1-\kappa_t)(DC-\Id)U\sigma,\\
 R_2&=\chi\kappa T_N\sigma,\\
 R_3&=\chi L(\kappa S_N)CU\sigma,\\
 R_4&=\chi L(1-\kappa)S_NCU\sigma,\\
 R_5&=-\chi\kappa T_NMVCU\sigma.
\end{align*}
Here $\chi\kappa=\chi$ and $\kappa_t(DC-\Id)=0$.
The global orders, before local cancellation, are
at most $1$ for $S_N=UMV-D$ and $0$ for $T_N=LU-\Id$.
Consequently the fixed-index mapping bounds give
\[
 U:L^2\longrightarrow H^{-1},\quad
 C:H^{-1}\longrightarrow H^{-1},\quad
 D:H^{-1}\longrightarrow H^{-2},
\]
while $M,V,L$ preserve each fixed Sobolev order.  The constants use
only finitely many low background seminorms.  We will specify one
common choice after listing the required paths.

First retain the two pieces of each localized defect:
\begin{equation}\label{v41:lin-local-defect-split}
 \begin{split}
 \kappa S_N&=\kappa S_N\eta+\kappa S_N(1-\eta),\\
 \kappa T_N&=\kappa T_N\eta+\kappa T_N(1-\eta).
 \end{split}
\end{equation}
The first terms have order $-4$ by the finite cancellation on the
common platform.  In the second terms the input and output are
separated.  Each is a finite sum of localized words in the primitive
pseudodifferential factors; neither contains $C$.

For these terms and for exterior outputs, use the exact collar
expansion of a word
$L_1\cdots L_a$ in $L$.  If $\eta_0=\chi$ and
$\eta_j=1$ near $\supp\eta_{j-1}$, with
$\eta_a(1-\kappa)=0$, then
\begin{align*}
 \eta_0L_1\cdots L_a(1-\kappa)
 =\sum_{j=1}^a &\eta_0L_1\eta_1\cdots
       L_{j-1}\eta_{j-1}L_j(1-\eta_j)
       L_{j+1}\cdots L_a(1-\kappa).
\end{align*}
In each summand a primitive pseudodifferential factor crosses two
separated supports.  Estimate the factors to its right at a fixed
negative index.  Lemma~\ref{v4:kernel-separated} maps that index to
the desired output index plus the sum of the orders on its left;
the remaining left factors then have fixed-index bounds.  There
are enough collars because their number was fixed using the largest
word length.  Apply the same expansion with $\kappa_t$, or with
$\kappa,\eta$ in~\eqref{v41:lin-local-defect-split}.  If the
separated primitive is a full $\Lambda^{\pm1}$ multiplier, its
input is not necessarily compact.  In this case
Lemma~\ref{v4:kernel-multiplier-tails} supplies the weighted kernel
bound before the remaining finite factors are applied.  Thus the
argument includes the full multiplier tails as well as compact
endpoint kernels.

We can now estimate the five terms. For $R_1$, the right factor $(DC-\Id)U$ maps
$L^2$ into $H^{-2}$, and
$\chi L(1-\kappa_t):H^{-2}\to H^{3/2}$ by the separated-word
estimate.  For $R_2$, both pieces of $\kappa T_N$ map
$L^2$ into $H^4$; the diagonal piece uses order $-4$, and the
separated piece uses the preceding estimate at these fixed indices.
For $R_3$, the path is
\[
 L^2\xrightarrow{\,U\,}H^{-1}
 \xrightarrow{\,C\,}H^{-1}
 \xrightarrow{\,\kappa S_N\,}H^3
 \xrightarrow{\,\chi L\,}H^3.
\]
For $R_4$, the global order-one bound for $S_N$ gives
$S_NCU:L^2\to H^{-2}$, after which
$\chi L(1-\kappa):H^{-2}\to H^{3/2}$ applies.
Finally, $MVCU:L^2\to H^{-1}$ and
$\kappa T_N:H^{-1}\to H^3$ give the bound for $R_5$.
All five outputs therefore lie in $H^{3/2}$.

For each of this finite set of paths, the fixed-index estimate
\eqref{v4:kernel-fixed-index}, the scalar inverse lemma and the
separated-kernel lemmas require finitely many symbol and kernel
seminorms.  Take their maximum, including the fixed orders of the
left and right subwords.  By~\eqref{v4:kernel-triangular}, one
finite low Sobolev index of $(h,z)$ controls this entire list.
Restrict the original low ball at that index once.  The five bounds
then share one constant, independent of $\varepsilon$ and of a later
high Sobolev order.  Summation proves~\eqref{v4:lin-low-gain} for
smooth $L^2$ sources, and density gives the stated extension.

For the high estimate, take a common finite $b,q_1$ and abbreviate
$\mathfrak H_s=1+\|h\|_{H^{s+q_1}}+\|z\|_{H^{s+q_1}}$.
The same five paths give, for all sufficiently large integers $s$,
\begin{equation}\label{v43:lin-error-high-terms}
 \begin{aligned}
 \|R_1f\|_{H^s}+\|R_4f\|_{H^s}
       &\le C_s\mathfrak H_s\|f\|_{H^b},\\
 \|R_2f\|_{H^s}
       &\le C_s\bigl(\|f\|_{H^{s-4}}
                         +\mathfrak H_s\|f\|_{H^b}\bigr),\\
 \|R_3f\|_{H^s}+\|R_5f\|_{H^s}
       &\le C_s\bigl(\|f\|_{H^{s-3}}
                         +\mathfrak H_s\|f\|_{H^b}\bigr).
 \end{aligned}
\end{equation}
Here is the high-norm accounting for these inequalities. In $R_1$ and
$R_4$, keep respectively $(DC-\Id)U\sigma$ and $S_NCU\sigma$ at
$L^2\to H^{-2}$, with uniform low constants. Expand the left word
across the collars as above. The separated factor maps this fixed
negative order to $H^s$ plus the sum of the orders to its left.
The tame product estimate for that left word puts the high norm on
one coefficient; all other input norms are fixed low norms of the
separated output. This proves the first line. It applies to complete
multiplier tails by Lemma~\ref{v4:kernel-multiplier-tails}.

For $R_2$, the diagonal block in
\eqref{v41:lin-local-defect-split} gains four derivatives, while its
separated block has the preceding all-orders bound. For $R_3$ put
$A=\kappa S_N\eta$ and $g=CU\sigma f$. Apply the tame estimate
for the diagonal operator of order $-4$ and then
\eqref{v44:physical-block-tame} and \eqref{v4:lin-word-bound}. For
fixed sufficiently large $b_1$ and $b\ge b_1+1$ they give
\[
 \begin{aligned}
 \|Ag\|_{H^s}
  &\le C_s\bigl(\|g\|_{H^{s-4}}+
                         \mathfrak H_s\|g\|_{H^{b_1}}\bigr),\\
 \|g\|_{H^{s-4}}
  &\le C_s\bigl(\|f\|_{H^{s-3}}+
                         \mathfrak H_s\|f\|_{H^b}\bigr),
 &\|g\|_{H^{b_1}}&\le C\|f\|_{H^b}.
 \end{aligned}
\]
Here the last inequality uses only fixed-order operator bounds, with
the required background indices included in the low neighborhood.
The separated part $\kappa S_N(1-\eta)$ uses $g$ only in $H^{-1}$.
The order-zero factor $\chi L$ preserves the resulting bound by
\eqref{v4:lin-word-bound}, using the preceding low-order paths for
its low source term. For $R_5$, take $A=\kappa T_N\eta$ and
$g=MVCU\sigma f$ in the same inequalities. The two additional
factors $M,V$ have order zero. Thus the loss in $g$ is again one,
and the separated part again uses $g$ in $H^{-1}$. Taking the maximum
of the finite offsets and low indices gives $q_1,b$ for all five terms.

These estimates imply \eqref{v4:lin-high-gain}; the remaining finite
low levels follow from the same fixed-index paths and interpolation.
All factors in the target equation are normalized: the identity
$\mathcal A\mathcal I=M$ has canceled the injection's
$\varepsilon^{-1}$. Hence the target-error bounds are uniform as
$\varepsilon$ tends to zero within the background neighborhood.
\end{proof}

\subsection{Eliminating the full target error}
The last localization uses the support of the output to make a
regularizing error small in $L^2$, as in \cite{Dencker}.
Fix all operators and collars in $R$, and the bound in
\eqref{v4:lin-low-gain}. Choose an output cutoff $\chi_{\rm ex}$ supported in a ball of radius
$2\delta_{\rm ex}$ contained in $\Omega_0$, equal to one on its
concentric ball of radius $\delta_{\rm ex}$.  Arrange
$\sigma_0=1$ near its entire support.  Set
\begin{equation}\label{v4:lin-final-G}
 K_{\rm ex}=\chi_{\rm ex}R\sigma_0,\qquad
 G=(\Id+K_{\rm ex})^{-1}_{L^2}.
\end{equation}
The inverse is justified by the two-dimensional Sobolev embedding:
\[
 \|K_{\rm ex}f\|_2
 \le |\supp\chi_{\rm ex}|^{1/2}\|R\sigma_0f\|_\infty
 \le C\delta_{\rm ex}\|f\|_2.
\]
Choose
\[
 2\delta_{\rm ex}<\operatorname{dist}(x_*,\partial\Omega_0),
 \qquad C\delta_{\rm ex}<\tfrac12.
\]
This bound uses the support area and the already fixed norm of $R$,
so the same radius works for every permitted background.

Inside the output plateau choose the physical domains and source
cutoff with the following inclusions, understood in the fixed chart:
\begin{equation}\label{v43:lin-domain-chain}
 \begin{gathered}
 \overline{\Omega_-}\subset\{\zeta=1\},\qquad
 \supp\zeta\Subset\Omega_+\Subset B(x_*,\delta_{\rm ex})
                         \subset\{\chi_{\rm ex}=1\},\\
 \supp\chi_{\rm ex}\Subset\Omega_0,\qquad
 \sigma_0=1\ \hbox{near }\overline{\Omega_0},\qquad
 \supp\sigma_0\Subset\{\sigma=1\}.
 \end{gathered}
\end{equation}
Thus an arbitrary smooth source on $\Omega_+$ is first multiplied by
$\zeta$ and extended by zero on the initial-value plane. The source fed to $G$ is supported in
the plateau of $\chi_{\rm ex}$.

Let $f$ have support in the plateau of $\chi_{\rm ex}$ and put
$g=Gf$.  The equation
\begin{equation}\label{v4:lin-G-support}
 g=f-\chi_{\rm ex}R\sigma_0g
\end{equation}
shows $\supp g\subset\supp\chi_{\rm ex}$.  Hence $\sigma_0g=g$,
and throughout the output plateau
\[
 \mathcal M_H\mathcal B^{\rm app}g=g+Rg=f.
\]
The approximate identity is valid for the entire support of $g$,
including the transition region of $\chi_{\rm ex}$.  This is why its
source collar was fixed before the final output localization.

We justify the smoothness of this same inverse.  The low estimate
gives $K_{\rm ex}:L^2\to H^1$.  At the one finite level $b+1$,
the high estimate gives $K_{\rm ex}:H^b\to H^{b+1}$ after the
finitely many needed coefficients have been included in the low ball.
Interpolation supplies the intermediate integer gains.  Starting
with the actual $L^2$ solution in~\eqref{v4:lin-G-support} therefore
gives $\|g\|_{H^b}\le C\|f\|_{H^b}$; further positive gains
give smoothness for smooth $f$.  Finally,
\[
 \|g\|_{H^s}\le\|f\|_{H^s}+C_s\|g\|_{H^{s-3/2}}
 +C_s(1+\|h\|_{H^{s+q_1}}+\|z\|_{H^{s+q_1}})\|g\|_{H^b}.
\]
For $s>b+3/2$, interpolate that middle norm between $H^b$ and $H^s$
with weight $\alpha=(s-3/2-b)/(s-b)<1$. Young's inequality gives
\[
 C_s\|g\|_{H^{s-3/2}}\le\tfrac12\|g\|_{H^s}
                   +C'_{s,b}\|g\|_{H^b}.
\]
For $b\le s\le b+3/2$ the middle norm is already bounded by its
$H^b$ norm. Inserting the low estimate for $g$ in both cases gives
\begin{equation}\label{v4:lin-G-tame}
 \|Gf\|_{H^s}\le C_s\left[\|f\|_{H^s}
 +(1+\|h\|_{H^{s+q_2}}+\|z\|_{H^{s+q_2}})\|f\|_{H^b}\right].
\end{equation}
The definition of $G$ and the domains have not changed with $s$.

\subsection{Control norms and completion}
It remains to estimate the physical injection.  In
\eqref{v4:lin-adjoint}, change variables $(a,q)=\widetilde\Psi_s(x)$.
With $X_s=\widetilde\Psi_s^{-1}$, the adjoint is
\begin{equation}\label{v4:lin-coarea}
 (\mathcal A^*v)(q,s)
 =\int\frac{\chi_A(X_s(a,q))w_s(X_s(a,q))\cdot v(X_s(a,q))}
              {|\det D_x\widetilde\Psi_s(X_s(a,q))|}\,\dd a.
\end{equation}
The base cutoff has interior support.  Its transported support
stays a positive distance from the outer image boundaries, so the
integrand extends by zero smoothly on fixed larger charts.  No
derivative of a moving integration boundary occurs.  On each time
slice the change of variables and its inverse have bounded Jacobians;
hence
\[
 \|(D^kv)\circ X_s\|_{L^2(\dd s\,\dd a\,\dd q)}
                         \le C\|v\|_{H^k}.
\]
Cauchy--Schwarz in the bounded $a$ interval controls the final integral.

In a derivative of total order $n$, the integrand is a sum of terms
$(D^{\alpha_0}a_0)((D^kv)\circ X)\prod_{j=1}^kD^{\alpha_j}X$,
with $|\alpha_j|\ge1$ for $j\ge1$ and
$|\alpha_0|+\sum_j|\alpha_j|=n$.  If $k=0$, the term is
$(D^{\alpha_0}a_0)(v\circ X)$ and is bounded directly by the
coefficient norm of order $s$ times $\|v\|_{L^2}$, using the
Jacobian bound.  For $k\ge1$ and $n\le s$, coefficient
interpolation between fixed $C^1$ and $C^s$ bounds has total weight
at most $(s-k)/(s-1)$.  Source interpolation between $H^b$ and $H^s$,
$b\ge2$, has weight $(k-b)_+/(s-b)$.  Their sum is at most one.
The product is therefore controlled by one high source norm or one
high coefficient norm times a low source norm.  Flow and inverse-flow
estimates give, with a fixed $q_3$,
\begin{equation}\label{v4:lin-adjoint-tame}
 \|\mathcal A^*v\|_{H^s(W)}\le C_s\left[\|v\|_{H^s}
 +(1+\|h\|_{H^{s+q_3}}+\|z\|_{H^{s+q_3}})\|v\|_{H^b}\right].
\end{equation}
The filters $C_\pm$ are bounded
on the joint Sobolev scale, $\mathcal D^*$ has order zero and
$Q_\pm$ has order one.  Thus $\mathcal I$ loses one source
derivative and contributes $|\varepsilon|^{-1}$.

Let $\mathsf E_\zeta$ denote multiplication by the cutoff in
\eqref{v43:lin-domain-chain}, followed by extension by zero from
$\Omega_+$ to the initial-value plane.
The final right inverse is
\begin{equation}\label{v4:lin-final-inverse}
 \mathcal B_{\varepsilon,h,z}
     =\operatorname{Re}\bigl(\mathcal IVCUG\mathsf E_\zeta\bigr).
\end{equation}
The physical comparison operator is real, so the real part preserves
its equation.  The leftmost velocity cutoff in $\mathcal I$ ensures
the fixed compact control support.  The gains and losses, read from
right to left, are $0,1,0,0,1$ for $G,U,C,V,\mathcal I$.
The finite-product estimate gives the net loss two and
\eqref{v4:lin-tame}.  The exact equation was proved on the final
output plateau, yielding~\eqref{v4:lin-exact}.

The target factors are uniformly bounded in the normalized background
scales, and $\mathcal I$ supplies the single factor
$|\varepsilon|^{-1}$ in \eqref{v4:lin-tame}. The conversion to physical
control norms will be made explicit in
Lemma~\ref{v43:newton-parameter-budget}.

The order of the choices is finite.  First choose the reference
characteristic and flow margins, the finite depths, and the coefficient
extensions.  Bound the low norm of the complete scalar $B$ and choose
the transport interval.  Next choose the finite collar chain, bound
the low target-error norm and choose the final output cutoff.  Finally
increase $b_0$ to include the finite indices for
Lemma~\ref{v43:schur-continuity}, the all-pairs twist, the scalar
$L^2$ bound, the scalar grading at $b+1+q_0$, the five low error
paths, and the bootstrap for $G$ through $b+1$. Restrict the original
low ball at this maximum. Each restriction preserves the previous
bounds. The domains and this one neighborhood then work at every
high Sobolev order and for every later finite formal truncation.
This completes the proof of Theorem~\ref{v4:local-inverse}.
\part{Exact return and geometric realization}\label{part:realization}
\section{Smoothing iteration for the actual return map}
\label{v4:newton-iteration}
\label{v4:newton-chapter}

The inverse in Theorem~\ref{v4:local-inverse} belongs to the periodic
comparison operator. The equation to be solved still concerns the original
oscillator. We write
\[
 N(H)=P_H-\Id,
 \qquad DN_H=DP_H\mathcal M_H^{\mathrm{tr}}.
\]
Every correction below is a real, smooth, $2\pi$-periodic control on one
fixed cylinder. Only the initial-value region on which the error is measured
shrinks. We use the uniform value estimates and fixed losses established
above. Smoothing in a rapidly convergent iteration goes back to
Nash and Moser \cite{Nash1956,Moser1966I}. For the estimates on scales
of spaces and the use of approximate inverses, see
\cite{Hamilton,Zehnder1975}.

\subsection{Fixed spaces and the linear estimate used in the iteration}
\label{v4:newton-scales}

Fix a seed in the class of Definition~\ref{v4:seed-class-B}, a corresponding
nonzero point $x_*$, and for the moment a nonzero parameter $\varepsilon$.
Let $\Omega_-\Subset\Omega_+$ be the common initial-value neighborhoods of
Theorem~\ref{v4:local-inverse}, chosen so that their closures avoid the
origin. A slightly larger flow neighborhood is fixed as well. All trajectories
under consideration exist for $0\le s\le2\pi$ there and their $q$ values
remain in a fixed compact subinterval of $(-V,V)$.
Put
\[
 \mathcal W=(-V,V)\times\R/(2\pi\mathbb Z).
\]
The compact velocity support of the inverse output lies strictly inside this
cylinder. All control norms below are taken on the entire cylinder.

Let $H_0=H_{0,\varepsilon,h}$ be a finite approximation, to be chosen in
Section~\ref{v4:newton-convergence}. Write $H=H_0+z$. A finite integer
$b\ge6$ is chosen once, large enough to include every low source index and
background index in the linear theorem and a prescribed final $C^k$ norm.
In particular it includes the low \emph{source} index plus the source loss
in that theorem. Define the two scales
\begin{equation}\label{v4:newton-norms}
 |z|_j=\|z\|_{C^{b+4+j}(\mathcal W)},
 \qquad \|f\|_{j,U}=\|f\|_{C^{b+j}(U)},\qquad j\ge0.
\end{equation}
An unlabelled control norm always means the full-cylinder norm.
The four additional derivatives in the control scale cover the first two
flow variations and the endpoint comparison at the low target level.
Increasing $b$ to cover a prescribed finite $k$ does not reselect the
inverse or its domains. It only restricts the iteration to a smaller ball
measured in a stronger finite norm inside the same inverse neighborhood.
The small parameter will place the chosen reference in that ball.

After converting the Sobolev estimates of Theorem~\ref{v4:local-inverse},
there is a fixed integer $d$ such that
\begin{align}
 \mathcal M_H\mathcal B_Hg&=g\quad\hbox{on }\Omega_-,
 \label{v4:newton-inverse-identity}\\
 |\mathcal B_Hg|_j
 &\le K_j(\varepsilon)
 \left(\|g\|_{j+d,\Omega_+}
       +(1+|z|_{j+d})\|g\|_{0,\Omega_+}\right)
 \label{v4:newton-inverse-bound}
\end{align}
for smooth sources supported in a fixed compact subset of $\Omega_+$.
The higher norms of the chosen $h$ and $H_0$ may enter $K_j$. The integer $d$
does not depend on $j$ or on the iteration number. To see the shift explicitly,
denote the unshifted Sobolev loss by $r_{\mathrm{lin}}$ and its background
offset by $q$. Apply the linear estimate at order $b+4+j+2$ and use
$H^{m+2}\hookrightarrow C^m$ on the two-dimensional control cylinder.
It is sufficient to take
\[
 d\ge\max\{r_{\mathrm{lin}}+6,q+2,1\},
\]
and then enlarge $d$ once to cover any remaining fixed offsets.
All low source Sobolev norms are covered by the original choice of $b$.

\begin{lemma}[Parameter bounds in physical control norms]
\label{v43:newton-parameter-budget}
For a fixed finite approximation of degree $M$, the admissible physical
ball can be taken to have radius $R_\varepsilon=c_0|\varepsilon|$,
where $c_0>0$ is fixed. On this ball, for every $j\ge0$,
\begin{equation}\label{v43:newton-inverse-power}
 K_j(\varepsilon)\le C_{j,h,M}|\varepsilon|^{-2}.
\end{equation}
The exponent and $c_0$ can be chosen before $M$; the constants and the
upper bound on $|\varepsilon|$ may depend on the chosen approximation.
\end{lemma}
\begin{proof}
Use $u$ for the normalized correction in the linear theorem, so that
$H=\varepsilon(hq^2+u)$. The finite approximation has $h_2=h$, $h_3=0$,
and
\begin{equation}\label{v43:newton-normalized-reference}
 \begin{split}
 H_0&=\sum_{m=2}^{M}\varepsilon^{m-1}\gamma_m(s)q^m,\\
 u&=u_{0,\varepsilon}+\varepsilon^{-1}z,
 \qquad u_{0,\varepsilon}
       =\sum_{m=4}^{M}\varepsilon^{m-2}\gamma_m(s)q^m.
 \end{split}
\end{equation}
These formulas hold on the control cylinder; the fixed outer velocity
cutoff equals one there. For each fixed $r$,
$\|u_{0,\varepsilon}\|_{H^r}\le C_{r,h,M}|\varepsilon|^2$.
Substitute this in \eqref{v4:lin-tame}. For $0<|\varepsilon|\le1$,
its right side is bounded by
\[
 C_{s,h,M}|\varepsilon|^{-2}
 \left[\|g\|_{H^{s+2}}+
       (1+\|z\|_{H^{s+q}})\|g\|_{H^{b_{\rm lin}+2}}\right],
\]
where $b_{\rm lin}$ denotes the low source index of that theorem.
The shift of scales just described proves
\eqref{v43:newton-inverse-power} with the same exponent for all $j$.

Let $E_b$ satisfy $\|z\|_{H^{b_0}}\le E_b|z|_0$ and let $\delta$
be the normalized background radius at the fixed seed. First decrease
$|\varepsilon|$ until
$\|u_{0,\varepsilon}\|_{H^{b_0}}<\delta/4$. Choose
$c_0\le\delta/(4E_b)$. Then $|z|_0<c_0|\varepsilon|$ places $u$
in the $\delta/2$ ball. Decreasing $c_0$ once also meets the fixed
physical flow-tube bounds. This choice depends on the seed and the
previously chosen indices, but not on $M$.
\end{proof}

The matrix $(DP_H(x))^{-1}$ in what follows is inverted pointwise in $x$.
The variational ODE and the matrix inverse formula give its high--low
product estimate on the same background ball. Multiplication by this
matrix, and by one fixed source cutoff, can therefore be included in
\eqref{v4:newton-inverse-bound} by increasing $d$ once. Since this
estimate uses the physical control $H$, this multiplication preserves
the parameter exponent in \eqref{v43:newton-inverse-power}.

\subsection{Differentiability and estimates for the true return}
\label{v4:newton-flow-estimates}

\begin{lemma}[Estimates on a common flow tube]\label{v4:newton-tame}
On the protected low background ball, the actual return map is twice
differentiable at every finite target order. In the scales
\eqref{v4:newton-norms}, for every integer $j\ge0$,
\begin{align}
 \|N(H_0+z)\|_{j,U}
 &\le C_j(1+|z|_j),\label{v4:newton-value}\\
 \|DN_{H_0+z}\eta\|_{j,U}
 &\le C_j\bigl(|\eta|_j+(1+|z|_j)|\eta|_0\bigr),
 \label{v4:newton-first}\\
 \|D^2N_{H_0+z}[\eta,\zeta]\|_{j,U}
 &\le C_j\bigl(|\eta|_j|\zeta|_0+|\zeta|_j|\eta|_0
                +(1+|z|_j)|\eta|_0|\zeta|_0\bigr).
 \label{v4:newton-second}
\end{align}
Here $U$ may be any of the fixed inner flow neighborhoods. The constants
require only the protected low norm of the variable background and
the indicated finite norms of $H_0$. For each fixed finite approximation
they are bounded uniformly as $\varepsilon\to0$. In addition,
\begin{equation}\label{v4:newton-low-variations}
 \|DN_H\eta\|_{1,U}\le C|\eta|_0,
 \qquad \|D^2N_H[\eta,\zeta]\|_{0,U}
       \le C|\eta|_0|\zeta|_0.
\end{equation}
\end{lemma}

\begin{proof}
The trajectory equation is the integral equation for
$V_H(p,q,s)=(q,-p+H(q,s))$. At each finite initial-value order $m$, its
linearization is a Volterra operator on
$C([0,2\pi],C^m(\overline U;\R^2))$. If its coefficient has finite
multiplication bound $M_m$, the $n$th Volterra iterate has norm at most
$(2\pi M_m)^n/n!$. The inverse of the linearized integral equation is
therefore the convergent Volterra series. The ordinary product rule and
Taylor formula for $H(q,s)$ show that the integral equation is twice
differentiable when the control has $m+2$ derivatives. This proves actual
control differentiability and identifies its derivatives with the following
variational equations.

Writing $Y_\eta=D_H\Psi_s[\eta]$, the first variation satisfies
\begin{equation}\label{v4:newton-variation-ode}
 \partial_sY_\eta
 =\begin{pmatrix}0&1\\-1&H_q(q(x,s),s)\end{pmatrix}Y_\eta
       +\binom0{\eta(q(x,s),s)},\qquad Y_\eta(x,0)=0.
\end{equation}
The mixed second variation has the same homogeneous coefficient and zero
initial value. Its forcing in the second component is
\begin{equation}\label{v4:newton-second-forcing}
 H_{qq}(Y_\eta)^2(Y_\zeta)^2
       +\eta_q(Y_\zeta)^2+\zeta_q(Y_\eta)^2.
\end{equation}
The superscript $2$ here denotes the second component, not its square;
all control coefficients are evaluated at $(q(x,s),s)$.

For clarity, we verify that these equations give estimates linear in the
high norm. Put $B=b+4$ and $[v]_r=\|v\|_{C^r(\mathcal W)}$.
Normalize up to two control directions $e_\nu$ by $[e_\nu]_B\le1$ and
adjoin constant parameters $a_\nu$ to the state. The augmented vector
field is
\[
 \mathcal V(p,q,a,s)
 =\left(q,-p+H(q,s)+\sum_\nu a_\nu e_\nu(q,s),0\right).
\]
All estimates are evaluated at $a=0$. If
$Q_K=1+[H]_K+\sum_\nu[e_\nu]_{K-1}$, interpolation between the fixed
low norms and order $K$ gives
\[
 |D^r\mathcal V|\le C_K Q_K^{(r-1)/(K-1)},\qquad 1\le r\le K.
\]
The first derivative is bounded independently of $Q_K$. Let $J_j$ be a
$j$th derivative of the augmented flow in the initial conditions and
parameters. Its highest derivative occurs in the term
$D\mathcal V\,J_j$. Every other term is
\[
 D^r\mathcal V[J_{i_1},\ldots,J_{i_r}],
 \qquad r\ge2,\quad i_\nu\ge1,\quad\sum_\nu i_\nu=j.
\]
The sum of the interpolation exponents is
\[
 \frac{r-1}{K-1}+\sum_\nu\frac{i_\nu-1}{K-1}
       =\frac{j-1}{K-1}.
\]
Induction and variation of constants with the same low coefficient bound
give $|J_j|\le C_{K,j}Q_K^{(j-1)/(K-1)}$. In particular the top-order
bound contains only one high norm, and its Gronwall factor uses only low
derivatives. Reversing the normalization of the directions yields
\begin{align*}
 \|\Psi\|_{C^m_x}&\le C_m(1+[H]_m),\\
 \|Y_\eta\|_{C^m_x}&\le C_m\bigl([\eta]_m+(1+[H]_{m+1})[\eta]_B\bigr),\\
 \|D_H^2\Psi[\eta,\zeta]\|_{C^m_x}
 &\le C_m\bigl([\eta]_{m+1}[\zeta]_B+[\zeta]_{m+1}[\eta]_B
             +(1+[H]_{m+2})[\eta]_B[\zeta]_B\bigr).
\end{align*}
These bounds are uniform over the finite time interval. Evaluation at
$s=2\pi$ and $m=b+j$ proves the displayed estimates. Taking $m=b+1$
for the first variation and $m=b$ for the second proves
\eqref{v4:newton-low-variations}. All needed low background derivatives
then have order at most $b+2<B$.

In particular, whenever its whole segment remains in the ball, Taylor's
formula holds in each finite target space as
\begin{equation}\label{v4:newton-taylor}
 N(H+\eta)=N(H)+DN_H\eta
       +\int_0^1(1-t)D^2N_{H+t\eta}[\eta,\eta]\dd t.
\end{equation}
The integral is an integral of continuous functions with values in that
Banach space. Its low norm is bounded by $C|\eta|_0^2$.
\end{proof}

The endpoint comparison from Proposition~\ref{v4:comparison-defect}, in these
scales, reads
\begin{equation}\label{v4:newton-comparison}
 \|DP_H(\mathcal M_H^{\mathrm{tr}}-\mathcal M_H)\eta\|_{0,U'}
       \le C\|N(H)\|_{1,U}|\eta|_0.
\end{equation}
The inverse interpolations defining the comparison must send $U'$ into
$U$. The additional derivative of the residual in this estimate is kept
explicit; it prevents replacing this step by a same-level Neumann series.

\subsection{Nested balls, extension, and smoothing}
\label{v4:newton-domains}

Choose $\rho_\infty>0$ and $c_{\rm gap}>0$ so that
\begin{equation}\label{v4:newton-balls}
 U_n=B(x_*,\rho_n),\qquad
 \rho_n=\rho_\infty+c_{\rm gap}\sum_{i=n}^\infty(i+1)^{-2},
 \qquad\overline U_0\subset\Omega_-.
\end{equation}
They can be chosen independently of $\varepsilon$ for each fixed seed
and its inverse domains. The gap is
$\delta_n=\rho_n-\rho_{n+1}=c_{\rm gap}(n+1)^{-2}$. Divide each gap into
eight equal collars. Choose a smooth cutoff $\chi_n$ supported in $U_n$
and equal to one on a neighborhood of
$B(x_*,\rho_{n+1}+6\delta_n/8)$. Its derivatives obey
$\|D^r\chi_n\|_\infty\le C_r\delta_n^{-r}$.
If $f$ is smooth on $U_n$, its cut-off zero extension is smooth on the
fixed target space, and
\begin{equation}\label{v4:newton-extension}
 \|\chi_nf\|_{j,\R^2}
       \le C_j(n+1)^{m_j}\|f\|_{j,U_n}.
\end{equation}
All derivatives of the cutoff lie where the original residual is defined.

Use fixed real Fourier smoothers $T_\theta$ on the periodic control
cylinder after zero extension in $q$, and $S_\theta$ on the target plane.
They may be chosen so that, for $j\le J$,
\begin{align}
 |T_\theta z|_J&\le C_{j,J}\theta^{J-j}|z|_j,&
 |(I-T_\theta)z|_j&\le C_{j,J}\theta^{j-J}|z|_J,
 \label{v4:newton-control-smoothing}\\
 \|S_\theta f\|_J&\le C_{j,J}\theta^{J-j}\|f\|_j,&
 \|(I-S_\theta)f\|_j&\le C_{j,J}\theta^{j-J}\|f\|_J.
 \label{v4:newton-source-smoothing}
\end{align}
The construction and the common low bound
$|T_\theta z|_0\le C_T|z|_0$ are given in
Appendix~\ref{v4:tools}. Every accumulated correction has its support
in the fixed compact velocity interval supplied by the inverse, so its zero
extension is smooth. A smoothed background need not retain this support;
it remains an allowed background on the whole control cylinder. The
reference $H_0$ is not smoothed.

Here is the domain test used for the comparison. Suppose $U'$ is at
distance $\Delta$ from the boundary of $U$, and put $R=N(H)$. If
\begin{equation}\label{v4:newton-collar-test}
 \|R\|_{C^0(U)}\le\Delta/4,
 \qquad\|DR\|_{C^0(U)}\le1/2,
\end{equation}
then $y\mapsto x-tR(y)$ is a contraction on
$\overline B(x,\Delta/2)$ for $x\in U'$ and $0\le t\le1$.
Its fixed point is $(\Id+tR)^{-1}(x)$ and lies in $U$. The segment
from that point to $x$ also lies in $U$. Thus the factor of $R$ or $DR$
in the comparison identity is measured only on the preceding residual
domain.

\subsection{The actual update and its four errors}
\label{v4:newton-update-section}

Starting with $z_0=0$, put
\begin{equation}\label{v4:newton-iterates}
 \begin{aligned}
 H_n&=H_0+z_n,& r_n&=N(H_n),&e_n&=\|r_n\|_{0,U_n},\\
 \overline H_n&=H_0+T_{\theta_n}z_n,&
 f_n&=S_{\theta_n}(\chi_nr_n).
 \end{aligned}
\end{equation}
Fix a cutoff $\kappa$ supported in $\Omega_+$ and equal to one near
$\overline\Omega_-$. Define
\begin{equation}\label{v4:newton-update}
 \eta_n=-\mathcal B_{\overline H_n}
       \bigl[\kappa(DP_{\overline H_n})^{-1}f_n\bigr],
 \qquad z_{n+1}=z_n+\eta_n.
\end{equation}
The inverse is applied to the displayed source, with its prescribed
compact control output. Consequently
\begin{equation}\label{v4:newton-linear-cancellation}
 DP_{\overline H_n}\mathcal M_{\overline H_n}\eta_n=-f_n
       \quad\hbox{on }\Omega_-.
\end{equation}
Assume provisionally that the backgrounds used here lie in the protected
ball. The source estimates, smoothing estimates and
\eqref{v4:newton-inverse-bound} imply
\begin{equation}\label{v4:newton-preliminary-increment}
 |\eta_n|_j\le C_j(\varepsilon)(n+1)^{m_j}
                         \theta_n^{j+d}e_n.
\end{equation}
Indeed, every high source factor is bounded by the low residual times
$\theta_n^{j+d}$ and a polynomial collar cost. A high background factor
is bounded by $C_j(1+\theta_n^{j+d}|z_n|_0)$.
The norms of the unsmoothed reference enter the constants.
The smoothed background is thus controlled by the protected low norm of $z_n$.

We next relate the residual of $\overline H_n$ to that of $H_n$.
Write $Z_{n,J}=1+|z_n|_J$. On intermediate collars contained in the
plateau of $\chi_n$, for $J\ge2$, one has
\begin{align}
 \|r_n-f_n\|_0
 &\le C_{J,h,M}(n+1)^{m_J}\theta_n^{-J}Z_{n,J},
 \label{v4:newton-source-tail0}\\
 \|r_n-f_n\|_1
 &\le C_{J,h,M}(n+1)^{m_J}\theta_n^{1-J}Z_{n,J},
 \label{v4:newton-source-tail1}\\
 \|N(\overline H_n)\|_1
 &\le C_{J,h,M}(n+1)^{m_J}
          \bigl(\theta_ne_n+\theta_n^{1-J}Z_{n,J}\bigr).
 \label{v4:newton-smoothed-residual}
\end{align}
For the first two inequalities, apply the smoothing tail estimate to
the zero extension of $\chi_nr_n$ and then use
\eqref{v4:newton-extension} and \eqref{v4:newton-value}.
On the collars in question $r_n=\chi_nr_n$. Also
$\|f_n\|_1\le C(n+1)^m\theta_ne_n$. Integrating the first variation
along the segment from $H_n$ to $\overline H_n$ gives
\[
 \|N(\overline H_n)-r_n\|_1
       \le C|(I-T_{\theta_n})z_n|_0
       \le C_J\theta_n^{-J}|z_n|_J.
\]
Together these prove the last inequality. The comparison error uses
this residual at the smoothed background. Its estimate gains the
one derivative from smoothing and places the unsmoothed tail in
$Z_{n,J}$; the induction therefore needs only the low residual $e_n$.

Taylor's formula and \eqref{v4:newton-linear-cancellation} now give the
following exact identity on $U_{n+1}$:
\begin{equation}\label{v4:newton-four-errors}
 \begin{split}
 r_{n+1}={}&r_n-f_n
       +DP_{\overline H_n}
          (\mathcal M_{\overline H_n}^{\mathrm{tr}}
                   -\mathcal M_{\overline H_n})\eta_n\\
 &+(DN_{H_n}-DN_{\overline H_n})\eta_n\\
 &+\int_0^1(1-t)D^2N_{H_n+t\eta_n}[\eta_n,\eta_n]\dd t.
 \end{split}
\end{equation}
The four terms arise from source smoothing, endpoint comparison,
background smoothing, and the nonlinear Taylor remainder, respectively.
Thus the comparison error remains part of the actual residual equation.

Suppose $e_n\le1$. Put $e=|\varepsilon|$ and suppress polynomial
factors in $n$ temporarily. The low norms of the four terms are bounded,
respectively, by
\begin{align}
 C_{J,h,M}\theta_n^{-J}Z_{n,J},&\qquad
 C_{J,h,M}e^{-2}\bigl(\theta_n^{d+1}e_n^2
               +\theta_n^{d+1-J}Z_{n,J}e_n\bigr),
 \label{v4:newton-errors12}\\
 C_{J,h,M}e^{-2}\theta_n^{d-J}Z_{n,J}e_n,&\qquad
 C_{h,M}e^{-4}\theta_n^{2d}e_n^2.
 \label{v4:newton-errors34}
\end{align}
The second estimate uses \eqref{v4:newton-comparison} and
\eqref{v4:newton-smoothed-residual}. For the third, integrate $D^2N$
between $H_n$ and $\overline H_n$ and use the low variation estimate;
the fourth follows directly from \eqref{v4:newton-taylor}.
The source tail uses no inverse. The second and third terms use one
increment, whose constant is bounded by $C_{J,h,M}e^{-2}$ by
Lemma~\ref{v43:newton-parameter-budget}; the Taylor term uses two.
The segment and collar hypotheses used in these estimates will be
verified in the simultaneous induction below, before the next residual
bound is invoked.

Set $\theta_{n+1}=\theta_n^{3/2}$, $\theta_0\ge2$, and choose the
fixed integer $\varpi=2d+4$. For each fixed level, every collar factor is a
polynomial in $n$ and is bounded by a level-dependent constant times
$\theta_n$. Products of finitely many such factors have the same property.
Thus the preceding estimates imply
\begin{align}
 |\eta_n|_j&\le C_j(\varepsilon)\theta_n^{j+\varpi}e_n,
 \label{v4:newton-increment}\\
 e_{n+1}&\le A_J(\varepsilon)
       \bigl[\theta_n^{-J+\varpi}(1+|z_n|_J)+\theta_n^\varpi e_n^2\bigr].
 \label{v4:newton-recurrence}
\end{align}
The loss $\varpi$ is independent of $J$. The constants may depend on $J$;
their explicit parameter bounds are
\begin{equation}\label{v43:newton-recurrence-powers}
 C_j(\varepsilon)\le C_{j,h,M}|\varepsilon|^{-2},
 \qquad A_J(\varepsilon)\le C_{J,h,M}|\varepsilon|^{-4}.
\end{equation}
Thus the exponents are independent of the truncation order and of the
regularity level; the coefficient constants may depend on both.
\section{Finite initialization and smooth exact return}
\label{v4:newton-convergence}
\label{v4:convergence-chapter}

We now close the simultaneous induction required by the preceding estimates.
The argument uses only finitely many initial smallness conditions. After
they have been met, the same sequence of controls converges in every finite
smooth norm. The ball on which exact return holds has a radius fixed before
the iteration begins.

\begin{theorem}[Exact return near a nonzero initial value]
\label{v4:exact-closure}
Fix a seed $h$ in the class of Definition~\ref{v4:seed-class-B} and a
corresponding nonzero initial value $x_*$. There is a ball
$U\Subset\R^2\setminus\{0\}$ containing $x_*$ with the following property.
Given a finite integer $k\ge2$ and a real number $a>1$, a sufficiently
accurate finite approximation $H_{0,\varepsilon,h}$ from
Theorem~\ref{v4:formal-approximation} can be selected so that, for every
sufficiently small nonzero $\varepsilon$, there is a real smooth
$2\pi$-periodic control $H_{\varepsilon,h}$ on the common control cylinder
such that
\begin{equation}\label{v4:newton-closure-conclusion}
 P_{H_{\varepsilon,h}}=\Id\quad\hbox{on }U,
 \qquad
 \|H_{\varepsilon,h}-H_{0,\varepsilon,h}\|_{C^k(\mathcal W)}
       \le C_{h,k,a}|\varepsilon|^a.
\end{equation}
The flow exists on the common full-period tube, and the control stays in
the background neighborhood of Theorem~\ref{v4:local-inverse}. In particular,
on that tube,
\begin{equation}\label{v4:newton-seed-preservation}
 H_{\varepsilon,h}(q,s)=\varepsilon h(s)q^2+o(|\varepsilon|)
       \quad\hbox{in }C^2.
\end{equation}
The completed control is selected separately for each nonzero parameter.
\end{theorem}

The proof combines the full right inverse of Theorem~\ref{v4:local-inverse},
the comparison identity of Proposition~\ref{v4:comparison-defect}, and the
finite approximation of Theorem~\ref{v4:formal-approximation}, all for the
same seed class.

\subsection{A finite induction with a protected background ball}
\label{v4:newton-base-induction}

Choose
\begin{equation}\label{v4:newton-exponents}
 \mu>2\varpi+8,
 \qquad J_0>6\varpi+3\mu+12,
 \qquad J_0\in\mathbb N.
\end{equation}
For this induction, keep $\varepsilon$ and the smooth reference $H_0$ unchanged.
Let $C_*$ and $A_*$ be the constants at level $J_0$ in
\eqref{v4:newton-increment} and \eqref{v4:newton-recurrence}, and put
\[
 \alpha=\frac{2J_0-10\varpi-5\mu}{6}>0,
 \qquad \beta=\frac\mu2-\varpi>0.
\]
The following three finite inequalities imply the residual induction,
once the domain and background conditions are checked:
\begin{equation}\label{v4:newton-finite-gates}
 e_0\le\theta_0^{-\mu},\qquad
 \theta_0^\alpha\ge2A_*(1+2C_*),\qquad
 \theta_0^\beta\ge2A_*.
\end{equation}
Indeed, suppose $e_i\le\theta_i^{-\mu}$ through step $n$. For every
previously constructed increment,
\[
 |\eta_i|_{J_0}\le C_*\theta_i^{J_0+\varpi-\mu}.
\]
Because $J_0+\varpi-\mu>2$, the sum up to $i=n-1$ is at most a constant
times its final term. More precisely, with $\theta_0\ge2$,
\begin{equation}\label{v4:newton-sharp-high-growth}
 1+|z_n|_{J_0}
       \le (1+2C_*)\theta_n^{2(J_0+\varpi-\mu)/3}.
\end{equation}
The same bound holds when $n=0$, since $z_0=0$.
Substitution into \eqref{v4:newton-recurrence} gives
\[
 e_{n+1}\le\theta_n^{-3\mu/2} A_*
   \left((1+2C_*)\theta_n^{-\alpha}+\theta_n^{-\beta}\right)
       \le\theta_{n+1}^{-\mu}.
\]
Only the finite level $J_0$ and the low indices of the construction have
entered this induction.

We verify the assumptions in the order in which they are used. Let
$C_T$ be the low-norm bound for the control smoother. Increase $\theta_0$
so that the prospective low increment sum satisfies
\begin{equation}\label{v4:newton-protection-budget}
 C_0(\varepsilon)\sum_{i=0}^\infty\theta_i^{\varpi-\mu}
       <\frac{R_\varepsilon}{8(1+C_T)}.
\end{equation}
At step $n$, the preceding increments place $z_n$ and $T_{\theta_n}z_n$
strictly inside the allowed ball. The linear inverse is consequently
available and gives the estimate of $\eta_n$, using the already known
$e_n$. That estimate and the same prospective sum show that
$z_n+t\eta_n$, $0\le t\le1$, remains in the ball. The segment between
$z_n$ and $T_{\theta_n}z_n$ also remains there by convexity. Thus the
Taylor and background-difference estimates are legitimate before the
bound for $e_{n+1}$ is proved. This closes the background part of the
simultaneous induction.

For later use, summation of the preceding increments at a level $J$
chosen independently of $n$ gives the more flexible estimate
\begin{equation}\label{v4:newton-rough-high-growth}
 1+|z_n|_J\le C_{J,h,M}|\varepsilon|^{-2}
       \theta_n^{\max(J+\varpi-\mu,0)/(3/2)+1}.
\end{equation}
If $J+\varpi-\mu>0$, the last summand dominates; if it is negative, the sum
is bounded. The extra exponent $1$ covers the zero case. This estimate
uses only increments with index less than $n$.

Apply \eqref{v4:newton-smoothed-residual} with $J=J_0$ and use
\eqref{v4:newton-rough-high-growth}. The polynomial collar factor can
be bounded by one additional power of $\theta_n$. On an intermediate
collar,
\begin{equation}\label{v4:newton-collar-residual}
 \|N(\overline H_n)\|_1
       \le C_{h,M}|\varepsilon|^{-2}\theta_n^{2-\mu}.
\end{equation}
For its high-order tail, the factor $\theta_n^{1-J_0}$ in
\eqref{v4:newton-smoothed-residual}, the exponent
$2(J_0+\varpi-\mu)/3+1$ in the rough growth bound, and one collar
factor give
\[
 1-J_0+\frac{2(J_0+\varpi-\mu)}3+1+1
       =3-\frac{J_0}3+\frac{2\varpi}3-\frac{2\mu}3<2-\mu.
\]
The last inequality follows from $J_0>2\varpi+\mu+3$, which is
implied by \eqref{v4:newton-exponents}. Increasing $\theta_0$ makes
\eqref{v4:newton-collar-residual} smaller than the required fraction of
every gap $c_{\rm gap}(n+1)^{-2}$ and makes its first derivative smaller than
$1/2$. This follows uniformly in $n$, since the reciprocal gaps grow
polynomially and $\theta_n^{2-\mu}$ decays faster than geometrically.
The contraction test \eqref{v4:newton-collar-test} then ensures that the
comparison at step $n$ samples only the preceding residual domain.
This check uses the previously available $e_n$ and $z_n$.

We have now justified all four errors and the residual induction. It follows
that every iterate exists, all full-period flows stay in the fixed flow
tube, and
\begin{equation}\label{v4:newton-base-decay}
 e_n\le\theta_n^{-\mu},
 \qquad
 \sum_{n=0}^\infty|\eta_n|_0
       \le C_{h,M}|\varepsilon|^{-2}\theta_0^{\varpi-\mu}.
\end{equation}
All initial requirements involve the finite levels fixed above.

\subsection{All smooth orders for the same sequence}
\label{v4:newton-smooth-limit}

The low-norm convergence is not yet sufficient for a smooth metric.
We improve the decay of the residual without changing the sequence.
Suppose that, for some $\nu>2\varpi$,
\begin{equation}\label{v4:newton-bootstrap-input}
 e_n\le C_\nu\theta_n^{-\nu}.
\end{equation}
Summing \eqref{v4:newton-increment} gives
\eqref{v4:newton-rough-high-growth} with $\nu$ in place of $\mu$.
Take a finite $J$ larger than both $\nu-\varpi$ and $2\varpi+4\nu+3$.
The exponent of the high-order tail in \eqref{v4:newton-recurrence}
then satisfies
\[
 -J+\varpi+\frac{J+\varpi-\nu}{3/2}+1
       \le \varpi-2\nu.
\]
Hence
\begin{equation}\label{v4:newton-bootstrap-output}
 e_{n+1}\le C'_{\nu}\theta_n^{-(2\nu-\varpi)}
       =C'_{\nu}\theta_{n+1}^{-(2\nu-\varpi)/(3/2)}.
\end{equation}
All constants at this new finite level may be large. They affect only the
constant in the improved decay estimate, or equivalently its starting index,
not the construction or its initial smallness requirements.

Starting with $\nu_0=\mu$, choose
\[
 \nu_{i+1}=\frac76\nu_i-\frac13\varpi.
\]
For $\nu_i>2\varpi$, this lies strictly between $\nu_i$ and
$(2\nu_i-\varpi)/(3/2)$. Moreover
$\nu_{i+1}-2\varpi=(7/6)(\nu_i-2\varpi)$, so $\nu_i$ tends to infinity.
Thus \eqref{v4:newton-bootstrap-input} holds with arbitrarily large finite
exponents. For each fixed $j$, take $\nu_i>j+\varpi+1$; then
\eqref{v4:newton-increment} proves
\[
 \sum_n|\eta_n|_j<\infty.
\]
It follows that the one sequence constructed above has a limit
\begin{equation}\label{v4:newton-limit}
 H_\infty=H_0+\sum_{n=0}^\infty\eta_n\in C^\infty(\mathcal W).
\end{equation}
It is real and $2\pi$-periodic, and its correction retains the fixed compact
velocity support. The protection budget leaves the limit in the common
flow neighborhood.

By \eqref{v4:newton-balls}, the intersection of the residual domains
contains $U=B(x_*,\rho_\infty)$. Finite-time continuous dependence of the
original ODE on its control and \eqref{v4:newton-base-decay} give
\begin{equation}\label{v4:newton-actual-return}
 P_{H_\infty}(x)=x\qquad(x\in U).
\end{equation}
Thus the original oscillator has exact return on the fixed open ball $U$.

\subsection{Selecting one finite approximation}
\label{v4:newton-parameter-selection}

It remains to choose the finite approximation and the parameter interval
so that all the preceding conditions hold and the correction has the
stated size.

First fix the seed, inverse domains, the finite linear losses and background
indices, and then $b,d,\varpi,\mu,J_0$. For the requested $C^k$ estimate,
choose $b+4\ge k$. Write $e=|\varepsilon|\le1$.
Lemma~\ref{v43:newton-parameter-budget} gives $R_\varepsilon=c_0e$.
By \eqref{v43:newton-recurrence-powers}, the two constant thresholds in
\eqref{v4:newton-finite-gates} satisfy
\begin{equation}\label{v43:newton-gate-powers}
 2A_*(1+2C_*)\le C_{h,M}e^{-6},
 \qquad 2A_*\le C_{h,M}e^{-4}.
\end{equation}
The coefficient in \eqref{v4:newton-collar-residual} is at most
$C_{h,M}e^{-2}$, as is the coefficient of the low increment sum.
Their exponents are independent of $M$.

Choose $A>0$ so that
\begin{equation}\label{v4:newton-parameter-margin}
 \begin{gathered}
 A\alpha>6,\qquad A\beta>4,\qquad A(\mu-2)>2,\\
 A(\mu-\varpi)>2+\max\{a,1\},
 \end{gathered}
\end{equation}
and set $\theta_0=e^{-A}$. The first two inequalities absorb the
thresholds in \eqref{v43:newton-gate-powers}. The third makes the
collar residual smaller than the required fractions of all gaps:
for $\theta_0\ge2$, the ratios
$(n+1)^2\theta_n^{2-\mu}/\theta_0^{2-\mu}$ are uniformly bounded.
The last inequality will place the total correction inside the
physical radius $c_0e$ and give its prescribed $C^k$ size.

Theorem~\ref{v4:formal-approximation} now supplies a \emph{finite}
approximation of order high enough that
\begin{equation}\label{v4:newton-initial-power}
 \|N(H_{0,\varepsilon,h})\|_{0,U_0}
       \le C_{h,M}|\varepsilon|^{M+1},
 \qquad M+1>A\mu.
\end{equation}
Here $M\ge5$ is the odd polynomial degree of the control truncation in
Theorem~\ref{v4:formal-approximation}; its scaled control is
$\varepsilon^{-1}H_M(\varepsilon q,s)$ and its residual order is
$M+1$. The approximation is extended with a fixed velocity cutoff
equal to one on the complete flow tube and its chosen margin. On the
control cylinder used by the inverse, it agrees with its uncut expression.
Formula~\eqref{v43:newton-normalized-reference} shows that its normalized
remainder is $O_{h,M}(e^2)$ at every fixed finite level. Hence the
chosen reference lies in the interior of the normalized background
ball for small $e$.

Finally decrease $|\varepsilon|$. The strict power margins absorb the
finite coefficient constants of this chosen approximation and give all
of \eqref{v4:newton-finite-gates}, \eqref{v4:newton-protection-budget}
and the collar conditions. This completes the choices in their required
order: the fixed losses determine the parameter exponents, these determine
$A$, and $A$ determines one finite $M$. The coefficient constants of that
chosen approximation are then absorbed by decreasing $e$.

The low increment bound then gives
\[
 \|H_\infty-H_{0,\varepsilon,h}\|_{C^k}
       \le C_{h,M}e^{-2}\theta_0^{\varpi-\mu}
       =C_{h,M}e^{A(\mu-\varpi)-2}
       \le C_{h,k,a}e^a,
\]
after shrinking the parameter interval once more. This proves
\eqref{v4:newton-closure-conclusion}. Since $k\ge2$ and $a>1$, the
initial leading term and this estimate imply
\eqref{v4:newton-seed-preservation}, and hence
\[
 \partial_q^2H_{\varepsilon,h}(q,s)
       =2\varepsilon h(s)+o(|\varepsilon|).
\]
This completes the proof of Theorem~\ref{v4:exact-closure}.

\begin{remark}[The scope of the completed control]
The completed control is real and $2\pi$-periodic, and retains the seed
through \eqref{v4:newton-seed-preservation}. Half-periodicity, reflection
symmetry, and prescribed jets at $q=0$ may be lost in the correction.
The origin need not remain an equilibrium; exact return and the ensuing
reconstruction are asserted on the fixed ball $U$ away from it.
\end{remark}
\section{Realization of the families and the scope of the construction}
\label{v4:realization-chapter}

The exact-return theorem now completes the families of indicatrices
constructed in the first part. Its role is to close each indicatrix while
retaining the curvature and angular-profile asymptotics of the prescribed
seed. We assemble the geometric conclusions and describe which members
are distinguished by these invariants.

\subsection{Proof of the geometric existence theorem}
\begin{proof}[Proof of Theorem~\ref{v4:main-theorem}]
Fix $h\in\mathcal B^\infty$ and choose a nonzero initial value $x_h$ in
Definition~\ref{v4:seed-class-B}. Choose a small closed ball $K$ about
$x_h$ disjoint from the origin, together with the larger flow and sampling
collars required in Part~\ref{part:analysis}. The strict characteristic
inequalities hold on a neighborhood of this ball. They determine the
local inverse and its finite low indices before the amplitude is decreased.

Theorem~\ref{v4:formal-approximation} supplies an actual finite-order
control with return residual as small a power of $|\varepsilon|$ as required.
Choose its order according to Theorem~\ref{v4:exact-closure}. That theorem
then gives a real smooth $2\pi$-periodic control $H_{\varepsilon,h}$ and a
fixed inner ball $U_h\Subset\operatorname{int}K$ on which
\[
 P_{H_{\varepsilon,h}}=\Id.
\]
The ball has positive radius and is independent of the iteration index.
All estimates concern the original return map, so this is the closure
identity needed in Theorem~\ref{v4:reconstruction}.

Take $k=2$ and $a=3$ in the closure theorem, so that its correction
is $O_h(|\varepsilon|^3)$ in $C^2$ on the sampling cylinder.
The finite initial approximation has leading term $\varepsilon h(s)q^2$;
all its remaining terms have higher powers of $\varepsilon$. A fixed
velocity cutoff is one on a cylinder containing the whole sampled flow,
so it does not alter this assertion there. Consequently
\[
 \|H_{\varepsilon,h}-\varepsilon h(s)q^2\|_{C^2}
      =O_h(|\varepsilon|^3)
\]
on that cylinder. After decreasing the amplitude threshold, the harmonic
comparison in Lemma~\ref{v4:geom-small-control} applies on $K$ and yields
winding number one. Theorem~\ref{v4:reconstruction} therefore constructs
a single Finsler metric on all of $TU_h\setminus0$, with positive
fundamental tensor, zero Landsberg tensor, and zero flag curvature.

Finally Theorem~\ref{v4:geometric-stability} applies to the same source
and the same base ball. Since $h$ has mean one, it is nonzero on some
open time interval. That theorem gives non-Berwald curvature at every
base point in some direction, nonreversibility of each tangent norm, and
the expansion \eqref{v4:intro-leading-tensor}. None of these conclusions
requires the final source to retain its original half-period symmetry.
This proves all assertions.
\end{proof}

\subsection{Explicit seed neighborhoods}
Every even single-frequency seed $h_n$, $n\ge4$, satisfies the balance
equations and belongs to $\mathcal B^\infty$ by
Propositions~\ref{v4:char-four} and~\ref{v4:char-high-frequency}.
For even $n\ge6$, any prescribed nonzero initial value $x_*$ may be used.
The preceding proof then constructs a control and a metric near that
chosen point. This is a separate local construction for each $x_*$; it
does not give one control closing on the whole punctured plane.

The strict gap away from the characteristic line and its nonzero transverse
derivative persist for seeds close in a fixed finite norm. As in
Corollary~\ref{v4:char-families}, intersecting such a neighborhood with the
balance conditions and $c>0$ gives a relative neighborhood in
$\mathcal A^\infty$ to which the preceding proof applies member by member.

For a more explicit parameter set, fix $J$ in
\eqref{v4:intro-ellipsoid}. Near the point
$a_1=\sqrt{70/13}$ and all other coefficients zero, the derivative of the
ellipsoid equation in $a_1$ is nonzero. It can therefore be solved for
$a_1$ as a positive smooth function of the other $2J-1$ coefficients.
Choose those coefficients in a sufficiently small open ball. Every
resulting seed is balanced, has positive drift, and satisfies the same
strict characteristic inequalities. This supplies parameter sets of
arbitrarily large finite dimension. Increasing $J$ is a new finite
choice; no uniform amplitude threshold over all $J$ is claimed.

\begin{remark}[Seed parameters and geometric distinction]
The dimensions above are dimensions of prescribed seed sets.
Proposition~\ref{v41:seed-separation} distinguishes two fixed seed orbits
when their amplitudes are sufficiently small. Its thresholds may depend
on the chosen pair. It does not assert that a whole continuous seed set
injects into a moduli space at one fixed amplitude, or that the completed
metrics can be selected jointly smoothly in all parameters.
\end{remark}

\subsection{Infinitely many distinct geometries}
The following choice of one metric for every even frequency proves
Corollary~\ref{v41:intro-nonhomothetic}. Fourier orthogonality gives
a single positive lower bound for the distance between any two seed
orbits. The amplitude can thus be chosen separately at each frequency.

\begin{corollary}[Infinitely many nonhomothetic surface germs]
\label{v41:infinite-nonhomothetic}
For every even $n\ge4$, one can choose a sufficiently small nonzero
amplitude $\varepsilon_n$ and a metric $F_n=F_{\varepsilon_n,h_n}$
provided by Theorem~\ref{v4:main-theorem} such that no nonempty base open
subset of $F_n$ is homothetic to a base open subset of $F_m$ when $n\ne m$.
In particular these metrics give countably infinitely many pairwise
nonhomothetic, and hence pairwise nonisometric, germs of regular non-Berwald
Landsberg surfaces. Each germ includes all nonzero tangent directions.
\end{corollary}
\begin{proof}
Use the seeds in \eqref{v4:seed-single-mode}. Their coefficients satisfy
\[
 \lambda_n^2=\frac23\left(n^2-7-\frac{12}{n^2-3}\right).
\]
This expression is increasing for $n\ge4$. For distinct positive integers
$n,m$, the trigonometric functions $\cos(ns)$ and
$\cos(m(\delta s+c))$ are orthogonal with respect to the probability measure
$ds/(2\pi)$, for every $\delta\in\{1,-1\}$ and $c\in S^1$.
Consequently
\[
 \begin{aligned}
 \|h_n-h_m(\delta\,\cdot+c)\|_{C^0}
 &\ge\left(\frac1{2\pi}\int_0^{2\pi}
       |h_n(s)-h_m(\delta s+c)|^2\,ds\right)^{1/2}\\
 &=\sqrt{\frac{\lambda_n^2+\lambda_m^2}{2}}
 \ge\lambda_4.
 \end{aligned}
\]
Thus the orbit distances in \eqref{v41:seed-orbit-distance} have one common
positive lower bound, independent of $n,m$.

For each even $n\ge4$, first choose an admissible nonzero initial value
and apply Theorem~\ref{v4:main-theorem}. Decrease its amplitude, separately
for that $n$, until the error in \eqref{v41:profile-asymptotic} satisfies
$\eta_{h_n}(\varepsilon_n)<\lambda_4/4$. This is possible because that error
tends to zero. If a local homothety between two selected metrics existed,
the proof of Proposition~\ref{v41:seed-separation} would give
\[
 \lambda_4\le d_{O(2)}(h_n,h_m)
 \le\eta_{h_n}(\varepsilon_n)+\eta_{h_m}(\varepsilon_m)
 <\lambda_4/2,
\]
a contradiction. The common separation bound makes these independent
amplitude choices sufficient for every pair at once.
\end{proof}

For a fixed seed there is also an isometry test which uses the fibre
average $A_\varepsilon(x)=\overline{\mathcal I_1}(x)$ itself.
By Proposition~\ref{v41:profile-invariant}, a local homothety of ratio
$\lambda>0$ carries it to $\lambda^{-1}A_\varepsilon$.
Since $A_\varepsilon=2\varepsilon+o(|\varepsilon|)$ uniformly in $x$,
it has the sign of $\varepsilon$ and satisfies
$|\varepsilon|<|A_\varepsilon(x)|<3|\varepsilon|$ for all sufficiently
small nonzero amplitudes. Opposite signs therefore exclude local
homothety. For amplitudes of the same sign with
$3|\varepsilon'|<|\varepsilon|$, the two ranges of $|A|$ are disjoint
and exclude local isometry. In particular, even a fixed seed admits
a sequence of pairwise locally nonisometric realizations, obtained by
taking successive positive amplitudes less than one third of the preceding one.

\subsection{Exact return at nonzero initial values}
The completed metric is local in the base, while its fibre directions
are complete. The role of the open base set is to allow differentiation
of the identity $P_H=\Id$; this gives a periodic variational frame and
smooth matching at the angular seam. Shrinking the base does not permit
an angular sector to replace the full indicatrix. Conversely, the
construction does not require the intermediate oscillator trajectories
to remain in this final base set.

\begin{remark}[The zero initial value of the reference]
The simple-characteristic inverse is constructed near a nonzero initial
value. Its reference symbol vanishes identically on the cotangent fibre
over zero, so the present argument makes no assertion of exact return on
a neighborhood of that point. Prescribing a Euclidean central norm and
an exact central Berwald tensor imposes additional conditions that are
not needed for the local existence theorem proved here.
\end{remark}

There is a restriction on complete extensions. At every base point
Theorem~\ref{v4:geometric-stability} and \eqref{v41:cartan-transverse-derivative}
give a unit direction with $\mathcal I_1\ne0$. In any forward complete
extension with vanishing flag curvature and Landsberg tensor, the unit
geodesic in that direction exists for all $t\ge0$. By
\eqref{v44:geom-affine-geodesic}, $\mathcal I_3$ is unbounded on its lift.
In particular there is no such extension with bounded $\mathcal I_3$.

Reversibility also has a consequence for the branch construction.
Consider any exactly closed oscillator of Theorem~\ref{v4:reconstruction}
in the small neighborhood where \eqref{v4:lin-all-pairs-twist} holds.
Fix a point $x$ where its indicatrix is centrally symmetric about zero.
Then $y\mapsto-y$ preserves the fibre energy and therefore the Hessian
arclength on the indicatrix. It preserves orientation and is a fixed-point-free
involution of a circle of length $2\pi$, so it acts by $s\mapsto s+\pi$.
Consequently $u(x,s+\pi)=-u(x,s)$ and, by angular differentiation,
$w(x,s+\pi)=-w(x,s)$. Duality gives $\beta(x,s+\pi)=-\beta(x,s)$.
For the actual closed flow the periodic comparison agrees with it on the
open set of exact return. At $x$ the opposite conormal branches in
\eqref{v4:lin-column-principal} therefore occur at $s$ and $s+\pi$,
and both give the same quotient $\Gamma$. For nonzero amplitude the
Schur principal symbol in \eqref{v4:lin-schur} vanishes in every
cotangent direction at $x$. The simple characteristic condition thus
requires the absence of such a fibre on the region of inversion.

\begin{remark}[Further geometric questions]\label{v43:family-questions}
These families provide a starting point for studying the local geometry
of non-Berwald Landsberg surfaces. The countable collection in
Corollary~\ref{v41:infinite-nonhomothetic} is obtained by selecting uniformly
separated seed orbits. Its countability describes this selection argument,
not a bound on the number of distinct geometries that the construction
may admit.

For a fixed seed and amplitudes of the same sign, the local homothety
classes of its completions remain to be determined. A related question
is whether completed metrics admit continuous or smooth choices on a
seed neighborhood with the profile distinguishing their geometric
classes. Fibrewise real analyticity of the realizations is also open.
One can further ask whether $\mathcal A^\infty=\mathcal B^\infty$.
The last question concerns the
characteristic condition used in this proof; a balanced seed outside
$\mathcal B^\infty$, if one exists, might still admit exact return by
another argument.
\end{remark}
\appendix
\section{Estimates on fixed function spaces}\label{v4:tools}

This appendix records the elementary estimates used in the local inverse
and the iteration. The domains, extensions, and low indices are fixed
before an arbitrary high regularity level is considered. For general
background on tame estimates and smoothing, see \cite[Part~II and
Part~III, Section~1.4]{Hamilton}.

\subsection{Products and interpolation}

On a fixed smooth bounded domain, or for compactly supported functions
after a fixed extension, the Sobolev product estimate is
\begin{equation}\label{v4:tools-product}
 \|uv\|_{H^r}\le C_r
  \bigl(\|u\|_{H^r}\|v\|_{H^b}
       +\|v\|_{H^r}\|u\|_{H^b}\bigr),
 \qquad r\ge b>\frac{n}{2}.
\end{equation}
For completeness, decompose both factors into inhomogeneous dyadic
frequency blocks. In a low--high product, bound the low frequency sum
in $L^\infty$ by its $H^b$ norm and place the weight $2^{jr}$ on the
high block. Reverse the factors for the high--low products. In the
comparable-frequency products, Bernstein's inequality and $b>n/2$
make the remaining geometric series summable. The square sums give
\eqref{v4:tools-product}. This proof explains the occurrence of only one
high norm in each term.

We also use interpolation between fixed integer $C^r$ norms. On a fixed
interval or rectangle with a fixed extension, for $0\le j\le J$,
\begin{equation}\label{v4:tools-interpolation}
 \|u\|_{C^{b+j}}
       \le C_{b,J}\|u\|_{C^b}^{1-j/J}
                     \|u\|_{C^{b+J}}^{j/J}.
\end{equation}
The harmless addition of the lower norm to the right is allowed for
extensions. One way to see the estimate is to split a convolution
representation at frequency $\theta$, obtaining
$\|u\|_{C^{b+j}}\le C(\theta^j\|u\|_{C^b}
+\theta^{j-J}\|u\|_{C^{b+J}})$, and optimize over $\theta\ge1$.
In a product arising from a derivative partition, the interpolation
exponents add. When their sum is at most one, the weighted arithmetic--
geometric mean inequality bounds the product by a sum with one high factor.
The partition calculation for the actual flow is given explicitly in
Lemma~\ref{v4:newton-tame}.

\subsection{Composition and inverse maps}

If $\phi=\Id+v$ and $\phi^{-1}$ have fixed low-order bounds and all
images lie in a fixed extension domain, the chain rule and
\eqref{v4:tools-interpolation} give a high--low composition estimate,
with a fixed offset when written in Sobolev norms. The exact identity
used for differences is
\begin{equation}\label{v4:tools-composition-difference}
 f(\phi_1(x))-f(\phi_0(x))
 =\int_0^1Df\bigl(\phi_0(x)+t(\phi_1-\phi_0)(x)\bigr)
                         (\phi_1-\phi_0)(x)\dd t.
\end{equation}
The segment must lie in the chosen domain. Thus one input derivative is
spent to retain the small difference factor. In the endpoint comparison
this factor is the actual return error.

For $Q=(\Id+R)^{-1}$, differentiation of $Q+R\circ Q=\Id$ gives
\[
 DQ=(I+DR\circ Q)^{-1}.
\]
Higher derivatives are obtained by differentiating this identity. The
highest derivative of $Q$ appears linearly with the same invertible
coefficient; the other terms are products of lower derivatives and
derivatives of $R$. The product and interpolation estimates therefore
give fixed-offset high--low bounds whenever $\|DR\|_\infty<1/2$.
This argument controls estimates on a fixed domain. The separate collar
test in \eqref{v4:newton-collar-test} guarantees that the required inverse
values and segments actually remain in that domain.

\subsection{One smoothing family at all finite orders}
\label{v4:tools-smoothing}

Choose a real even function $\vartheta\in C_c^\infty(\R^n)$ equal to one
near the origin, and let $K$ be its inverse Fourier transform. Then $K$
is a Schwartz function,
\[
 \int K=1,\qquad\int x^\alpha K(x)\dd x=0\quad(|\alpha|>0).
\]
For a smooth compactly supported input set
$S_\theta f=K_\theta*f$, where
$K_\theta(x)=\theta^nK(\theta x)$, $\theta\ge1$.
Derivatives moved to the kernel prove
\[
 \|S_\theta f\|_{C^{b+J}}
       \le C_{b,j,J}\theta^{J-j}\|f\|_{C^{b+j}}
       \qquad(j\le J).
\]
For the tail, apply Taylor's formula of order $J-j-1$ to
$D^{b+j}f(x-z)$; the vanishing moments annihilate its polynomial part.
Integrating the remainder against $K_\theta$ gives
\[
 \|(I-S_\theta)f\|_{C^{b+j}}
       \le C_{b,j,J}\theta^{j-J}\|f\|_{C^{b+J}}.
\]
If $j=J$, use the $L^1$ norm of $K$ instead. The same multiplier is
used for every order; no order-dependent smoothing operator is selected.

For controls, use the product kernel in $(q,s)$ and periodize it in $s$.
The input is first extended by zero in $q$, which is smooth because every
inverse output is supported inside a fixed compact subinterval. Integrating
the absolute value of the periodized kernel over one period is bounded by
the absolute integral of the original kernel. It follows that the control
smoother has a bound $|T_\theta z|_0\le C_T|z|_0$ with one constant
$C_T$ independent of $\theta$. For targets, it is the already cut-off
source $\chi_nf$ that is extended by zero. Thus the smoothing construction
does not require extension of arbitrary boundary data on a shrinking ball.

\subsection{The two sequence estimates used in the iteration}

Let $\theta_n=\theta_0^{(3/2)^n}$, $\theta_0\ge2$. For each fixed
$m\ge0$, there is $C_m$ with
\[
 (n+1)^m\le C_m\theta_n\qquad(n\ge0).
\]
Indeed, the logarithm of the left side grows like $\log(n+1)$, whereas
the logarithm of the right side grows geometrically. Consequently every
fixed collection of cutoff derivative costs is absorbed in a fixed power
of $\theta_n$ with a level-dependent constant.

For $t>0$, the sequence of positive powers increases sufficiently rapidly
that
\[
 \sum_{i<n}\theta_i^t\le C_t\theta_{n-1}^t
       =C_t\theta_n^{2t/3}.
\]
For $t<0$, the infinite sum is finite and at most a constant times
$\theta_0^t$. These follow by comparing the ratios of consecutive terms;
the ratios move monotonically farther from one. When $t=0$, the sum is
$n$, which is bounded by a constant times $\theta_n$. These elementary
facts account for both the rough high-norm growth and the summability of
the increments. They are statements about the single chosen sequence,
so their use at a higher regularity order does not restart the iteration.
\section{Finite calculus for the complete operators}\label{v4:kernel-appendix}

This appendix records the quantitative facts used in
Sections~\ref{v4:lin-normalization-section}--\ref{v4:lin-inverse-section}.
The operators are the complete operators defined there.  An error
called negative order below includes its finite-frequency part; when
the order statement holds only locally, its exterior pieces are
estimated separately.  The usual amplitude, composition and adjoint
calculus is used in the convention of
\cite[Chapter XVIII]{HormanderIII}. For the fourfold seed, the
fixed-diagonal amplitudes and the conjugated scalar remainder occur
in \cite[Sections~5.2 and~6.3--6.4]{Mao}. Here the bounds must hold
throughout a neighborhood of an arbitrary admissible reference seed,
with the same loss at every Sobolev level.

\subsection{Fixed domains and one high derivative factor}
Fix the outer flow charts, the compact set of control evaluation
points, the branch margins, the finite reduction depths and the
finite collar chain of the main construction.  Controls are smooth
on a slightly larger open cylinder.  All coefficient estimates are
on the chosen compact interior charts with their extensions.  They may
equivalently be obtained after multiplying the control by one fixed
cutoff which is one near the entire evaluation set.  This does not
alter the constructed operators or require an extension of an
arbitrary smooth function across the boundary of its open domain.

Let $b_*$ be a fixed low index.  Suppose a product of jets with
orders $a_i$ satisfies
\begin{equation}\label{v4:kernel-jet-budget}
 \sum_i(a_i-b_*)_+\le n.
\end{equation}
Interpolation between the low norm and the norm of order $b_*+n$
has weights $(a_i-b_*)_+/n$, whose sum is at most one.  Weighted
arithmetic--geometric mean therefore bounds the product by one high
norm times bounded low norms.  With labeled variation directions
$v_1,\ldots,v_j$, each occurring once, this gives
\begin{equation}\label{v4:kernel-marked-product}
 C_n\left[(1+\|a\|_{C^{b_*+n}})\prod_i\|v_i\|_{C^{b_*}}
 +\sum_i\|v_i\|_{C^{b_*+n}}
                         \prod_{k\ne i}\|v_k\|_{C^{b_*}}\right].
\end{equation}
Here $a$ denotes the joint background, with bounded low norm.
The assertion follows by first dividing each nonzero direction by
its low norm and applying the preceding arithmetic.  A zero low norm
means a zero direction.  Fixed compact parameter integrals preserve
the estimate.

For reference, a derivative of $f\circ X$ of order $n$ is a sum of
terms
\[
 (D^kf)\circ X[D^{\nu_1}X,\ldots,D^{\nu_k}X],
 \qquad |\nu_i|\ge1,\quad\sum_i|\nu_i|=n.
\]
At threshold one, the total excess of the map factors is $n-k$;
together with the outer derivative this is $n-1$.  Products split
derivatives without increasing their sum.  Differentiating a reciprocal
has the same property after its undifferentiated denominator is
bounded below.  For an implicit equation $G(X,a)=0$, its highest
derivative has the form
\[
 G_XD^nX=-D_aG[D^na]-\mathcal Q_n,
\]
where every term in $\mathcal Q_n$ partitions the $n$ differentiations
among lower derivatives.  The uniform bound on $G_X^{-1}$ leaves the
same derivative budget.

The flow calculation is equally explicit.  A variation $V$ of a
solution to $\dot X=F_a(X,s)$ satisfies
\[
 \dot V=D_XF_a(X,s)V+D_aF_a(X,s)[\dot a],
 \qquad V(0)=0.
\]
After further differentiation, the highest derivative again solves
this linear equation, with a forcing made of lower jet products.
The finite-time propagator is bounded by the fixed low norm.
Variation of constants, the chain-rule budget and induction give
one high background or one high direction factor.  Inverse flow
maps follow by differentiating their inverse identities.  Time
derivatives are recovered from the differential equation.
These facts apply to the periodization inverse on the same margins
as well as to the true flow.

We apply these rules to the actual amplitude
\eqref{v4:lin-full-amplitude}.  Write $\mathbf x=(x,y)$,
$e=\xi/|\xi|$, and
$\Delta_s(\mathbf x)=\det(n_s(\mathbf x),\partial_s n_s(\mathbf x))$.
The branch times satisfy
\[
 \Xi(\mathbf x,e,s)=\det(e,n_s(\mathbf x))=0,
 \qquad e=\pm\frac{n_s}{|n_s|},
\]
with the sign specifying the branch.  At a root,
\begin{equation}\label{v41:kernel-branch-derivative}
 |\partial_s\Xi|=\frac{\Delta_s}{|n_s|}\ge c,
 \qquad
 \partial_{\mathbf x_j}s_\pm
 =-\frac{\det(e,\partial_{\mathbf x_j}n_s)}
          {\det(e,\partial_sn_s)}\bigg|_{s=s_\pm}.
\end{equation}
The upper bound for $|n_s|$ and the lower bound for $\Delta_s$
are uniform on the low neighborhood.  Angular derivatives have
the same formula with the corresponding derivative of $e$ in the
numerator.  At total endpoint and angular order $n$, the highest
derivative $D^ns_\pm$ is multiplied by this same $\partial_s\Xi$;
all remaining terms partition the $n$ differentiations among lower
branch derivatives and derivatives of $n_s$.  Induction with the
product budget therefore bounds these derivatives by one background
norm of order $n+c_0$.  The branch degree and twist give finitely
many angular charts with uniform margins.  A branch time is a circle
coordinate, and every coefficient evaluated at it is periodic, so
the estimates agree on overlaps.

The branch times are homogeneous of degree zero in $\xi$, and
\[
 \sigma_\pm=\frac{\pm|\xi|}{|n_{s_\pm}|},\qquad
 \frac1{|\sigma_\pm|\Delta_{s_\pm}}
       =\frac1{|\xi|}\frac{|n_{s_\pm}|}{\Delta_{s_\pm}}.
\]
Thus $k_{\pm,ij}$ is $|\xi|^{-1}$ times a smooth function of
the endpoints and direction, multiplied by
$\gamma_\pm(\sigma_\pm)$.  Each frequency differentiation of a
homogeneous factor supplies a factor $|\xi|^{-1}$ and at most one
additional angular derivative.  If a differentiation hits the
frequency cutoff, then $1\le|\sigma_\pm|\le2$; the uniform upper
and lower bounds for $|n_s|$ place $\xi$ in one fixed compact
annulus.  All derivatives there obey the same total derivative
budget.  The amplitude vanishes on a fixed neighborhood of zero.
Consequently its complete order-$-1$ estimates, including the
transition frequencies, use the coefficient index
$|\alpha|+|\beta|+|\gamma|+c_0$ and the frequency weight
$\langle\xi\rangle^{-1-|\gamma|}$.  Derivatives of the fixed
endpoint and velocity cutoffs change the constants only.

Here $n_s$ already contains one spatial flow derivative, while
$\Delta_s$ contains one further time derivative.  The factors $w_s$
are obtained by a matrix reciprocal with a uniform determinant margin.
These are fixed costs, as is the single $\varepsilon$ derivative
used in Hadamard's formula~\eqref{v4:lin-hadamard}.
Let $p_1(x,\xi)$ be the first-order principal symbol factored
in~\eqref{v4:lin-factor}.  Its elliptic factor is obtained near
the simple zero from
\[
 a_0(x,(t,1))=\int_0^1\partial_{\xi_1}p_1
                   (x,(\beta(x)+r(t-\beta(x)),1))\,\dd r.
\]
This formula has a fixed extra derivative cost.  Away from the zero
ordinary division has fixed margins.  Hadamard division in
\eqref{v4:lin-hadamard} likewise costs one fixed parameter derivative.
Thus, for the primitive normal amplitudes and every fixed principal
symbol prescription, there is an integer $c$ such that
\begin{equation}\label{v4:kernel-triangular}
 \sup_{x,y,\xi}\langle\xi\rangle^{-m+|\gamma|}
 |\partial_x^\alpha\partial_y^\beta\partial_\xi^\gamma a|
 \le C_{\alpha,\beta,\gamma}
   (1+\|h\|_{C^{|\alpha|+|\beta|+|\gamma|+c}}
      +\|z\|_{C^{|\alpha|+|\beta|+|\gamma|+c}}).
\end{equation}
For example, in differentiating $|n_{s_\pm}|/\Delta_{s_\pm}$,
write $\mathbf u=(x,y,e)$ and count the $\mathbf u$ derivatives before
restoring the powers of $|\xi|$. A term in $D_{\mathbf u}^n(G(\mathbf
u,s_\pm(\mathbf u)))$ has a derivative of $G$ of order $b+k$ and
$k$ derivatives of $s_\pm$ of degrees $\nu_i\ge1$, with
$b+\sum_i\nu_i=n$. The extra derivative orders to assign to these
factors add to $b+k+\sum_i(\nu_i-1)=n$. In the implicit equation for
$s_\pm$ the term with one factor of degree $n$ is isolated on the left.
Induction gives the same distribution for the remaining terms, with
the reciprocal in~\eqref{v41:kernel-branch-derivative} bounded at low
order. The derivative of $\Delta$ and the evaluation of the
velocity account for fixed initial costs. Taking one low index above
those costs, \eqref{v4:kernel-jet-budget} bounds every product by one
high norm at the index in~\eqref{v4:kernel-triangular}.
The one marked derivative needed for the normalized columns uses
\eqref{v4:kernel-marked-product}.  The offset $c$ absorbs this
derivative and the fixed flow, inverse and implicit-function costs;
it is independent of $\alpha,\beta,\gamma$.

\subsection{Fixed-index mapping and the tame estimate}
For a compactly localized amplitude $a$ of order $m$, define
\[
 |a|^m_{r,k}=\max_{|\alpha|+|\beta|\le r,\,|\gamma|\le k}
 \sup_{x,y,\xi}\langle\xi\rangle^{-m+|\gamma|}
 |\partial_x^\alpha\partial_y^\beta\partial_\xi^\gamma a|.
\]
Expand in the endpoint variables on fixed larger boxes,
$a(x,y,\xi)=\sum_{k,l}e^{ik\cdot x}e^{il\cdot y}a_{kl}(\xi)$.
Integration by parts with $(1-\Delta_x)^L(1-\Delta_y)^L$ gives
\[
 \sup_\xi\langle\xi\rangle^{-m}|a_{kl}(\xi)|
 \le C_L\langle k\rangle^{-2L}\langle l\rangle^{-2L}|a|^m_{4L,0}.
\]
If $M_kf(x)=e^{ik\cdot x}f(x)$, the corresponding operator identity is
\[
 T_a=\sum_{k,l}M_k\,a_{kl}(D)\,M_l.
\]
Peetre's inequality gives
$\|M_k\|_{H^t\to H^t}\le C_t\langle k\rangle^{|t|}$, while
Plancherel's identity gives
$\|a_{kl}(D)\|_{H^{t+m}\to H^t}
 =\sup_\xi\langle\xi\rangle^{-m}|a_{kl}(\xi)|$.
Consequently the operator series is absolutely convergent if
$2L>|t|+2$ and $2L>|t+m|+2$.  In particular,
\begin{equation}\label{v4:kernel-fixed-index}
 \|T_a\|_{H^{t+m}\to H^t}
 \le C_{m,t}|a|^m_{4L,0},\qquad
 L=\lceil|t|+|m|\rceil+3.
\end{equation}
The same bound applies to the adjoint amplitude.  It includes
negative and fractional fixed indices.

For the high-order tame estimate use this bound only at $t=0$.
Integration by parts in the input gives the exact commutator
\[
 [\partial_{x_j},T_a]=T_{(\partial_{x_j}+\partial_{y_j})a}.
\]
There is no boundary term because the complete amplitude includes
its cutoffs.  Applying the fixed bound after $s$ commutations gives
\[
 \|T_af\|_{H^s}\le C_s\sum_{j=0}^s
 (1+\|h\|_{C^{j+c_1}}+\|z\|_{C^{j+c_1}})
                                    \|f\|_{H^{s-j+m}}.
\]
Choose a fixed low source index $b\ge m$ and a fixed Sobolev offset
$q>c_1+1$ on the two-dimensional control cylinder.  Interpolate the
coefficient with weight $j/s$ between its bounded low norm and
$H^{s+q}$.  The source weight between $H^b$ and $H^{s+m}$ is
$(s-j+m-b)_+/(s+m-b)$, at most $(s-j)/s$.  Their sum is at most
one, giving
\begin{equation}\label{v4:kernel-tame-mapping}
 \|T_af\|_{H^s}\le C_s\left[\|f\|_{H^{s+m}}
 +(1+\|h\|_{H^{s+q}}+\|z\|_{H^{s+q}})\|f\|_{H^b}\right].
\end{equation}
The finitely many low values of $s$ use
\eqref{v4:kernel-fixed-index} directly.  This proof explains why
the coefficient index is $s+q$, rather than a multiple of $s$.

For a finite product, the same form of estimate is preserved.
For two factors of orders $m_1,m_2$, the first high source term is
estimated by the second factor at order $s+m_1$.  The first high
coefficient term only requires the second factor at the fixed order
$b_1$, where its coefficient constant is included in the low bound.
Thus the coefficient offset can be chosen to dominate
$q_1$ and $m_1+q_2$, and the low source order to dominate
$b_2$ and $b_1+m_2$.  This gives net source order $m_1+m_2$ and
one high coefficient.  Induction proves the assertion for a fixed
finite word.  The one $\varepsilon$ derivative in the normalized
columns is taken in their primitive formulas before this argument.
The word length and this derivative cost are fixed before $s$ is varied.

\subsection{Finite compositions and their retained errors}
All symbol reductions are finite.  For left symbols $a,b$, the
exact oscillatory composition is
\[
 (a\#b)(x,\xi)=(2\pi)^{-2}\iint e^{-iy\cdot\eta}
                 a(x,\xi+\eta)b(x+y,\xi)\,\dd y\,\dd\eta.
\]
Taylor expansion of the second factor to order $d-1$ gives
\begin{equation}\label{v4:kernel-finite-composition}
 a\#b=\sum_{|\alpha|<d}
       \frac{\partial_\xi^\alpha a\,D_x^\alpha b}{\alpha!}
             +r_d,\qquad r_d\in S^{m_a+m_b-d}.
\end{equation}
The remainder is the actual integral Taylor remainder: its integrand
is a finite sum of
\[
 \frac{d}{\alpha!}(1-t)^{d-1}e^{-iy\cdot\eta}
 (\partial_\xi^\alpha a)(x,\xi+\eta)
 (D_x^\alpha b)(x+ty,\xi),\qquad |\alpha|=d,\quad0\le t\le1,
\]
integrated in $t,y,\eta$.  Endpoint cutoffs are retained in the
symbols or in separate complete smooth kernels.  Integrate by parts
using $(1-\Delta_y)$ and $(1-\Delta_\eta)$ sufficiently many
times to make this integral absolutely convergent after the fixed
symbol derivatives required for a mapping estimate.  The number
depends on the orders, the chosen finite depth and those fixed
frequency seminorms.  Ordinary high endpoint derivatives distribute
by the product rule; the number of additional integrations needed
does not grow with their order.  Consequently each such derivative
requires its own order plus one fixed coefficient offset, by
\eqref{v4:kernel-triangular} and
\eqref{v4:kernel-marked-product}.
More precisely, for a prescribed frequency index $k$ and depth $d$,
the remainder made from primitive factors satisfies
\begin{equation}\label{v43:kernel-remainder-budget}
 |r_d|^{m_a+m_b-d}_{n,k}
 \le C_{n,k,d}\bigl(1+\|h\|_{C^{n+c_{k,d}}}
                         +\|z\|_{C^{n+c_{k,d}}}\bigr).
\end{equation}
Here $c_{k,d}$ also includes the fixed orders of the factors.  It is
independent of $n$.  To see the distinction, first choose the
integrations by parts for the requested order and the $k$ frequency
derivatives.  Each of the subsequent $n$ endpoint derivatives falls
on one factor of that same absolutely convergent integrand; it does
not change its frequency order.  The sum of the additional primitive
derivative counts is $n$.  Applying the product estimate at the fixed
low threshold proves the displayed bound.  No estimate uniform in
all $k$ is inferred for this completed remainder.

Amplitude reduction to a left symbol uses the identical finite
Taylor argument in $y-x$.  Adjoints interchange endpoints and
conjugate the complete amplitude before the same finite reduction.
Thus the difference between an actual word and its first $d$
symbol terms remains an actual negative-order operator with the
stated seminorm estimates.  If its principal terms cancel only on
a platform, the resulting lower order is asserted only there.
In particular, the global difference $P-A_0V_0$ in the main text
is not replaced by an unspecified global order-zero operator.

The finite geometric sums used for the pivot and elliptic factor
are exact polynomials in the already fixed operators.  The
definitions of $S_N=UMV-D$ and $T_N=LU-\Id$ therefore specify
all their low-frequency and smooth pieces.  The orders of their
local diagonal errors are obtained from the finite products in
\eqref{v4:lin-triangular} and~\eqref{v4:lin-scalar-finite}.
They do not require an infinite asymptotic summation or a choice of
parametrix depending on the requested Sobolev order.

\subsection{Separated supports and the complete multiplier tails}
\begin{lemma}[Separated primitive kernels]\label{v4:kernel-separated}
Let a primitive amplitude of order $m$ satisfy
\eqref{v4:kernel-triangular}.  On fixed endpoint sets separated by
$|x-y|\ge d_0>0$, its kernel obeys
\[
 \|K\|_{C_x^sC_y^t}\le C_{s,t,d_0}
 (1+\|h\|_{C^{s+t+c_2}}+\|z\|_{C^{s+t+c_2}}),
\]
with $c_2$ independent of $s,t$.  For a fixed negative input index
$-k$, the resulting operator has the bound
\begin{equation}\label{v4:kernel-separated-map}
 \|Kf\|_{H^s}\le C_{s,k}
 (1+\|h\|_{H^{s+k+c_3}}+\|z\|_{H^{s+k+c_3}})\|f\|_{H^{-k}}.
\end{equation}
\end{lemma}
\begin{proof}
Distribute $s+t$ endpoint derivatives.  If $j$ hit the phase, they
produce a polynomial of degree $j$ in $\xi$; only $s+t-j$ hit
the amplitude.  Integrating by parts in $\xi$ with
$(x-y)\cdot\partial_\xi/(i|x-y|^2)$ a total of
$\max(0,j+\lceil m\rceil+4)$ times makes the frequency integral
absolutely convergent in dimension two.  If some derivatives hit
the polynomial, they decrease its degree.  In every term the total
coefficient derivative count is at most
$s+t+\max(0,\lceil m\rceil+4)+c$.  The separation denominators
are bounded.  This proves the kernel estimate.  Pairing in the input
variable at a fixed negative order, followed by a fixed Sobolev
embedding offset, proves~\eqref{v4:kernel-separated-map}.
\end{proof}

For a finite word, insert the nested identity cutoffs as in the
proof of Lemma~\ref{v4:lin-complete-error}.  At least one primitive
factor crosses two separated supports.  Apply the preceding lemma
to that factor, estimate the right factors at a fixed negative
level, and use the finite-product estimate on its left.  The
number of inserted collars is bounded by the fixed word length.
This proves the same kind of estimate for separated portions of
the actual finite parametrices and their defects.

The complete multipliers require an additional check because their
outer input need not be compact.  On $\R^2$ the kernel of
$(1-\Delta)^{-1/2}$ is
\[
 p_E(x)=\frac{e^{-|x|}}{2\pi|x|}.
\]
On the cylinder of circumference $L$, put
$r_j=(t^2+(y+jL)^2)^{1/2}$.  Its periodized kernel is
\begin{equation}\label{v4:kernel-Bessel}
 p_{\mathcal X}(t,y)=\frac1{2\pi}\sum_{j\in\mathbb Z}
                                  \frac{e^{-r_j}}{r_j}.
\end{equation}
One verification is to use
\[
 (1-\Delta)^{-1/2}=\frac1{\sqrt\pi}
       \int_0^\infty a^{-1/2}e^{-a}e^{a\Delta}\,\dd a.
\]
For the heat kernel $h_a$ on the cylinder, positivity and
$\int_{\mathcal X}h_a=1$ give
\[
 \frac1{\sqrt\pi}\int_0^\infty
          a^{-1/2}e^{-a}\|h_a\|_{L^1}\,\dd a=1.
\]
Thus Fourier transformation can be interchanged with this integral.
At $\lambda=\xi_t^2+(2\pi n/L)^2$ its value is
$\pi^{-1/2}\int_0^\infty a^{-1/2}e^{-(1+\lambda)a}\,\dd a
=(1+\lambda)^{-1/2}$ for every $\lambda\ge0$.
This identifies~\eqref{v4:kernel-Bessel} with the full inverse
multiplier.  Applying $1-\Delta$ in distributions gives the kernel
of $\Lambda_{\mathcal X}$.  For $r>0$, direct radial differentiation
gives
$ (1-\partial_r^2-r^{-1}\partial_r)(e^{-r}/r)
 =-e^{-r}(r^{-2}+r^{-3})$.
Hence the off-diagonal restriction of that distribution is
\[
 -\frac1{2\pi}\sum_j e^{-r_j}(r_j^{-2}+r_j^{-3}).
\]
The diagonal distribution remains part of the operator; the last
formula is used only when the two supports are separated.

\begin{lemma}[Separated noncompact multiplier tails]\label{v4:kernel-multiplier-tails}
If $\chi,\eta$ have compact support and $\eta=1$ near
$\supp\chi$, then $\chi\Lambda_{\mathcal X}^{\pm1}(1-\eta)$
maps $H^{-k}$ to $H^s$ for every fixed $s,k$.  The same assertion
holds in the physical plane.
\end{lemma}
\begin{proof}
On a region separated from the diagonal, every fixed derivative of
the kernels above is bounded by $Ce^{-|t|/\sqrt2}$.  Indeed, for
$y\in[-L/2,L/2]$ use
\[
 \sum_je^{-r_j}\le e^{-|t|/\sqrt2}
                  \sum_je^{-|y+jL|/\sqrt2},
\]
whose second factor is uniformly bounded.  Differentiation adds
fixed inverse powers of $r_j$, bounded by the separation.
Let $f=(1-\Delta)^av$ with $v\in L^2$ and $2a\ge k$.
Move this full source differential operator onto the separated
kernel.  Its differentiated kernel is square integrable in both
variables, because its output is compact and its input tail
decays exponentially.  Cauchy--Schwarz proves the desired bound.
Density identifies it with the original localized multiplier.
\end{proof}

This tail estimate is used before forming a compact-endpoint
representation for the scalar remainder and the left parametrix.
Mere off-diagonal smoothness would not have justified those
estimates for a noncompact input.

\subsection{Coordinate changes and the common background neighborhood}
We now estimate the particular conjugate
$B=\mathsf P B_{\mathcal X}\mathsf P^{-1}$ used in
Lemma~\ref{v4:lin-scalar-lemma}, where
$B_{\mathcal X}=\mathsf E B_0\mathsf J$ transfers the compact-endpoint
kernel of the planar operator $B_0$ by
\eqref{v44:scalar-chart-transfer}. The list of original factors and
their raw orders is fixed in Lemma~\ref{v43:finite-target-words}.
In particular, $P$ is a finite sum of words of raw order at most
one; $A_0,E_0,E_N$ have raw order zero.  We keep the exact difference
$C_0=P-A_0V_0$ throughout this calculation.  A bound for the symbol
of a completed composition will be needed only at a fixed number
of frequency derivatives.

Use centered plane coordinates $(t,y)$ in the chart.
Lift the cylinder flow $X$ to the diffeomorphism $\widehat X$ of the
whole plane determined by $\widehat Y(0,y)=y$ and
$\partial_t\widehat Y=-\widehat\beta(t,\widehat Y)$, where
$\widehat\beta$ is the periodic lift of $\beta_{\mathcal X}$.
It satisfies $\widehat Y(t,y+l)=\widehat Y(t,y)+l$, where $l$ is
the circumference. The displacement in either time direction is at most
$\int_{\mathbb R}\|\beta_{\mathcal X}(t,\cdot)\|_\infty\,\dd t$.
The transverse width chosen in Section~\ref{v4:lin-inverse-section}
therefore places both transported endpoint supports inside a single
lift with a fixed margin. On this lift the kernel of $B$ is the
kernel of $\widehat{\mathsf P}B_0\widehat{\mathsf P}^{-1}$, where
$\widehat{\mathsf P}f=f\circ\widehat X$. Indeed, $\rho=1$ on both
endpoint supports in~\eqref{v44:scalar-chart-transfer}, and the
input Jacobian is precisely the Jacobian of $\widehat X$.
The kernel vanishes outside the transported supports on the cylinder.
Thus all intermediate integrations in $B_0$ can be performed on the
plane before transferring the resulting kernel.

Use $X,Y,\beta$ for the lifted flow, its transverse component and
the lifted coefficient in the following calculations. The equation for its first transverse
derivative is
\[
 \partial_yY(t,y)=
 \exp\left(-\int_0^t\beta_y(a,Y(a,y))\,\dd a\right)>0.
\]
The low bound on $\beta$ and its compact time support give uniform
bounds for this derivative and its reciprocal.  The same bounds
hold globally on the auxiliary cylinder.  Differentiating the flow
equation further gives a linear equation for the highest derivative,
with the same bounded propagator.  If a term in its forcing has
$k$ lower flow jets of positive degrees $d_i$ and $a$ explicit
derivatives of $\beta$, then $a+\sum_i d_i=n$ and its derivative
excess at threshold one is
\[
 (a+k-1)+\sum_{i=1}^k(d_i-1)=n-1.
\]
Variation of constants and induction give the one-high-factor
bound at every order $n$.  Time derivatives follow from the original
equation with one fewer differentiation.  The inverse map has the
same bound by its inverse identity.  Composing these estimates with
the fixed-cost construction of $\beta$ from the characteristic root
uses one fixed additional background index.

\begin{lemma}[Conjugation of the original factors]
\label{v43:kernel-conjugated-primitive}
Let $a$ be one of the primitive amplitudes or prescribed coefficient
symbols in the finite construction, with order $m$.  After pullback
by $X$, its near-diagonal amplitude $a^X$ satisfies, for all endpoint
multi-indices $\alpha,\beta$ and frequency multi-indices $\gamma$,
\begin{equation}\label{v43:kernel-conjugated-triangular}
 \sup_{u,v,\xi}\langle\xi\rangle^{-m+|\gamma|}
 |\partial_u^\alpha\partial_v^\beta
                       \partial_\xi^\gamma a^X|
 \le C_{\alpha,\beta,\gamma}
 \bigl(1+\|h\|_{C^{|\alpha|+|\beta|+|\gamma|+c_X}}
        +\|z\|_{C^{|\alpha|+|\beta|+|\gamma|+c_X}}\bigr),
\end{equation}
where $c_X$ is independent of these multi-indices.  The complete
separated kernels obey the corresponding total endpoint derivative
bound, including the multiplier tails.
\end{lemma}
\begin{proof}
Use fixed lifted coordinate charts near the diagonal and set
\[
 A_X(u,v)=\int_0^1DX(v+r(u-v))\,\dd r,
 \qquad X(u)-X(v)=A_X(u,v)(u-v).
\]
This matrix is triangular.  Its second diagonal entry is an average
of the positive derivative $Y_y$, so it and $A_X^{-1}$ have common
low bounds.  For an original kernel $K$, the exact conjugated kernel
with respect to product Lebesgue densities is
\begin{equation}\label{v43:kernel-exact-pullback}
 K^X(u,v)=K(X(u),X(v))J_X(v),\qquad J_X=\det DX>0.
\end{equation}
The frequency substitution $\xi=A_X(u,v)^{\mathsf T}\zeta$ therefore
gives
\begin{equation}\label{v43:kernel-primitive-pullback}
 a^X(u,v,\xi)=
 a\bigl(X(u),X(v),A_X(u,v)^{-\mathsf T}\xi\bigr)
             \frac{J_X(v)}{|\det A_X(u,v)|}.
\end{equation}
All fixed cutoffs are included in this formula.

Here is the derivative count for this expression.  Put
$n=|\alpha|+|\beta|$ and $k=|\gamma|$.  Suppose $n_1$ of the
endpoint derivatives hit the composed amplitude and the rest hit
the density quotient.  A chain-rule term whose outer amplitude
receives $r$ argument derivatives contains inner jets of positive
degrees $d_1,\ldots,d_r$ with sum $n_1$, before the $k$ frequency
derivatives are distributed.  Including the density derivatives,
its total excess is bounded by
\[
 r+k+\sum_{i=1}^r(d_i-1)+(n-n_1)=n+k.
\]
A positive derivative of $A_X$ uses one extra derivative of $X$;
matrix inversion and the density use only products and reciprocals.
Raising the fixed low threshold absorbs that extra derivative and
the fixed cost of obtaining $X$ from $(h,z)$.  The primitive estimate
\eqref{v4:kernel-triangular}, the flow count above, and the product
estimate~\eqref{v4:kernel-jet-budget} consequently leave one high
background norm at index $n+k+c_X$.

There is also a frequency-order count.  Whenever an endpoint
derivative differentiates $A_X^{-\mathsf T}\xi$, it supplies one
factor linear in $\xi$ and one additional frequency derivative of
$a$.  These have opposite order changes.  The $k$ external frequency
derivatives lower the net order by $k$.  The uniform bounds on
$A_X^{\pm1}$ compare $\langle A_X^{-\mathsf T}\xi\rangle$ with
$\langle\xi\rangle$.  This proves
\eqref{v43:kernel-conjugated-triangular}, including the bounded
frequency region.  A conjugated differential factor can either be
read from its distribution kernel or from the ordinary chain rule;
it has the same assigned order and count.

For a separated primitive kernel use~\eqref{v43:kernel-exact-pullback}
directly.  The inverse Lipschitz bound preserves a positive separation
of the two endpoint sets.  Applying Lemma~\ref{v4:kernel-separated}
and the ordinary composition rule gives a coefficient index equal
to the total endpoint order plus a fixed offset.  For a fixed
multiplier tail, derivatives of the original kernel contain no
background factors.  Its time decay is preserved exactly because
the first component of $X(t,y)$ is $t$.  The Jacobian and its inverse
are uniformly bounded, so change of variables preserves the
integrability used in Lemma~\ref{v4:kernel-multiplier-tails}.
Differentiation then distributes between this fixed kernel, the
coordinate jets and $J_X$.  The same total derivative count applies.
\end{proof}

For a product of planar operators, the lifted global inverse pair gives
the exact identity
\begin{equation}\label{v43:kernel-conjugated-product}
 \widehat{\mathsf P}(A_1\cdots A_r)\widehat{\mathsf P}^{-1}
   =(\widehat{\mathsf P}A_1\widehat{\mathsf P}^{-1})\cdots
     (\widehat{\mathsf P}A_r\widehat{\mathsf P}^{-1}).
\end{equation}
We use this identity before applying any symbol reductions.  The
compact integral in~\eqref{v4:lin-hadamard} commutes with this
conjugation; its one prescribed parameter derivative is already
included in the primitive budget.  Thus every factor to which
Lemma~\ref{v43:kernel-conjugated-primitive} is applied has the
all-frequency estimate just proved.

\begin{proposition}[The scalar grading for the finite construction]
\label{v43:kernel-conjugated-finite-word}
There are fixed finite indices $q_*,b_0$ such that the complete
operators $B_0$ on the plane and $B=\mathsf PB_{\mathcal X}\mathsf P^{-1}$
on the cylinder, and their
adjoints, have representations by finitely many compact-output left
symbols of order zero and retained compact-endpoint smooth kernels.
For these representations, the grading $|B|_j$ of
Lemma~\ref{v4:lin-scalar-lemma} satisfies
\begin{equation}\label{v43:kernel-B-grading}
 |B|_j+|B_0|_j
 \le C_j\bigl(1+\|h\|_{H^{j+q_*}}
                    +\|z\|_{H^{j+q_*}}\bigr),\qquad j\ge0.
\end{equation}
The same estimate holds for the chosen adjoint representations.
The low constants are uniform on one $H^{b_0}$ background
neighborhood, independently of $j$.
\end{proposition}
\begin{proof}
Expand the finite polynomial $E_N$ and the finite words defining
$P$, retaining $P-A_0V_0$ as their difference.  Apply
\eqref{v43:kernel-conjugated-product} to each word.  All its raw
orders are unchanged.  Reduce each resulting word near the diagonal
by the finite Taylor formula~\eqref{v4:kernel-finite-composition}.
For a word of raw order at most one, keep its leading diagonal
product and use the integral first remainder at each multiplication
or amplitude reduction.  Every term other than that leading product
then has order at most zero.  The factor prescriptions here are
classical with integer orders: their homogeneous leading terms
differ from their complete symbols by one lower order.  In the
fixed multiplier and frequency filters this includes their entire
bounded-frequency difference.

The leading diagonal symbol after conjugation is evaluated under
the cotangent substitution
$\kappa_X(u,\xi)=(X(u),DX(u)^{-\mathsf T}\xi)$.
In particular,
\[
 p_1\circ\kappa_X
 -(a_0\circ\kappa_X)
       \bigl((\xi_1-\beta\xi_2)\circ\kappa_X\bigr)=0
\]
on the transported coefficient platform.  The input density and
the frequency Jacobian in~\eqref{v43:kernel-primitive-pullback}
cancel at $u=v$; hence no additional principal factor occurs.
This is precisely the order-one cancellation in the actual
conjugate of $C_0$.  The two endpoint cutoffs of $B_0$ restrict every
nonzero diagonal symbol to that platform.  Multiplication by the
order-zero words of $E_N$ leaves a global order-zero diagonal
operator.  The derivatives of $\beta$ and of the cutoffs occur in
the retained first remainders and are not removed by this
principal calculation.

Fix the exact oscillatory reduction of each localized diagonal
amplitude to a left symbol, retaining its integral remainder.
Choose the finite frequency derivative list needed to estimate
this reduction and to apply the fixed-index mapping estimate.
Apply~\eqref{v43:kernel-remainder-budget} only
for this list.  The number of integrations by parts and each finite
word length are fixed.  If $j$ additional base derivatives are
taken, they distribute among the conjugated primitives with total
count $j$.  The product estimate and
\eqref{v43:kernel-conjugated-triangular} therefore bound every
resulting left symbol $b_\ell$ by
\begin{equation}\label{v43:kernel-left-symbol-grading}
 \max_{|\alpha|\le j}\sup_{u,\xi}
       |\partial_u^\alpha b_\ell(u,\xi)|
 \le C_j\bigl(1+\|h\|_{C^{j+c_*}}
                    +\|z\|_{C^{j+c_*}}\bigr).
\end{equation}
This bound applies to the retained symbol remainders as well as
to their displayed finite Taylor terms.  Their frequency derivative
list was fixed after conjugation; it does not grow with $j$.

For clarity, the smooth kernels in this representation are estimated
as complete kernels.  Insert the fixed nested cutoffs along each
finite word.  A term leaving the local diagonal region contains a
separated primitive factor.  Move an endpoint derivative towards
that factor using
\[
 \partial_iA=A\partial_i+
                       \operatorname{ad}_{\partial_i}A.
\]
The corresponding identity at the input end is obtained by
integration by parts.  After $j$ derivatives, each term assigns a
total of $j$ derivatives to primitive commutators or to the separated
kernel.  The commutators preserve the assigned orders.  If $a$ of
the derivatives on the separated kernel hit its phase, the
integration-by-parts proof of Lemma~\ref{v4:kernel-separated}
uses $a$ additional frequency derivatives; only the other derivatives
have already hit its amplitude.  The sum is again $j$, plus the
fixed cost determined by the surrounding factor orders.  Those
surrounding factors are estimated at fixed Sobolev levels after
their derivatives have been moved in this way.  This proves the
joint endpoint estimate for a separated finite word with one high
background norm.  More explicitly, after distributing the endpoint
derivatives, estimate the resulting operator from $H^{-2}$ to $H^2$.
In dimension two the point masses have uniformly bounded $H^{-2}$
norms and $H^2$ embeds into continuous functions. Evaluation at the
two endpoints therefore bounds the differentiated kernel, at a
fixed additional derivative cost. The argument also covers the complete
multiplier tails by Lemma~\ref{v43:kernel-conjugated-primitive}.

Write
\[
 \|K\|_{\mathrm{Sch}}=
 \max\left\{\sup_u\int|K(u,v)|\,\dd v,
             \sup_v\int|K(u,v)|\,\dd u\right\}.
\]
For every retained kernel $K_\ell$, the preceding calculation gives
\begin{equation}\label{v43:kernel-Schur-grading}
 \max_{|\alpha|+|\beta|\le j}
 \|\partial_u^\alpha\partial_v^\beta K_\ell\|_{\mathrm{Sch}}
 \le C_j\bigl(1+\|h\|_{C^{j+c_*}}
                    +\|z\|_{C^{j+c_*}}\bigr).
\end{equation}
The final endpoint cutoffs of $B_0$ and their transported versions
make these retained kernels compact at both ends.  Intermediate
noncompact multiplier integrations were estimated before this last
localization.  Bounded-frequency pieces are covered by the same
calculation with a compact frequency integral.

The sum of the left sides of
\eqref{v43:kernel-left-symbol-grading} and
\eqref{v43:kernel-Schur-grading} is exactly the symbol-and-kernel
grading used in the scalar lemma.  These are the same functions
$b_\ell(u,\xi)$ and $K_\ell(u,v)$ for all $j$.  The symbol grading
is the ordinary $C^j$ grading with values in $L^\infty_\xi$.
For the kernels, translation in either endpoint preserves the Schur
norm; derivative interpolation in that Banach norm gives the
corresponding inequality.  All functions use the fixed chart
extensions.  Thus the interpolation in the scalar lemma applies
to this representation without making a new symbol choice at a
higher index.  Taking a fixed Sobolev embedding margin on the
two-dimensional control cylinder proves
\eqref{v43:kernel-B-grading}. For the adjoint the localized amplitude is
$\overline{a^X(v,u,\xi)}$ and the remaining kernel is
$\overline{K^X(v,u)}$. Swapping the endpoints preserves their total
derivative count and interchanges the two integrals in the Schur norm.
Reduction of the adjoint amplitude to a left symbol therefore uses
the same fixed list of frequency derivatives. This proves the bound
for the actual adjoint. Choose $b_0$ to include the finite low norms used for these
reductions, the flow margins, and the particular low grading index
in the scalar lemma.
\end{proof}

For example, the exact commutator formula for a left symbol is
$[\partial_i,\Op(b)]=\Op(\partial_i b)$; for a kernel it is
$(\partial_{u_i}+\partial_{v_i})K$.  The fixed-index estimate and
Schur's test now give, with one fixed $q_0$,
\begin{equation}\label{v43:kernel-B-commutator}
 \|\operatorname{ad}_\partial^\alpha B\|_{L^2\to L^2}
 \le C_\alpha |B|_{|\alpha|+q_0}
 \le C_\alpha\bigl(1+\|h\|_{H^{|\alpha|+q_0+q_*}}
                         +\|z\|_{H^{|\alpha|+q_0+q_*}}\bigr).
\end{equation}
Thus the grading interpolated in
Lemma~\ref{v4:lin-scalar-lemma} belongs to the constructed operator
itself.  That lemma gives the high-order estimates for the same
$L^2$ inverse chosen there.  Pullback and inverse pullback have the
corresponding one-high-factor bounds by the flow calculation above
and duality at fixed negative indices.  No domain, cutoff, or finite
reduction depth changes when the high Sobolev index increases.

There is a simple way to include the seed in the background chart.
On the fixed velocity interval write
\[
 H_{\varepsilon,h,z}
 =\varepsilon\bigl(h^\circ(s)q^2+\widetilde z(q,s)\bigr),\qquad
 \widetilde z=(h-h^\circ)q^2+z.
\]
After the chosen velocity cutoff the map $(h-h^\circ,z)\mapsto
\widetilde z$ is linear and bounded at every Sobolev order with
the same order.  Consequently the estimates for the full control
background apply to seed variation as well.  The seed's special
role is in supplying a balanced finite approximation and the
strict reference characteristic margin, rather than in imposing
an additional restriction on the correction space.

Finally, every derivative budget used above is finite before the
formal truncation is chosen.  One may take a conservative bound
by listing all primitive words, summing their absolute orders,
and taking the maximum over this finite list, together with the
fixed-index seminorm requirements.  This bound need not be optimal.
The essential source orders in the final inverse are the explicitly
computed orders of $\mathcal IVCUG$, and the essential low gain
is the separate $L^2\to H^{3/2}$ calculation of
Lemma~\ref{v4:lin-complete-error}.
\par\medskip
\noindent\textbf{Acknowledgments:} This study was funded by the National Key R \& D Program of China (Grant No. 2024YFA1013201) and the Science and Technology Research Program of Chongqing Municipal Education Commission (Grant No. KJZD-K202301110). The authors used ChatGPT 6 Astra to draft and revise text, to develop and check proofs, and to compare the paper's statements about the cited literature with the sources. The authors take full responsibility for the content of this paper.


\begin{thebibliography}{99}

\bibitem{Bao2007}
D. Bao.
On two curvature-driven problems in Riemann--Finsler geometry.
In \emph{Finsler Geometry, Sapporo 2005---In Memory of Makoto Matsumoto},
Advanced Studies in Pure Mathematics 48, Mathematical Society of Japan, Tokyo, 2007, 19--71.

\bibitem{BCS}
D. Bao, S.-S. Chern, Z. Shen.
\emph{An Introduction to Riemann--Finsler Geometry}.
Springer, 2000.

\bibitem{Bryant1995}
R. L. Bryant.
\emph{Finsler surfaces with prescribed curvature conditions}.
Unpublished manuscript, Duke University, July 8, 1995.
\url{https://sites.math.duke.edu/~bryant/FinslerSurface_ppt.pdf}.

\bibitem{Crampin2024}
M. Crampin.
The current state of play in the Landsberg--Berwald problem of Finsler geometry.
\emph{Extracta Mathematicae} 39 (2024), no.~1, 57--95.
\url{https://doi.org/10.17398/2605-5686.39.1.57}.

\bibitem{Dencker}
N. Dencker.
\emph{Local solvability of quasilinear pseudodifferential operators of real principal type}.
\href{https://arxiv.org/abs/2409.19023v2}{arXiv:2409.19023v2}, 2025.

\bibitem{Elgendi2026}
S. G. Elgendi.
\emph{On Landsberg unicorns of general $(\alpha,\beta)$-type}.
\href{https://arxiv.org/abs/2610.37241v1}{arXiv:2610.37241v1}, 27 October 2026.

\bibitem{ElgendiMuzsnay2024}
S. G. Elgendi, Z. Muzsnay.
The geometry of geodesic invariant functions and applications to Landsberg surfaces.
\emph{AIMS Mathematics} 9 (2024), no.~9, 23617--23631.
\url{https://doi.org/10.3934/math.20241148}.

\bibitem{ElgendiMuzsnay2025}
S. G. Elgendi, Z. Muzsnay.
New insights on the unicorn Landsberg problem in Finsler geometry.
\emph{Mediterranean Journal of Mathematics} 22 (2025), no.~4, article 91.
\url{https://doi.org/10.1007/s00009-025-02859-6}.

\bibitem{ElgendiYoussef2021}
S. G. Elgendi, N. L. Youssef.
A note on L. Zhou's result on Finsler surfaces with $K=0$ and $J=0$.
\emph{Differential Geometry and its Applications} 77 (2021), article 101779.
\url{https://doi.org/10.1016/j.difgeo.2021.101779}.

\bibitem{Hamilton}
R. S. Hamilton.
The inverse function theorem of Nash and Moser.
\emph{Bulletin of the American Mathematical Society (N.S.)}
7 (1982), 65--222.

\bibitem{HomanZhou2017}
A. Homan, H. Zhou.
Injectivity and stability for a generic class of generalized Radon transforms.
\emph{Journal of Geometric Analysis} 27 (2017), no.~2, 1515--1529.
\url{https://doi.org/10.1007/s12220-016-9729-4}.

\bibitem{HormanderIII}
L. H\"ormander.
\emph{The Analysis of Linear Partial Differential Operators III: Pseudo-Differential Operators}.
Classics in Mathematics, Springer, Berlin, 2007.
\url{https://doi.org/10.1007/978-3-540-49938-1}.

\bibitem{LiLiuAnalytic}
M. Li, K. Liu.
\emph{Fibrewise Real-Analytic Landsberg Metrics in Dimension at Least Four Are Berwald}.
\href{https://arxiv.org/abs/2611.05067v1}{arXiv:2611.05067v1}, 3 November 2026.

\bibitem{Mao}
J. Mao.
\emph{The construction of a regular non-Berwald Landsberg surface}.
arXiv:2611.26015v1, 19 November 2026.

\bibitem{MaoZhou2026}
J. Mao, L. Zhou.
\emph{The rigidity of the closed three dimensional regular Landsberg metrics}.
\href{https://arxiv.org/abs/2609.11718v1}{arXiv:2609.11718v1}, 10 September 2026.

\bibitem{Matveev2009}
V. S. Matveev.
On ``All regular Landsberg metrics are always Berwald'' by Z. I. Szab\'o.
\emph{Balkan Journal of Geometry and Its Applications} 14 (2009), no.~2, 50--52.
\url{https://www.emis.de/ft/28745}.

\bibitem{Moser1966I}
J. Moser.
A rapidly convergent iteration method and non-linear partial differential equations. I.
\emph{Annali della Scuola Normale Superiore di Pisa, Classe di Scienze} (3)
20 (1966), no.~2, 265--315.
\url{https://www.numdam.org/item/ASNSP_1966_3_20_2_265_0/}.

\bibitem{MuMiaoHuang2025}
X. Mu, H. Miao, L. Huang.
Isochronous centers and flat Finsler metrics (I).
\emph{Canadian Journal of Mathematics} 77 (2025), no.~4, 1294--1314.
\url{https://doi.org/10.4153/S0008414X24000336}.

\bibitem{Nash1956}
J. Nash.
The imbedding problem for Riemannian manifolds.
\emph{Annals of Mathematics} (2) 63 (1956), no.~1, 20--63.
\url{https://doi.org/10.2307/1969989}.

\bibitem{Quinto1980}
E. T. Quinto.
The dependence of the generalized Radon transform on defining measures.
\emph{Transactions of the American Mathematical Society} 257 (1980), no.~2, 331--346.
\url{https://doi.org/10.1090/S0002-9947-1980-0552261-8}.

\bibitem{SSS2010}
S. V. Sabau, K. Shibuya, H. Shimada.
On the existence of generalized unicorns on surfaces.
\emph{Differential Geometry and its Applications} 28 (2010), 406--435.

\bibitem{Shen2009}
Z. Shen.
On a class of Landsberg metrics in Finsler geometry.
\emph{Canadian Journal of Mathematics} 61 (2009), no.~6, 1357--1374.

\bibitem{Szabo2008}
Z. I. Szab\'o.
All regular Landsberg metrics are Berwald.
\emph{Annals of Global Analysis and Geometry} 34 (2008), no.~4, 381--386.
\url{https://doi.org/10.1007/s10455-008-9115-y}.

\bibitem{SzaboCorrection}
Z. I. Szab\'o.
Correction to ``All regular Landsberg metrics are Berwald''.
\emph{Annals of Global Analysis and Geometry} 35 (2009), no.~3, 227--230.
\url{https://doi.org/10.1007/s10455-008-9131-y}.

\bibitem{Zehnder1975}
E. Zehnder.
Generalized implicit function theorems with applications to some small divisor problems, I.
\emph{Communications on Pure and Applied Mathematics} 28 (1975), no.~1, 91--140.
\url{https://doi.org/10.1002/cpa.3160280104}.

\bibitem{Zhou2014}
L. Zhou.
The Finsler surface with $K=0$ and $J=0$.
\emph{Differential Geometry and its Applications} 35 (2014), suppl., 370--380.
\url{https://doi.org/10.1016/j.difgeo.2014.02.003}.

\end{thebibliography}
\end{document}